\documentclass{amsbook}
\usepackage{custom}
\usepackage{pxfonts}
\usepackage[frozencache]{minted}
\usepackage[hidelinks]{hyperref}
\usepackage{parskip}

\newenvironment{centeredpage}
    {%
       \newpage
       \null\vfill
    }%
    {
	   \vfill\null
	   \newpage
    }

\makeatletter
\newcommand*{\@cryear}{}
\newcommand*{\cryear}[1]{\renewcommand*{\@cryear}{#1}}

\newcommand\copyrightpage{\begin{centeredpage}%
\vspace{-\baselineskip}
\begin{center}
Copyright \copyright\ \@cryear\\ 
     \begin{tabular}[t]{c}
        \@author
      \end{tabular}
\end{center}
\end{centeredpage}%
\global\let\cryear\relax
}

\newcommand*{\@dedication}{}
\newcommand*{\dedication}[1]{\renewcommand*{\@dedication}{#1}}
\newcommand\dedicationpage{%
\begin{centeredpage}
\begin{center}
	\@dedication
\end{center}
\end{centeredpage}
}

\newcommand*{\@acknowledgements}{}
\newcommand*{\acknowledgements}[1]{\renewcommand*{\@acknowledgements}{#1}}
\newcommand\acknowledgementspage{%
\newpage
\vspace{-\baselineskip}
\begin{center}
	{\huge ACKNOWLEDGEMENTS\par}
	 \vspace{20\p@}
	\@acknowledgements
\end{center}
\newpage
}

\newcommand*{\@preface}{}
\newcommand*{\preface}[1]{\renewcommand*{\@preface}{#1}}
\newcommand\prefacepage{%
\newpage
\vspace{-\baselineskip}
\begin{center}
	{\huge PREFACE\par}
	 \vspace{20\p@}
	\@preface
\end{center}
\newpage
}
\makeatother

\dedication{I dedicate this dissertation to my wife, Rachel, and my children, Claire and Will. Work hard; I love you.}

\acknowledgements{I wish to acknowledge and thank my advisor, Chris Judge, for his support, guidance, advice, instruction, and infinite patience as I meandered my way through the graduate program. I also wish to thank the members of my committee, Chris Connell, Peter Sternberg, and Dylan Thurston, for their insight and guidance (and willingness to attend a Saturday morning defense). In my time in graduate school, I have met many interesting and helpful people --- too many to name individually with any hope of a complete list. Among names that spring to mind: Erik Wallace, Pawan Patel, and especially Matt Drury, who has shown the way professionally as well as academically. Of course, these acknowledgments would not be complete without mentioning my fellow troublemakers and entrepreneurs David Sprunger, Ranjan Rohatgi, and Tristan Tager. I would also like to acknowledge the professors who put up with me in their classes; in the intervening years, I have learned more from your classes than I thought possible. I wish to thank my coworkers, especially Matt Best, Meg Walters, and Chris Prokop, for their encouragement. Thank you to my parents, both natural and in-law, whose unyielding faith in my eventual completion I often felt I did not deserve. Last, in the position of honor, thanks to my children Claire and Will, and my wife Rachel, for your unconditional love, and for sacrificing so many evenings and weekends on the altar of writing.
}

\title{Laplace Subspectrality}
\author{Neal Coleman}
\date{December 2017}

\cryear{2017}

\begin{document}

\begin{abstract}
Exploring the relationship between geometry and the resonant frequencies of a shape is of interest to pure and applied mathematicians. These resonant frequencies are related to the spectrum of the Laplacian, a partial differential operator. A long-standing research program asks: What geometric information can one deduce from these harmonics?

This dissertation explores a related question. Say one shape is subspectral to another provided each successive resonant frequency of the one is less than the corresponding frequency of the other. What information can be deduced about the relationship between the shapes' geometries?

We use variational arguments to study subspectrality of self-adjoint operators. We develop analytical tools to study Laplace subspectrality. We then study subspectrality in rectangles, construct counterexamples in more general classes of domains, and use the heat trace to relate length subspectrality to Laplace subspectrality. Finally, we discuss a more general question relating eigenvalues to functions, stemming from a 1961 conjecture of Polya.
\end{abstract}

\maketitle

\begin{centeredpage}
\begin{center}
Ah! well a-day! what evil looks

Had I from old and young! 

Instead of the cross, the Albatross 

About my neck was hung.
\end{center}
\end{centeredpage}

\dedicationpage
\acknowledgementspage

\tableofcontents
\newpage
\listoffigures
\chapter*{Introduction}%
\section*{Background and motivation}

The Laplace operator, Laplacian for short, is a second-order partial differential operator defined to act on smooth functions on a compact Riemannian manifold with boundary by mapping a smooth function $f$ divergence of the gradient vector field of $f$. When the Riemannian manifold is a bounded subdomain of $\mb{R}^n$, then the Laplacian is given by
\[ \Delta f = -\sum_i \frac{\partial^2f}{\partial x_i^2}. \]

We consider the eigenvalue problem
\[ \Delta f = \lambda f \]
where either the manifold has no boundary or we have imposed mixed Dirichlet and Neumann boundary conditions on $\partial M$. 

The Riemannian metric induces a measure, which gives the space $L^2$ of square-integrable functions on the manifold. Imposing boundary conditions and restricting the Laplacian to those smooth functions satisfying the boundary conditions allows us to define a nonnegative, symmetric linear operator, densely defined in $L^2$. Using the method of Friedrichs we construct a self-adjoint extension of this operator. If the boundary condition imposed is that each function in the domain take the value of zero on the boundary, we call the resulting self-adjoint extension the Dirichlet Laplacian. If the boundary condition imposed is each function have vanishing normal derivative on the boundary, then we call the resulting self-adjoint extension the Neumann Laplacian. If we have encoded mixed boundary conditions with a function $\nu$, then we call the resulting self-adjoint extension the $\nu$-Laplacian.

If the boundary of the manifold is sufficiently regular, this self-adjoint extension is compactly resolved, yielding a discrete point spectrum of eigenvalues $0\leq\lambda_1\leq\lambda_2\leq\cdots$ of finite multiplicity, accumulating only at infinity, and a decomposition of $L^2$ into orthogonal finite-dimensional eigenspaces. The counting function of this sequence is defined as $N(\lambda) = \#\{k\ |\ \lambda_k \leq \lambda\}$

We motivate the sequel with the following two results tying the sequence of eigenvalues of the Laplacian to the geometry of the manifold. The first is due to Weyl \cite{Weyl1911}. If $M$ is a finite-volume Riemannian manifold, then the function $w^M(\lambda) = \frac{\omega_n}{(2\pi)^n}|M|\lambda^{n/2}$, where $\omega_n$ is the volume of the $n$-dimensional unit disk, is called the Weyl function of $M$. See Definition \ref{weylfunc}.
\begin{thm*}[Weyl's law] If $\Omega$ is a compact domain in $\mb{R}^n$ with piecewise smooth boundary upon which we have imposed Dirichlet boundary conditions, and the Laplacian has eigenvalue counting function $N$, then
\[ \lim_{\lambda\to\infty} \frac{w^\Omega(\lambda)}{N(\lambda)} = 1.\]\end{thm*}

The second is inherent in the work of Weyl and was articulated by Polya \cite{Polya1961}.
\begin{thm*}[Dirichlet domain monotonicity]Suppose $\Omega$ and $M$ are two $n$-dimensional manifolds with piecewise smooth Lipschitz boundary. If $\Omega$ isometrically embeds in $M$, we impose boundary conditions on $\partial M$, and we impose Dirichlet conditions on $\partial\Omega\cap\mbox{interior}(M)$ on $\partial\Omega\cap\partial M$ we match the boundary conditions on $M$, then for each $k$, we have
\[ \lambda_k(M) \leq \lambda_k(\Omega). \]\end{thm*}

Considering the converse yields a natural question: If for all $k$ we have $\lambda_k(M)\leq\lambda_k(\Omega)$, then does $\Omega$ isometrically embed in $M$? More generally, when does this condition on the eigenvalues of $M$ and $\Omega$ hold, and when it does, what conclusions can be drawn?

\section*{Outline of this dissertation}

The subject of this dissertation is pointwise comparison of Laplace spectra. We say one domain and boundary condition is \emph{Laplace subspectral} to another, provided that for all $k$ the $k^{th}$ eigenvalue of the first is no greater than the $k^{th}$ eigenvalue of the second. This is equivalent to the condition that the counting function of the eigenvalues of the first is no less than the counting function of the eigenvalues of the second; see Lemma \ref{subsp_func_lemma}.  In the context of counting functions, we can also compare spectra to functions defined on $[0,\infty)$ and we say that a domain is Laplace subspectral to a function $F:[0,\infty)\to\mb{R}$ provided its counting function is no less than $F$.

In Chapter 1, we precisely define subspectrality for self-adjoint, nonnegative, compactly resolved operators defined in Hilbert spaces:
\begin{def*}[\ref{subspectral_definition}]
Suppose $H$ and $H'$ are two Hilbert spaces. Suppose $T$ and $T'$ are self-adjoint, compactly resolved, nonnegative operators defined in $H$ and $H'$ respectively, with counting functions $N_T$ and $N_{T'}$ respectively. We say that $T$ is subspectral to $T'$ provided for all $x\in\mb{R}$ we have $N_T(x) \geq N_{T'}(x)$. If $T$ is subspectral to $T'$, then we say that $T'$ is superspectral to $T$.
\end{def*}
We also define notions of asymptotic isospectrality, asymptotic and eventual subspectrality, and subspectrality to a function. We then prove sufficient conditions for subspectrality of self-adjoint operators that obey certain conditions:
\begin{thm*}[\ref{subspectr1}]
Suppose $T$ and $T'$ are self-adjoint nonnegative compactly resolved operators defined in a Hilbert space $H$, with associated quadratic forms $\mf{t}$ and $\mf{t}'$ and form domains $V$ and $V'$ respectively.

If $V'\subset V$ and for all $u\in V$ we have the inequality $\mf{t}(u,u)\leq \mf{t}'(u,u)$, then $T$ is subspectral to $T'$.
\end{thm*}
We then make definitions for subspectrality to functions and we prove some propositions relating subspectrality to behavior of the Laplace transform of the counting function of the spectrum, which will be used in Chapter 3.

In Chapter 2, we make a precise definition for subspectrality in the specific context of the Laplace operator densely defined in the space of square-integrable functions of a compact Riemannian manifold with sufficiently regular boundary. We unify all the domain monotonicity theorems as corollaries to a result of Chapter 1 by proving
\begin{thm*}[\ref{partitionthm}]
Let $M$ be a normal manifold. Let $\{\Gamma_i\}_{i=1}^N$ be a finite partition of $M$ by normal, codimension zero manifolds. Impose boundary conditions $\nu$ on $M$. On $\partial M\cap\cup_i\partial\Gamma_i$, impose boundary conditions by restricting $\nu$.

The internal boundaries of the $\Gamma_i$ are $\cup_i \partial\Gamma_i - \partial M$. If we impose Dirichlet conditions on the internal boundaries of the $\Gamma_i$, then $M$ is subspectral to $\sqcup_i\Gamma_i$. If we impose Neumann conditions on the internal boundaries of the $\Gamma_i$, then $M$ is superspectral to $\sqcup_i\Gamma_i$.
\end{thm*}

We establish two tools for studying Laplace subspectrality: 
\begin{thm*}[\ref{quantitative_weyl}]Let $\Omega$ be a normal domain in $\mb{R}^n$. For any $\eta > 0$ denote by $\Omega^{-\eta}$ the set of points in $\Omega$ of distance greater than $\eta$ from $\partial\Omega$ and denote by $\Omega^{\eta}$ the set of points in $\mb{R}^n$ of distance no more than $\eta$ from some point of $\Omega$.

Then $\Omega$ is Dirichlet-superspectral to 
\[ \lambda\mapsto \frac{\omega_n|\Omega^{\eps\sqrt{n}}|}{(2\pi)^n}\bigg(1 + \frac{\pi\sqrt{n}}{\eps\sqrt{\lambda}}\bigg)^n\lambda^{\frac{n}{2}} \]
and Dirichlet-subspectral to 
\[ \lambda\mapsto \frac{\omega_n|\Omega^{-\eps\sqrt{n}}|}{(2\pi)^n}\bigg(1 - \frac{\pi\sqrt{n}}{\eps\sqrt{\lambda}}\bigg)^n\lambda^{\frac{n}{2}}\]
for any $\eps > 0$.
\end{thm*}
and
\begin{prop*}[\ref{metric_subspectrality}]
Let $M$ be a normal manifold with specified boundary conditions $\nu$. Let $g$ and $h$ be two Riemannian metrics. If we have chosen Neumann conditions, let $g$ and $h$ be boundary-conformal (Definition \ref{bdry_conformal}). Denote by $\lambda_k(g)$, resp. $\lambda_k(h)$ the $k^{th}$ eigenvalue of the $\nu$-Laplacian for the metric $g$, resp. $h$. Then
\[ \frac{1}{\delta_+}\bigg(\frac{\delta_-}{\delta_+}\bigg)^{n/2}  \leq \frac{\lambda_k(h)}{\lambda_k(g)} \leq \frac{1}{\delta_-}\bigg(\frac{\delta_+}{\delta_-}\bigg)^{n/2}. \]
\end{prop*}
We then apply these tools to some classes of Euclidean domains. Finally we perform some numerical experiments to study the case of equal-area domains, where the tools proven earlier in the chapter do not apply.

In Chapter 3, we study necessary conditions for subspectrality. In particular, we study consequences for Weyl functions, the heat trace, and the length spectrum on torii and on hyperbolic surfaces. We note that it is not clear that several of the theorems have content; we conjecture that they do.

In Chapter 4, we revisit the result of Polya \cite{Polya1961} and prove a generalization of his theorem:
\begin{thm*}[\ref{PolyaGeneral}]
If a normal domain $\Omega$ has packing constant $\delta$ and Weyl function $w$, then $\Omega$ is Dirichlet-subspectral to $w/\delta$.
\end{thm*}

We also prove a result about the $n$-sphere related to Polya's theorem:
\begin{thm*}[\ref{polya_sphere_thm}]
The $n$-sphere is not superspectral to its Weyl function.
\end{thm*}

We also prove an analogous result for the heat trace and a result regarding the asymptotic behavior of the counting function of certain sequences of Euclidean domains. We numerically study Polya's conjecture in the case of Neumann boundary conditions, and we make a conjecture regarding a Polya-type theorem for isometric domains tiling an unbounded Riemannian manifold.

Finally, in the Appendices, we recap the spectral theorem for compactly resolved self-adjoint operators and the discreteness of the spectrum of the Laplacian for the sake of completeness, and we provide the computer code used to conduct numerical investigation and generate figures.

\section*{Brief history of relevant advances}

The study of the eigenvalues of the Laplacian on a domain or manifold $\Omega$ and their relation to the geometry of $\Omega$ is known as ``spectral geometry'' or ``global harmonic analysis.'' Standard monographs include Berger-Gauduchon-Mazet \cite{Berger1971} and Chavel \cite{Chavel1984}; more recent introductions to the subject include Rosenberg \cite{Rosenberg1997} and Labl\'ee \cite{Lablee2015}. See also the survey by Grebenkov-Nguyen \cite{Grebenkov2013}.

We now give a brief history of highlights in the subject pertaining to the topic of the present work. Working in the first decades of the twentieth century, researchers studied relationships between the spectrum of the Laplacian on a domain and asymptotic formulas which involved geometric quantities. Lord Rayleigh in 1900 \cite{Rayleigh1900} described an asymptotic formula for the vibratory modes of a cubical container. In 1911 \cite{Weyl1911} and 1912 \cite{Weyl1912}, Weyl showed that for a bounded domain in $\mb{R}^n$ the counting function of the spectrum, that is, the function mapping a real number $\lambda$ to the number of eigenvalues no greater than a given $\lambda$, is asymptotically proportional to $\lambda^{n/2}$, where the constant of proportionality is a product of the volume of the domain, and a constant depending only on $n$. Weyl's argument is described in Courant-Hilbert \cite{Courant1953} Chapter 8.

In the mid twentieth century, researchers continued to study and refine comparisons between the Laplace spectrum and functions with quantities of geometric significance. Quantities which can be deduced from the Laplace spectrum of a domain are said to be ``audible.'' For example, the results of Weyl show that the volume of a domain is audible. Polya \cite{Polya1961} showed that, for planar domains given Dirichlet conditions, the counting function of the spectrum exceeds the value of the polynomial to which Weyl showed it is asymptotic.

Researchers also studied whether the isometry class of a given manifold is audible, that is, determined by the Laplace spectrum. Two manifolds or domains are said to be ``isospectral'' provided their Laplace spectra, including multiplicity, are equal.

The first examples of nonisometric, isospectral manifolds were provided by Milnor in 1964 \cite{Milnor1964}. In a seminal talk and paper two years later, Kac \cite{Kac1966} posed the question, ``Can one hear the shape of a drum?'' That is, when restricted to planar domains, are two isospectral domains necessarily isometric?

A series of results over next several decades continued progress in understanding the asymptotic distribution of eigenvalues, establishing additional audible quantities, and sharpening error estimates on comparisons between the eigenvalue sequence and functions to which the sequence is asymptotic. For example, McKean and Singer \cite{Mckean1967} showed the Euler characteristic of a domain is audible; Van den Berg and Srisatkunaraja \cite{VandenBerg1988} showed that in a Euclidean polygon, a certain function of the reciprocals of the angles is audible; and Ivrii \cite{Ivrii1980} provided the sharpest asymptotic error estimates to date, for manifolds with a certain condition on the geodesic flow.

Researchers also continued to construct examples of isospectral, nonisometric manifolds. Sunada \cite{Sunada1985} provided a method of constructing such isospectral manifolds using covering spaces. In 1992 Gordon, Webb, and Wolpert \cite{Gordon1992} adapted Sunada's technique to construct examples of simply-connected planar domains which are isospectral but not isometric. However, researchers also discovered that when one considers more restrictive classes of domain or manifold, isometry classes are audible; for example, this is trivial for Euclidean disks and an exercise for Euclidean rectangles; for Euclidean triangles, this was the subject of Durso's thesis \cite{Durso1989}. The result was reproven in a different fashion by Grieser and Maronna \cite{Grieser2013}. For many classes of domain, such as hyperbolic triangles or convex Euclidean domains, whether isometry class is audible remains unknown.

\mainmatter
\chapter{Subspectral operators in Hilbert spaces}%
\section{Preliminaries}

Suppose $H$ is a Hilbert space with Hermitian inner product $(\cdot,\cdot)$ and associated norm $\|u\|^2 = (u,u)$ for all $u\in H$. Suppose that $T$ is an unbounded symmetric operator with domain $\dom T = D$. Suppose further that $T$ is densely defined, that is, $D$ is dense in $H$, and that $T$ is bounded below, that is, there is some real $K$ so for all $u\in D$ we have $(Tu,u) \geq K$. Define on $D\times D$ the quadratic form $\mf{t}(u,v) = (Tu,v)$. For $\sigma > -K$, the quadratic form $\mf{t}_\sigma(u,v) = \sigma (u,v) + \mf{t}(u,v)$ is a Hermitian inner product. We define the Hilbert space $V$ to be the completion of $D$ with respect to this norm, as in the construction of the Friedrichs extension, Theorem \ref{friedrichs}.

\begin{defn}[Form domain]
We say that $V$ is the form domain of $T$.
\end{defn}

The inclusion $D\hookrightarrow H$ extends to a bounded injection $V\hookrightarrow H$. If the inclusion is compact, then the resolvent operator of $T$ is compact when it is defined. We call such operators compactly resolved. In the sequel, we will consider unbounded, self-adjoint, nonnegative, compactly resolved operators $T$. We will consider form domains constructed using $\sigma=1$. These operators and form domains will be constructed from symmetric operators via the Friedrichs extension.

\section{Spectrum}

By the spectral theorem, in Appendix 2 Theorem \ref{spectral_theorem}, the spectrum of an unbounded compactly resolved operator $T$ is discrete and comprises eigenvalues of finite multiplicity. See also the spectral theorem for unbounded self-adjoint operators, e.g. Rudin \cite{Rudin1987} ch 13.33 or Gilbarg-Trudinger \cite{Gilbarg1998} ch 5 and 8. Because we consider positive operators, the spectrum of $T$ is a subset of $[0,\infty)$. Denote by $E(\lambda)$ the eigenspace of $T$ associated to the eigenvalue $\lambda$, and denote by $\widehat{E}(\lambda)$ the span of the eigenspaces $E(\lambda')$ for all $\lambda' \leq \lambda$.

If we have ordered the eigenvalues $\lambda_1\leq\cdots$ and chosen an orthonormal $T$-eigenbasis $v_1,v_2\cdots$ of $H$, then we shall write 
\[E(\lambda_k) = \on{span}(v_k) \mbox{ and } \hat{E}_k = \on{span}(v_1,\ldots,v_k).\]
We set the convention that $E(\lambda) = 0$ for $\lambda$ not an element of the spectrum of $T$. By the spectral theorem, the set $\hat{E}(\lambda)$ is dense in $H$ as $\lambda\to\infty$. Because $T$ is symmetric, for any two distinct eigenvalues $\lambda\neq\lambda'$ we have $E(\lambda)\perp E(\lambda')$.

\section{Min-max and max-min principles}

Let $T$ be an operator as above. Recall we have defined the quadratic form $\mf{t}$ on its form domain $V$. Define the Rayleigh quotient
\[\mc{R}(u) = \frac{\mf{t}(u,u)}{(u,u)}.\]
This is a real-valued nonnegative functional defined on $V-\{0\}$.

The variational characterization of eigenvalues given by the following theorems is a consequence of linear algebraic consideration of dimension and codimension. We record the proofs here for the sake of completion, as we shall use similar techniques in the proof of Theorem \ref{subspectr1}. The proofs can also be found in Chavel \cite{Chavel1984} Chapter 1 Section 5 and in Polya \cite{Polya1961} Lemmas 1 and 2.

\begin{thm}\label{maxmin}[Max-min theorem]
Let $T$ be a self-adjoint, unbounded, compactly resolved operator with spectrum $\lambda_1\leq\lambda_2\leq\cdots$. Choose an orthonormal $T$-eigenbasis $u_k$ of $H$, $k=1,2,\ldots$. For any $k$, and any $(k-1)$-dimensional subspace $F_{k-1}$ of $H$, the eigenvalue $\lambda_k$ satisfies
\[ \lambda_k \geq \min{\{\mc{R}(u) \ |\ u\in V, u\perp F_{k-1}\}} \]
where orthogonality is with respect to the inner product in $H$. Equality occurs when $F_{k-1} = E(\lambda_{k-1})$.
\end{thm}
\begin{proof}
The vectors $u_1,\ldots,u_k$ form an orthonormal eigenbasis of $\widehat{E}_k$. Let $v_1,\ldots,v_{k-1}$ denote an orthonormal basis of $F_{k-1}$. The linear system $(u,v_j) = 0$ for $u\in \widehat{E}_k$ and $j=1,\ldots,k$ is underdetermined, hence has a nontrivial solution $u_p\in \wh{E}_k\cap F_{k-1}^\perp$.

We show $\mc{R}(u_p) \leq \lambda_k$. Because $u_p\in \wh{E}_k$, we have $u_p = \sum_{i\leq k} c_iu_i$. Then we compute
\begin{align*}
\mc{R}(u_p) &= \mc{R}\bigg(\sum_{i\leq k}c_iu_i \bigg)\\
	&= \frac{\mf{t}(\sum c_iu_i, \sum c_iu_i)}{(\sum c_iu_i, \sum c_iu_i)} \\
	&= \frac{\sum_{i,j} c_ic_j \mf{t}(u_i, u_j)}{\sum c_i^2} \\
	&= \frac{\sum_{i\leq k} c_i^2\lambda_i}{\sum_{i\leq k} c_i^2} \\
	&\leq \lambda_k
\end{align*}
Therefore, we must have $\min_{u\in F_{k-1}^\perp} \leq \lambda_k$. This establishes the inequality. Equality occurs when $u_p$ is a $\lambda_k$ eigenvector.
\end{proof}

Notice that we do not assume $F_{k-1}\subset V$. This is relevant in the proof of Theorem \ref{subspectr1}.

\begin{thm}\label{minmax}[Min-max theorem]
Let $T$ be a self-adjoint, unbounded, compactly resolved operator with spectrum $\lambda_1\leq\lambda_2\leq\cdots$ and form domain $V$. For any $k$, and any $k$-dimensional subspace $F_k$ of $V$, the eigenvalue $\lambda_k$ satisfies
\[ \lambda_k \leq \max_{u\in F_k} \mc{R}(u). \]
Equality occurs when $F_k = \hat{E}_k$.
\end{thm}
\begin{proof}
Let $u_1,\ldots,u_k$ denote an orthonormal eigenbasis of $\widehat{E}_k$. 

Let $v_1,\ldots,v_k$ denote an orthonormal basis of $F_k$. The linear system $(v,u_j) = 0$ for $v\in F_k$ and $j=1,\ldots,k-1$ is underdetermined, hence has a solution $v_p$ which is perpendicular to $\hat{E}_{k-1}$. Thus we may write $v_p = \sum_{i\geq k} c_iu_i$. 

The Rayleigh quotient of $v_p$ satisfies $\mc{R}(v_p)\geq \lambda_k$:
\begin{align*}
\mc{R}(v_p) &= \mc{R}\bigg(\sum_{i\geq k}c_iu_i \bigg)\\
	&= \frac{\mf{t}(\sum c_iu_i, \sum c_iu_i)}{(\sum c_iu_i, \sum c_iu_i)} \\
	&= \frac{\sum_{i,j} c_ic_j \mf{t}(u_i, u_j)}{\sum c_i^2} \\
	&= \frac{\sum_{i\geq k} c_i^2\lambda_i}{\sum_{i\geq k} c_i^2} \\
	&\geq \lambda_k
\end{align*}

This establishes the inequality. Equality occurs when $v_p$ is a $\lambda_k$-eigenvector. 
\end{proof}

\section{Counting function and subspectrality}

For $T$ a self-adjoint, nonnegative, compactly resolved operator, we make the following definition:
\begin{defn}[Spectral counting function]
Define the counting function $N_T$ of $T$ on $[0,\infty)$ by
\[N_T(x) = \dim \widehat{E}(x).\]
\end{defn}
Because eigenspaces are finite-dimensional, the function $N_T$ is defined on $[0,\infty)$ and is monotone increasing. If $\lambda$ is an eigenvalue of $T$, then because the spectrum of $T$ is discrete, for any sufficiently small $\eps > 0$ we have 
\[N_T(\lambda) = N_T(\lambda - \eps) + \on{multiplicity}(\lambda).\]

We now define subspectrality.
\begin{defn}[Subspectrality; superspectrality]\label{subspectral_definition}
Suppose $H$ and $H'$ are two Hilbert spaces. Suppose $T$ and $T'$ are self-adjoint, nonnegative operators with discrete spectra defined in $H$ and $H'$ respectively. We say that $T$ is subspectral to $T'$ provided for all $x\in\mb{R}$ we have $N_T(x) \geq N_{T'}(x)$. If $T$ is subspectral to $T'$, then we say that $T'$ is superspectral to $T$.
\end{defn}

Note the following.
\begin{lemma}
For operators $T$ and $T'$, denote by $\lambda_k(T)$ and $\lambda_k(T')$ the respective $k^{th}$ eigenvalues of $T$ and $T'$, counted with multiplicity. Then $T$ is subspectral to $T'$ if and only if for all $k$ we have $\lambda_k(T)\leq\lambda_k(T')$. 
\end{lemma}
\begin{proof}
Recall $N_T(x)$ is equal to the cardinality of the set $I(x) = \{k\ |\ \lambda_k(T)\leq x\}$, and likewise $N_{T'}(x)$ is equal to the cardinality of the set $I'(x) = \{k\ |\ \lambda_k(T')\leq x\}$. 

Suppose $\lambda_k(T)\leq\lambda_k(T')$ for all $k$. For any $x\in\mb{R}$, if $\lambda_k(T') \leq x$ then we must also have $\lambda_k(T)\leq \lambda_k(T') \leq x$, so $I(x)\subset I'(x)$, hence $N_T(x) \geq N_{T'}(x)$. 

To prove the converse, suppose $N_T \leq N_{T'}$. For any $k\in\mb{N}$ we have 
\[N_{T'}(\lambda_k(T')) = k = N_T(\lambda_k(T)) \geq N_{T'}(\lambda_k(T))\]
and as $N_{T'}$ is monotone increasing, we have $\lambda_k(T') \geq \lambda_k(T)$.
\end{proof}

\section{Sufficient conditions for subspectrality of operators}

We show a condition sufficient for concluding one operator is subspectral to another. This theorem generalizes the ideas contained in existing proofs of domain monotonicity theorems.

\begin{thm}\label{subspectr1}
Suppose $T$ and $T'$ are self-adjoint nonnegative compactly resolved operators defined in a Hilbert space $H$, with associated quadratic forms $\mf{t}$ and $\mf{t}'$ and form domains $V$ and $V'$ respectively.

If $V'\subset V$ and for all $u\in V$ we have the inequality $\mf{t}(u,u)\leq \mf{t}'(u,u)$, then $T$ is subspectral to $T'$.
\end{thm}
\begin{proof}
Denote by $\mc{R}$ and $\mc{R}'$ resp. the Rayleigh quotients of $T$ and $T'$ resp. Notice first that if $\mf{t}(u,u)\leq \mf{t}'(u,u)$ for all $u\in D'$, then $\mc{R}(u)\leq \mc{R}'(u)$ for all $u\in D'$.

Choose an orthonormal $T$-eigenbasis $u_k$ of $H$ and an orthonormal $T'$-eigenbasis $u_k'$ of $H$. Denote by $E_k'$ the span of $u_k'$ and by $\wh{E}_k'$ the span of $u_1',\ldots,u_k'$.

We apply the maxmin principle Theorem \ref{maxmin} to the $(k-1)$-dimensional subspace $\wh{E}_{k-1}$ of $H$ to conclude that
\[\lambda_k(T') \geq \min \{ \mc{R}'(u)\ |\ u\in V'\cap \wh{E}_{k-1}^\perp \}.\]

For any $u\in V'\cap \wh{E}_{k-1}^\perp$, we have in particular that $\mc{R}'(u) \geq \mc{R}(u)$. As $V'\subset V$, we have
\begin{align*}
\min \{ \mc{R}'(u)\ |\ u\in V',\ u\perp \wh{E}_{k-1} \} &\geq \min \{ \mc{R}(u)\ |\ u\in V,\ u\perp \wh{E}_{k-1} \}\\
	&= \lambda_k(T)
\end{align*}

Equality follows from Theorem \ref{maxmin}, where we have that $\lambda_k(T) \geq \min \{ \mc{R}(u)\ |\ u\in V,\ u\perp F_{k-1} \}$ achieves its minimum when $F_{k-1} = \widehat{E}_{k-1}$.

Thus we have
\begin{align*}
\lambda_k(T') &\geq \min_{u\perp \hat{E}_{k-1}} \mc{R}'(u) \\
	&\geq \min_{u\perp \hat{E}_{k-1}} \mc{R}(u) \\
	&= \lambda_k(T)
\end{align*}
as desired.
\end{proof}

We wish to compare operators on different Hilbert spaces. We prove a lemma for when we have an expansion from one Hilbert space to another.

\begin{lemma}[Pullback subspectrality]\label{pullback_subsp}
Suppose $F:H\to H'$ is a continuous map such that $\|u\|_H^2 \leq \|Fu\|_{H'}^2$ for all $u\in H$. Suppose $T'$ is a self-adjoint, compactly resolved, nonnegative operator with form $\mf{t}'$ and form domain $V'$. Let $T$ be the operator constructed from $\mf{t} = F^*\mf{t}'$ in the Friedrichs extension Theorem \ref{friedrichs}. Then $T$ is compactly resolved and $T'$ is subspectral to $T$.
\end{lemma}
\begin{proof}
By construction, $T$ is self-adjoint. We first show that if $T'$ is compactly resolved, then $T$ is compactly resolved. Denote by $V$ the form domain of $T$, and denote by $\mf{t} = F^*\mf{t}'$ the form associated to $T$. Note that $V = \{u\in H\ |\ Fu\in V'\} = F^{-1}(V'\cap\on{ran}(F))$. Because $F$ is continuous on $H$ and preserves $\mf{t}$, it is continuous on $V$. If $u_j$ is a bounded sequence in $V$, then $F(u_j)$ is a bounded sequence in $V'$. Because $V'$ compactly embeds in $H'$, the sequence $F(u_j)$ has a limit point in $H'$. Because $F$ is continuous, the preimage of that limit point is a limit point for $u_j$ in $H$. Thus the inclusion of $V$ into $H$ is compact and so by Theorem \ref{spectral_theorem} $T$ is compactly resolved.

As $T$ is compactly resolved and self-adjoint, its spectrum is discrete. Because $T'$ is bounded below, so is $T$. Denote by $\lambda_k(T)$ the spectrum of $T$ and by $\lambda_k(T')$ the spectrum of $T'$. Denote by $\mc{R}$ the Rayleigh quotient on $V$ and by $\mc{R}'$ the Rayleigh quotient on $V'$. Choose an orthonormal $T$-eigenbasis for $H$ and denote by $\hat{E}_k$ the span of the first $k$ eigenvectors of $T$.

Because $F$ always increases vectors' norms, it must be injective. By the min-max principle, and using that $F$ is injective to conclude that $\dim F(\hat{E}_k) = \dim \hat{E}_k$,
\begin{align*}
\lambda_k(T) &= \max_{u\in \hat{E}_k} \mc{R}(u) \\
	&= \max_{u\in\hat{E}_k} \frac{\mf{t}(u,u)}{\|u\|_H^2} \\
	&\geq \max_{u\in\hat{E}_k} \frac{\mf{t}(u,u)}{\|Fu\|_{H'}^2} \\
	&= \max_{u\in\hat{E}_k} \frac{\mf{t'}(Fu, Fu)}{\|Fu\|_{H'}^2} \\
	&= \max_{u\in F(\hat{E}_k)} \mc{R}'(u) \\
	&\tag{*} \geq \min_{E\subset V} \max_{u\in E} \mc{R}'(u) \\
	&= \lambda_k(T')
\end{align*}
where the minimum in $(*)$ is taken over all linear subspaces $E$ of $V'$ where $\dim E = \dim\hat{E}_k$.
\end{proof}

In fact, if $F$ is an injective isometry, we have the following observation:
\begin{lemma}\label{inj_isom_subsp}
Suppose $F:H\to H'$ is a map such that $\|u\|_H^2 = \|Fu\|_{H'}^2$ for all $u\in H$. Suppose $T'$ is a self-adjoint, compactly resolved, nonnegative operator with form $\mf{t}'$ and form domain $V'$. Let $T$ be the operator constructed from $\mf{t} = F^*\mf{t}'$ in the Friedrichs extension Theorem \ref{friedrichs}. Then $T$ is compactly resolved and $T'$ is subspectral to $T$. If $V'\subset \on{ran}(F)$, then $T'$ is isospectral to $T$.
\end{lemma}
\begin{proof}
That $T$ is compactly resolved and $T'$ is subspectral to $T$ follows from Lemma \ref{pullback_subsp}. Denote the inner product of $H$ by $n$ and the inner product of $H'$ by $n'$. Because $F$ is an isometry, for all $u,v\in H$, we have $n(u,v) = n'(Fu, Fv)$. Suppose $V'\subset \on{ran}(F)$. We show that $u\in H$ is an eigenvector of $T$ with eigenvalue $\lambda$ if and only if $Fu\in H'$ is an eigenvector of $T'$ with eigenvalue $\lambda$. In the following argument we use that if $u\in V$ has $\mf{t}(u,v) = \lambda n(u,v)$ for all $v\in V$, then in fact $u$ is an eigenvector of $T$ with eigenvalue $\lambda$. This follows from the construction of the inverse operator in the Friedrichs extension \ref{friedrichs}.

Suppose $u\in H$ is an eigenvector of $T$ with eigenvalue $\lambda$. Then $u\in V$ and we have for all $v\in V$ that $\mf{t}(u,v) = \lambda n(u,v)$. As $\mf{t}(u,v) = \mf{t}'(Fu, Fv)$ and $n(u,v) = n'(Fu, Fv)$ we have that $\mf{t}'(Fu, Fv) = \lambda n'(Fu, Fv)$. Because $V'\subset\on{ran}(F)$ we have that as $v$ runs over $V$, $Fv$ runs over $V'$, establishing that $Fu$ is an eigenvector of $T'$ with eigenvalue $\lambda$.

Conversely, suppose $Fu$ is an eigenvector of $T'$ with eigenvalue $\lambda$. Then for all $v\in V$ we have $\mf{t}(u, v) = \mf{t}'(Fu, Fv) = \lambda n'(Fu, Fv) = \lambda n(u, v)$. This establishes that $u$ is an eigenvector of $T$ with eigenvalue $\lambda$.
\end{proof}

We now study the situation of two operators and two norms on the same space.

Now let $H$ be vector space. Suppose $n$ and $n'$ are two positive inner products such that the pairs $(H,n)$ and $(H,n')$ are both Hilbert spaces with the same topology. Then there exist constants $c_H,C_H$ such that
\[ 0 < c_Hn \leq n' \leq C_Hn \]

Let $T$ be unbounded, self-adjoint, and nonnegative on $(H,n)$ with form domain $V$. Let $T'$ be unbounded, self-adjoint, and nonnegative on $(H, n')$ with form domain $V'$.

We make the following definition:
\begin{defn}[Comparable operators]\label{def_comparable}
If the images of $V$ and $V'$ in $H$ are equal and there exist constants $c_V,C_V$ such that 
\[ c_V\mf{t}\leq \mf{t}' \leq C_V\mf{t}, \]
then we say that the operators are comparable.
\end{defn}

\begin{lemma}[Comparing Rayleigh quotients]
Suppose we have two nonnegative, self-adjoint, comparable operators $T$ and $T'$ defined on $(H,n)$ and $(H,n')$ respectively. Denote by $\mc{R}$ (resp $\mc{R}'$) the Rayleigh quotients of $\mf{t}$ (resp $\mf{t}'$) with respect to $n$ (resp $n'$). Then the null spaces of $\mf{t}$ and $\mf{t}'$ coincide, and away from their null spaces we have 
\[\frac{c_V}{C_H} \leq \frac{\mc{R}'}{\mc{R}} \leq \frac{C_V}{c_H}.\]
\end{lemma}
\begin{proof}
Recall $\mc{R} = \frac{\mf{t}}{n}$ and $\mc{R'} = \frac{\mf{t}'}{n'}$. Further recall the estimates $c_Hn \leq n' \leq C_Hn$ and, because the operators are comparable, we have $c_V\mf{t}\leq \mf{t}' \leq C_V\mf{t}$. Then $\mf{t}(u) = 0$ if and only if $\mf{t}'(u) = 0$, and for $u$ such that $\mf{t}(u) \neq 0$, we have
\[ \frac{c_V}{C_H} \leq \frac{n}{n'}\frac{\mf{t}'}{\mf{t}} \leq \frac{C_V}{c_H}  \]
as claimed.
\end{proof}

\begin{prop}[Subspectrality of comparable operators]\label{metr_subsp}
Suppose $T$ and $T'$ are non-negative, compactly resolved, self-adjoint, comparable operators defined in Hilbert spaces $(H,n)$ and $(H,n')$, respectively. Denote by $\lambda_k$ (resp $\lambda_k'$) the eigenvalues of $\mf{t}$ (resp $\mf{t}'$) with respect to $n$ (resp $n'$). Then for all $k$ we have that $\lambda_k = 0$ if and only if $\lambda_k' = 0$, and if either is not equal to zero,
\[ \frac{c_V}{c_H}\leq \frac{\lambda_k'}{\lambda_k} \leq \frac{C_V}{C_H}. \]
\end{prop}

\begin{proof}As $T$ and $T'$ are comparable, we may identify their form domains as a single subspace $V$ of $H$. 

By Theorem \ref{minmax}, we have 
\[\lambda_k' = \min_{F_k\subset V} \max_{u\in F_k} \mc{R}'(u) \]
where $F_k$ ranges over $k$-dimensional subspaces of $V$. Likewise,
\[ \lambda_k = \min_{F_k\subset V}\max_{u\in F_k}\mc{R}(u)\]
where again $F_k$ ranges over $k$-dimensional subspaces of $V$.

By the previous lemma, the null spaces of $\mf{t}$ and $\mf{t}'$ coincide. Because $T$ and $T'$ are both compactly resolved, their null spaces are finite dimensional, so we have $\lambda_i = \lambda_i'$ for all $i = 1,\ldots,\dim\ker \mf{t}$.

On $V - \ker \mf{t}$ we have
\[ \frac{c_V}{C_H}\mc{R} \leq \mc{R}' \leq \frac{C_V}{c_H}\mc{R} \]
for any $u\in V$.

For all subspaces in the remainder of the proof, we restrict to those with dimension greater than $\dim\ker\mf{t}$.

Let $F$ be an arbitrarily chosen finite-dimensional linear subspace of $V$. Let $\mc{R}$ obtain its maximum on $F$ at $u_m$. Likewise, let $\mc{R}'$ obtain its maximum at $u_m'$. Then 
\[ \max_F \mc{R} = \mc{R}(u_m) \leq \frac{C_H}{c_V}\mc{R}'(u_m) \leq \frac{C_H}{c_V}\max_F \mc{R}'. \]
Similarly,
\[ \max_F\mc{R}' = \mc{R}'(u_m') \leq \frac{C_V}{c_H}\mc{R}(u_m') \leq \frac{C_V}{c_H}\max_F \mc{R} \]
and so we have
\[ \frac{c_V}{C_H}\max_F \mc{R} \leq \max_F\mc{R}' \leq \frac{C_V}{c_H}\max_F\mc{R} \]
This holds for every finite-dimensional subspace of $V$.

Let $k$ be an arbitrary positive integer greater than $\dim\ker\mf{t}$. We continue to denote the $k^{th}$ eigenspace of $T'$ by $E_k'$ and the $k^{th}$ eigenspace of $T$ by $E_k$, and we let $F_k$ range over $k$-dimensional subspaces of $V$. 

We use Theorem \ref{minmax} to estimate $\lambda_k'$ from below in terms of $\lambda_k$:
\begin{align*}
\lambda_k' &= \max_{E_k'}\mc{R}' \\
	&\geq \frac{c_V}{C_H}\max_{E_k'}\mc{R}\\
	&\geq \frac{c_V}{C_H}\min_{F_k\subset V}\max_{F_k}\mc{R}\\
	&= \frac{c_V}{C_H} \lambda_k.
\end{align*}
To estimate $\lambda_k$ from below in terms of $\lambda_k'$, we make a similar computation. 
\begin{align*}
\lambda_k &= \max_{E_k}\mc{R} \\
	&\geq \frac{c_H}{C_V} \max_{E_k}\mc{R}' \\
	&\geq \frac{c_H}{C_V} \min_{F_k\subset V}\max_{F_k}\mc{R}_h \\
	&= \frac{c_H}{C_V} \lambda_k'.
\end{align*}
Thus we have
\[ \frac{c_V}{C_H}\lambda_k \leq \lambda_k' \leq \frac{C_V}{c_H}\lambda_k \]
as desired.
\end{proof}

We now discuss the consequences of these inequalities for subspectrality, via the minimax characterization of eigenvalues in terms of the Rayleigh quotient.

\begin{cor}[Sufficient conditions for subspectrality] Suppose we have the assumptions present in Theorem \ref{metr_subsp}. Then $T$ is subspectral to $T'$ provided $c_V\geq C_H$ and $T$ is superspectral to $T'$ provided $C_V \leq c_H$.
\end{cor}
\begin{proof}
If $c_V\geq C_H$ then 
\[\lambda_k' \geq \frac{c_V}{C_H}\lambda_k \geq \lambda_k\] for all $k$. If $C_V\leq c_H$ then
\[\lambda_k'\leq \frac{C_V}{c_H}\lambda_k \leq \lambda_k.\]
\end{proof}

We also have the following necessary condition:
\begin{cor}
If $T$ is subspectral to $T'$ then $C_V \geq c_H$. If $T$ is superspectral to $T'$ then $c_V\leq C_H$.
\end{cor}
\begin{proof}
If $T$ is subspectral to $T'$ then
\[ \lambda_k \leq \lambda_k' \leq \frac{C_V}{c_H}\lambda_k \]
and if $T$ is superspectral to $T'$ then
\[ \frac{c_V}{C_H}\lambda_k \leq \lambda_k' \leq \lambda_k\]
The result follows by dividing out $\lambda_k$.
\end{proof}

This pair of corollaries can be phrased in terms of intervals in $\mb{R}$. Away from the nullspace of $\mf{t}$ and $\mf{t}'$, the ratio $\mf{t}'/\mf{t}$ maps $V$ into $I_v = [c_V,C_V]$ and the ratio $n'/n$ maps $H$ into $I_h = [c_H,C_H]$. If $I_v$ falls to the left of $I_h$ then $T$ is superspectral to $T'$. If $I_v$ falls to the right of $I_h$ then $T$ is subspectral to $T'$.

If $T$ is superspectral to $T'$, however, then we must have only that the left endpoint of $I_v$ is less than the right endpoint of $I_h$, and if $T$ is subspectral to $T'$ then we must have only that the right endpoint of $I_v$ is greater than the left endpoint of $I_h$.

\section{Other notions of subspectrality}

\subsection{Subspectrality to functions}
It is useful to compare spectral counting functions to other functions.

\begin{defn}[Sub/superspectral to a function]\label{SubspFunc}
Let $T$ be a self-adjoint, nonnegative, compactly resolved operator on a Hilbert space $H$ with counting function $N_T:[0,\infty)\to\mb{R}$. Let $F:[0,\infty)\to\mb{R}$ be given. If $N_T\leq F$ then we say that $T$ is superspectral to $F$, or equivalently, $F$ is subspectral to $T$. If $N_T + 1\geq F$ then we say that $T$ is subspectral to $F$, or equivalently, $F$ is superspectral to $T$. 
\end{defn}

For example, a self-adjoint operator on an $n$-dimensional vector space is superspectral to the constant function $F:x\mapsto n$.

We note that if $N_T \geq F$, then $N_T + 1 \geq N_T \geq F$ implies that $T$ is subspectral to $F$.

The following lemma relates subspectrality to a function with comparison of inequalities and motivates the use of the shifted counting function.

\begin{lemma}\label{subsp_func_lemma}[Comparing eigenvalues to functions]
Suppose $F$ is a continuous, monotone increasing function defined on $[0,\infty)$. Suppose $T$ is a self-adjoint, nonnegative, compactly resolved operator with counting function $N_T$. Denote the $k^{th}$ eigenvalue of $T$ by $\lambda_k$.

Then:
\begin{itemize}
\item  $T$ is subspectral to $F$ if and only if $\lambda_k \leq F^{-1}(k)$ for each $k\in\mb{N}$. 
\item $T$ is superspectral to $F$ if and only if $\lambda_k\geq F^{-1}(k)$ for each $k\in\mb{N}$.
\end{itemize}
\end{lemma}
\begin{proof}
As $F$ is continuous and monotone increasing, it has a continuous, monotone increasing inverse $F^{-1}$.

Fix an arbitrary positive integer $k$. If $T$ is subspectral to $F$, then for each $x$, we have $F(x) \leq N_T(x) + 1$. Let $\epsilon$ be less than the difference between $\lambda_k$ and the greatest eigenvalue less than $\lambda_k$. Then
\[F(\lambda_k - \epsilon) \leq N_T(\lambda_k - \epsilon) + 1 = N_T(\lambda_k) - \on{mult}(\lambda_k) + 1 \leq k.\]
Letting $\epsilon\to 0$ we have $F(\lambda_k) \leq k$ and so we have $\lambda_k \leq F^{-1}(k)$. 

If $T$ is superspectral to $F$, then for each $x$, we have $F(x) \geq N_T(x)$. In particular, $F(\lambda_k) \geq N_T(\lambda_k) = k$ and so $\lambda_k \geq F^{-1}(k)$.

As $k$ is arbitrary, these hold for all $k\in\mb{N}$ and one direction of the proposition is proven.

Conversely, suppose for each $k$ we have $F^{-1}(k) \leq \lambda_k$. Then we have $k \leq F(\lambda_k)$.  Let $x$ be an arbitrary positive real number in the complement of the spectrum of $T$. Let $\lambda_j$ and $\lambda_{j+1}$ be the largest eigenvalue of $T$ smaller than $x$ and the smallest eigenvalue of $T$ greater than $x$, respectively. They exist because the spectrum of $T$ is discrete with finite multiplicity. As $F$ is monotone increasing, we have $F(\lambda_j) \leq F(x) \leq F(\lambda_{j+1})$. Because on the interval $[\lambda_j, \lambda_{j+1})$ the function $N_T$ is constant, we have $N_T(x) = N_T(\lambda_j) \leq F(\lambda_j) \leq F(x)$. Thus $T$ is superspectral to $F$.

Similarly, suppose for each $k$ we have $F^{-1}(k) \geq \lambda_k$. Then we have $k\geq F(\lambda_k)$. Let $x$ be an arbitrary positive real number in the complement of the spectrum of $T$ and let $\lambda_j$ and $\lambda_{j+1}$ be the largest eigenvalue of $T$ smaller than $x$ and the smallest eigenvalue of $T$ greater than $x$, respectively. Because on the interval $[\lambda_j, \lambda_{j+1})$ the function $N_T$ is constant, and $j = \max\{k\ |\ \lambda_k = \lambda_j\}$, we have $N_T(x) = N_T(\lambda_j) = j$. As $F$ is monotone increasing, we have that
\[F(x) \leq F(\lambda_{j+1}) \leq j+1 = N_T(\lambda_j) + 1 = N_T(x) + 1\]
as desired.
\end{proof}

\subsection{Asymptotic and eventual subspectrality}

\begin{defn}[Eventual subspectrality]
Let $T$ be a self-adjoint nonnegative compactly resolved operator on a Hilbert space $H$ with counting function $N_T$. For any function $F:[0,\infty)\to\mb{R}$ we say that $T$ is sub(resp super)spectral to $F$ beyond $x_0$ provided for all $x\geq x_0$ we have $N_T(x)\geq (\mbox{resp }N_T^s(x) \leq) F(x)$. If there exists some $x_0$ such that $T$ is sub(resp super)spectral to $F$ beyond $x_0$, then we say that $T$ is eventually sub(resp super)spectral to $F$.
\end{defn}

One way to show that one positive operator is subspectral to another is to show first that the operator is subspectral to the other beyond some $x_0$, and then to check that subspectrality holds for all $x\in[0,x_0]$.

\begin{defn}[Asymptotic subspectrality]
Let $T$ be a self-adjoint nonnegative compactly resolved operator on a Hilbert space $H$ with counting function $N_T$. For any function $F:[0,\infty)\to\mb{R}$ we say that $T$ is asymptotically sub(resp super)spectral to $F$ provided $\lim_{x\to\infty} N_T(x)/F(x) \geq (\mbox{resp }\leq) 1$. 
\end{defn}

\subsection{Isospectrality}

It is useful to relate the notion of subspectrality to that of isospectrality.
\begin{defn}[Isospectrality]
Let $T,T'$ be self-adjoint nonnegative compactly resolved operators on a Hilbert spaces $H,H'$ resp with counting functions $N_T,N_{T'}$ resp. If $N_T = N_{T'}$ then we say that: $T$ is isospectral to $T'$, equivalently, $T$ and $T'$ are isospectral, equivalently, $T$ and $T'$ form an isospectral pair.
\end{defn}

If $T, T'$ are as above and $T$ is both subspectral and superspectral to $T'$ then $T$ is isospectral to $T'$.

We also extend the definitions of eventual subspectrality and asymptotic subspectrality to eventual isospectrality and asymptotic isospectrality in a natural way. We say that a function $F$ is $o(x^\alpha)$ provided $F(x)/x^\alpha\to 0$ as $x\to\infty$. We say $F = G + o(x^\alpha)$ provided $F - G$ is $o(x^\alpha)$.
\begin{defn}[eventual, asymptotic isospectrality]
Let $T$ be as above with counting function $N_T$ and let $F$ be a real valued function on $[0,\infty)$. If for all $x\geq x_0$ we have $N_T(x) = F(x)$ then we say $T$ is isospectral to $F$ beyond $x_0$. If there exists $x_0$ so that $T$ is isospectral to $F$ beyond $x_0$, then we say that $T$ is eventually isospectral to $F$.

If $N_T(x)/F(x) \to 1$ as $x\to\infty$, we say $T$ is asymptotically isospectral to $F$.
\end{defn}

These definitions simplify statements about eigenvalues. Take as an example this result due to Weyl \cite{Weyl1911}, \cite{Weyl1912}, see proof in Courant-Hilbert \cite{Courant1953} Vol I Ch VI:
\begin{thm}[Weyl's law]\label{WeylLaw}
Let $\Omega\subset\mb{R}^n$ be a compact domain with Lipschitz boundary. Let $w$ be the single-term Weyl function of $\Omega$ (see Definition \ref{weylfunc}).
Then the Dirichlet and Neumann Laplace operators are asymptotically isospectral to $w$.
\end{thm}

As another example, consider this conjecture, due to Polya \cite{Polya1961}:
\begin{conj}[Polya's conjecture]\label{PolyaConj}
Let $\Omega\subset\mb{R}^n$ be compact with Lipschitz boundary. Let $w$ be the single-term Weyl function of $\Omega$ (see Definition \ref{weylfunc}). 
Then the Dirichlet Laplace operator is superspectral to $w$ and the Neumann Laplace operator is subspectral to $w$.
\end{conj}

For discussion of boundary conditions and Weyl's law, see Chapter 2. For discussion of Polya's conjecture, see Chapter 4.

\section{Necessary conditions for subspectrality of operators}

Now equipped with the definition of asymptotic isospectrality, we prove the following lemmas regarding necessary conditions for subspectrality. These results will be of use in Chapter 3.

\begin{lemma}[Necessary condition if asymptotic isospectrality is known]\label{subspectr_sum}
Suppose $T$ is an operator in $H$ and $T'$ is an operator in $H'$, both nonnegative self-adjoint with discrete spectrum. Let $0 \leq r_0 < r_1 < \cdots < r_n$ and $0 \leq s_0 < s_1 < \cdots < s_m$ be increasing sequences of real numbers.

Suppose $N_T = c_nx^{r_n} + \cdots + c_1x^{r_1} + o(x^{r_0})$ and $N_{T'} = c_m'x^{s_m} + \cdots + c_1'x^{s_1} + o(x^{s_0})$.

If $T$ is subspectral to $T'$, then $s_m\leq r_n$. In the case that $s_m = r_n$, there is some $k$ so $s_{m-j} = r_{n-j}$ and $c_{m-j}' = c_{n-j}$ for each $j=0,\ldots,k$, and either $c_{n-k} - c_{m-k}' > 0$; or, if $r_{n-k-1} > s_{m-k-1}$, we have $c_{n-k-1} > 0$; or, if $r_{n-k-1} < s_{m-k-1}$, we have $c_{m-k-1}' < 0$.
\end{lemma}
\begin{proof}
By the hypothesis that $T$ is subspectral to $T'$, we have that $N_T - N_{T'} \geq 0$. In particular, the function
\[ \sum_{j=1}^n c_nx^{r_n} - \sum_{j=1}^m c_m'x^{s_m}\]
must have a positive leading term.

Because $N_T$ and $N_{T'}$ are both nonnegative functions, we have that $c_m' > 0$ and $c_n > 0$. So $s_m \leq r_n$. If $s_m = r_n$, then for all $j=0,\ldots,k$ such that $s_{m-j} = r_{n-j}$ and $c_{m-j}' = c_{n-j}$, we have each term vanishing in the sum above, hence the leading term is $(c_{n-k} - c_{m-k}')x^{r_{n-k} - s_{m-k}}$ and so $c_{n-k} - c_{m-k}' \geq 0$. If $c_{n-k} = c_{m-k}'$ then the next term has power either $r_{n-k-1}$ or $s_{m-k-1}$. If the leading power is $r_{n-k-1}$ then $c_{n-k-1} > 0$. If the leading power is $s_{m-k-1}$ then $c_{m-k-1}' < 0$ to ensure the leading term is positive.
\end{proof}

Information about the counting function is often derived from the following integral transform. If $F$ is a function of bounded variation on $[0, \infty)$, we consider the integral $\hat{dF}$ defined as the Riemann-Stieltjes integral
\[ \hat{dF}(s) = \int_0^\infty e^{-st}dF(t). \]
If $F$ is differentiable and $\hat{dF}$ exists, then $\hat{dF}$ is the Laplace transform of the derivative of $F$. The Riemann-Stieltjes integral is defined in Rudin \cite{Rudin1976} Chapter 6.
\begin{lemma}[Necessary condition on Laplace transform]\label{subspectr_lapl}
Suppose $T$ is an operator in $H$ and $T'$ is an operator in $H'$, both nonnegative self-adjoint with discrete spectrum. Denote by $Z$ and $Z'$ the functions $\hat{dN_T}$ and $\hat{dN_{T'}}$, respectively. Suppose $T$ is subspectral to $T'$. Then $Z \geq Z'$. 

Suppose further that $Z(s) = p(s) + e(s)$ and $Z'(s) = p'(s) + e'(s)$ where $p$ is a polynomial of order $m$ with lowest order $r$, the error term $e$ has $e(s)s^r\to 0$ as $s\to0$, the function $e'$ is a polynomial of order $m'$ with lowest order $r'$, and the error term $e'$ has $e'(s)s^{r'}\to 0$ as $s\to 0$.

Then $m \leq n$. In the case that $m=n$, if $c_j = c'_j$ for all $j = m, m-1, \ldots, k+1 > \max{s,s'}$, then $c_k \geq c'_k$.
\end{lemma}
\begin{proof}
Notice:
\[ Z(s) = \int_0^\infty e^{-st}dN_T(t) = \sum_0^\infty e^{-\lambda_k(T)s} \]
where eigenvalues are counted with multiplicity. Likewise $Z'(s) = \sum e^{-\lambda_k(T')s}$. (This is the trace of the solution kernel of the heat equation for $T$. The theory in general may be found in Rudin \cite{Rudin1987} 13.34-13.38; see also Taylor, \cite{Taylor2011a} Ch. 6 and \cite{Taylor2011b} Ch 8.) In particular if $\lambda_k(T) \leq \lambda_k(T')$ then for all $s>0$ we have $e^{-\lambda_k(T)s}\geq e^{-\lambda_k(T')s}$ and so $Z(s) \geq Z'(s)$.

The inequality in the second part of the theorem is proven in a fashion identical to the previous theorem's proof.
\end{proof}

The following lemma is a well-known argument (see for example the proof of Huber's theorem in Buser \cite{Buser1992} 9.2.9) and forms part of the proof of several propositions in Chapter 3.

\begin{lemma}\label{heat_ineq_lemma}
Suppose we have
\[ F(t) = \sum_{i=1}^\infty a_ie^{-r_it} \]
and
\[ G(t) = \sum_{i=1}^\infty b_ie^{-s_it} \]
where $c_i, b_i, r_i, s_i$ are such that $F$ and $G$ converge for all $t>0$ and the sequences $r_i$ and $s_i$ are decreasing with finite multiplicity.

If $F \leq G$ then $r_1 \geq s_1$.
\end{lemma}
\begin{proof}
Suppose $F \leq G$. That is, for arbitrary $t>0$, we have
\[ \sum_{i=1}^\infty a_i e^{-r_it} \leq \sum_{i=1}^\infty b_i e^{-s_it}. \]
Factor $e^{-r_1t}$ from the left and $e^{-s_1t}$ from the right and take the logarithm of each side, then divide through by $t$:
\[ -r_1 + \frac{1}{t} \log\bigg( \sum a_i e^{(r_1-r_i)t}  \bigg) \leq \frac{1}{t} \log\bigg( \sum b_i e^{(s_1 - s_i)t} \bigg) \]
Notice that as $t\to\infty$, as $r_1$ is less than all but finitely many $r_1$ and $s_1$ is less than all but finitely many $s_i$, the argument of the logarithm on the left hand side of the inequality tends to the product of $a_1$ and the multiplicity of $r_1$, and the argument of the logarithm on the right hand side of the inequality tends to the product of $a_1$ and the multiplicity of $s_1$.

Thus the left hand side tends to $-r_1$ and the right hand side tends to $-s_1$; multiplying through by $-1$ yields the result.
\end{proof}

We have the immediate corollary
\begin{cor}\label{heat_ineq_cor}
Suppose $F$ and $G$ are as in the statement of Lemma \ref{heat_ineq_lemma}. Suppose $F \leq G$ If $c_i = b_i = 1$ and $r_i = s_i$ for all $i= 1,\ldots, k-1$, then $r_k \geq s_k$.
\end{cor}
\begin{proof}
Cancel the first $k-1$ terms from each side of the inequality $F \leq G$, then apply Lemma \ref{heat_ineq_lemma}.
\end{proof}

\chapter{Laplace subspectrality}%
\section{Riemannian geometry}

Suppose $M$ is a smooth oriented manifold of dimension $n$ possibly with Lipschitz, piecewise-smooth boundary. Recall that $M$ is equipped with a tangent bundle $TM$ and a cotangent bundle $T^*M$. A function $f:M\to\mb{R}$ is said to be differentiable provided in every local coordinate expression $x_i$ of $f$, the derivatives $\partial_if$ exist. The set of continuous functions $M\to\mb{R}$ is denoted by $C^0(M)$ and the set of $k$-times continuously differential functions is denoted by $C^k(M)$. This condition is defined inductively: for $k>0$ function is $k$-times continuously differentiable provided it is differentiable and in every set of local coordinates all of its derivatives are $(k-1)$-times continuously differentiable. A function is of class $C^\infty(M)$ provided it is of class $C^k$ for all $k\geq 0$. The exterior derivative $d$ acts on smooth functions by the local coordinate expression $df = \sum\partial_i f dx^i.$

If $M$ is a Riemannian manifold with a smooth positive definite metric $g$ we have a volume form $dv_g$. In local coordinates $x_i$ the metric $g$ has the expression $g_{ij} = g(\partial_i,\partial_j)$ and the volume form has expression $dv_g = \sqrt{|g|}dx^1\wedge\cdots\wedge dx^n$. Gram-Schmidt provides the existence of orthonormal frame fields. An orthonormal frame field is a set of locally defined vector fields $e_i$ such that $g(e_i,e_j) = \delta_{ij}$ pointwise and the $e_i$ span $T_pM$ for all $p$ in the neighborhood.

The metric gives a bundle isomorphism between $TM$ and $T^*M$ defined by $X\mapsto (V\mapsto g(X,V))$. By composing this with the exterior derivative we have the gradient operator $\nabla: C^1(M)\to\mf{X}(M)$ defined by $g(X,\nabla f) = df(X)$. The metric also gives a covariant derivative operator on vector fields mapping a vector field $X$ to the linear operator $v\mapsto D_vX$ on each tangent space. In local coordinates, the gradient is given by $\nabla f = g^{ij}\partial_jf$ where $g^{ij}$ denotes the inverse of the matrix $g_{ij}$. Given a smooth vector field $X$, the divergence of $X$ is defined to be the trace of the covariant derivative of $X$ on each tangent space.

The Laplace operator, or Laplacian, $\Delta$ is defined as the divergence of the gradient. In local coordinates it is expressed as 
\[ \Delta f = -\frac{1}{\sqrt{|g|}}\sum_{i,j} \partial_i (\sqrt{|g|}g^{ij}\partial_jf) \]
In a domain in Euclidean space $\mb{R}^n$ this reduces to the expression $\sum_i \partial_i^2$.

When the metric is clear from context, we refer to $\Delta$, $dv$, and $\nabla$.

We make the following useful definition:
\begin{defn}[Normal manifold]
We say a compact Riemannian manifold with (possibly empty) piecewise smooth, Lipschitz boundary is a normal manifold. If the manifold is a codimension zero submanifold of $\mb{R}^n$, we say it is a normal domain.
\end{defn}
We often refer to manifolds as $M$ and to domains as $\Omega$.

\section{Laplace operator eigenvalue problem}

Suppose $M$ is a normal manifold. We consider the eigenvalue problem
\[ \Delta f = \lambda f \]
If $\partial M = \emptyset$ then we consider the closed eigenvalue problem. If $M$ has boundary, we consider the following different boundary conditions:
\begin{itemize}
\item Dirichlet (function restricted to boundary is zero),
\item Neumann (outward normal derivative is zero), 
\item mixed (partition $\partial M$ into finitely many subsets which are open in $\partial M$ and impose Dirichlet or Neumann conditions on each of those subsets).
\end{itemize}

If $M$ has boundary, we define boundary conditions in the following fashion. Denote by $\nu$ the indicator function of the interior of a normal (not necessarily connected) submanifold of $\partial M$. We set $D^\nu$ to be those smooth functions on $M$ whose outward normal derivative is equal to zero on $\nu^{-1}(1)$ and whose support does not intersect $\nu^{-1}(0)$, and $\Delta^\nu$ the Laplace operator restricted to the domain $D^\nu$.

In particular, Dirichlet boundary conditions are described by the function $\nu \equiv 0$ and Neumann boundary conditions are described by $\nu\equiv 1$. When convenient and clear from context, we shall conflate the function $\nu$ with the appropriate integer.

\subsection{Quadratic forms and the spectrum of the Laplacian}

By Green's theorem (see Chavel \cite{Chavel1984} I.10), in the case that $M$ has boundary, $\Delta^\nu$ is symmetric on $D^\nu$. Applying the Friedrichs extension Theorem \ref{friedrichs} to $\Delta^\nu$ with $\sigma=1$, we construct the form domain $V^\nu$. By Proposition \ref{spectral_thm_laplacian} the form domain $V^\nu$ compactly embeds in $L^2(M)$ and so the Friedrichs extension of the Laplacian $\Delta^\nu$ is compactly resolved and has discrete spectrum.

In the case that $M$ does not have boundary, it is a complete Riemannian manifold and thus, by a result of Gaffney \cite{Gaffney1951}, the Laplacian defined on smooth functions is essentially self-adjoint, and its Friedrichs extension is compactly resolved and has discrete spectrum.

Note that if $M$ has boundary, our definition of ``normal'' requires $M$ have piecewise smooth boundary. In fact, for manifolds which isometrically embed in Euclidean space, the Neumann Laplacian is compactly resolved provided $M$ has the ``segment property.'' For discussion see Reed \& Simon \cite{Reed1978} Vol IV Ch XIII Section 14. For manifolds which do not have the segment property, the Neumann Laplacian need not be compactly resolved; c.f. ibid, the immediately preceding section, for an example.

The form associated to the Laplacian $\Delta^\nu$ is $\mf{q}_\nu$, defined by
\[ \mf{q}_\nu(u,v) = \int g(\nabla u, \nabla v)\ dv \]
and the inner product on $V^\nu$ is $(u,v)_\nu = (u,v) + \mf{q}_\nu(u,v)$.

We therefore have a Rayleigh quotient defined on $V^\nu$:
\[ \mc{R}^\nu(u) = \frac{\mf{q}(u,u)}{\|u\|^2} \]
and we may apply to $\mf{q}_\nu$ and $\mc{R}^\nu$ the theory of Chapter 1.

When the boundary conditions are specified in context, we will drop the subscript or superscript $\nu$ and refer to $\mf{q}$, $\mc{R}$, and $V$.

We may now make the following definitions to apply the theory of Chapter 1.

\begin{defn}[Sub/superspectrality of manifolds]
Suppose $M$ and $M'$ are normal manifolds. We say that $M$ is Dirichlet-(resp. Neumann-) subspectral to $M'$ provided $\Delta^0(M)$(resp. $\Delta^1(M)$) is subspectral to $\Delta^0(M')$(resp. $\Delta^1(M')$).

When the same boundary conditions are applied to the Laplacian on $M$ and $M'$ and they are clear from context, or if $M$ and $M'$ are both closed, we may simply say that $M$ is subspectral, or Laplace subspectral, to $M'$.
\end{defn}

The other definitions from Chapter 1 carry over as well by substituting $M$ for $\Delta^\nu$.

To illustrate the definition of subspectrality, we prove a simple statement:
\begin{prop}[Finite Riemannian covers]
Suppose $\hat{M}\to M$ is a finite-sheeted Riemannian cover of closed manifolds. Then $\hat{M}$ is subspectral to $M$.
\end{prop}
\begin{proof}
Every eigenfunction on $M$ lifts to an eigenfunction on $\hat{M}$. Thus the spectrum of $M$ is a subset of the spectrum of $\hat{M}$ which immediately implies that $\hat{M}$ is subspectral to $M$.
\end{proof}

We make a definition for the Weyl functions of a manifold. These are polynomials in $\sqrt{\lambda}$.
\begin{defn}\label{weylfunc}[Weyl functions]
Let $M$ be a normal $n$-manifold possibly with boundary. Denote by $\omega_n$ the volume of the unit ball in $\mb{R}^n$. We define the one-term Weyl function of $M$ to be the following:
\[ w_1^M:\lambda \mapsto \frac{\omega_n|M|}{(2\pi)^n}\lambda^{n/2}. \]
If $M$ has boundary and we have chosen boundary conditions denoted by $\nu=0,1$ for Dirichlet and Neumann, respectively, we define the two-term Weyl function to be
\[ w_2^M: \lambda \mapsto w_1^M(\lambda) + (-1)^{\nu+1}\frac{1}{4}\frac{\omega_{n-1}|\partial M|}{(2\pi)^{n-1}}\lambda^{\frac{n-1}{2}} \]
When $M$ is clear from context we may omit the superscript or subscript in $w_k^M$.
\end{defn}

\section{Domain monotonicity}

There are three well-known monotonicity theorems for the Laplacian. Neumann and Dirichlet monotonicity can be found in Chavel \cite{Chavel1984} Chapter I Section 5, immediately following the Max-Min Theorem. Generalizing the functional analytic ideas in their proofs to Theorem \ref{subspectr1} allows us to prove them all as corollaries.

\begin{cor}[Neumann is subspectral to Dirichlet]
Let $M$ be a normal manifold with boundary. Then the Neumann Laplacian $\Delta^1$ is subspectral to the Dirichlet Laplacian $\Delta^0$.
\end{cor}
\begin{proof}
All compactly supported smooth functions have vanishing normal derivative. Therefore we have inclusion of domains $V^0\subset V^1$ in $L^2(M)$ and on $V^0$ we have $\mf{q}_0 = \mf{q}_1$. The result follows from Theorem \ref{subspectr1}.
\end{proof}

The following result concisely expresses the Dirichlet and Neumann domain monotonicity theorems.
\begin{thm}[Partition theorem]\label{partitionthm}
Let $M$ be a normal manifold. Let $\{\Gamma_i\}_{i=1}^N$ be a finite partition of $M$ by normal, codimension zero manifolds. Impose boundary conditions $\nu$ on $M$. On $\partial M\cap\cup_i\partial\Gamma_i$, impose boundary conditions by restricting $\nu$.

The internal boundaries of the $\Gamma_i$ are $\cup_i \partial\Gamma_i - \partial M$. If we impose Dirichlet conditions on the internal boundaries of the $\Gamma_i$, then $M$ is subspectral to $\sqcup_i\Gamma_i$. If we impose Neumann conditions on the internal boundaries of the $\Gamma_i$, then $M$ is superspectral to $\sqcup_i\Gamma_i$.
\end{thm}

\begin{proof}
Let $H = L^2(M)$, and $H_i = L^2(\Gamma_i)$. Give $\oplus_i H_i$ the natural Hilbert space structure. Let $\nu_i$ denote the boundary conditions imposed on $\Gamma_i$. Let $V$ be the form domain of the $\nu$-Laplacian on $M$ and let $D^{\nu_i}$ and $V_i$ be the domain and form domain, respectively, of the $\nu_i$-Laplacian on $\Gamma_i$. Denote by $\mf{q}_\oplus$ the energy form on $\oplus_i V_i$ and denote by $\mf{q}$ the energy form on $V$.

Let $R: H\to \oplus_iH_i$ be the restriction map $Ru = u|_{\Gamma_1}\oplus\cdots\oplus u|_{\Gamma_N}$. Let $W$ be the domain in $H$ of $R^*\mf{q}_\oplus$. The key observation in the proof is that if we define the $\nu_i$ by imposing Dirichlet conditions on the internal boundaries of the $\Gamma_i$, then $W\subseteq V$, while if we define the $\nu_i$ by imposing Neumann conditions on the internal boundaries of the $\Gamma_i$< then $V\subseteq W$.

First we show that $R$ is a bijective isometry of Hilbert spaces. Then we show that on $W\cap V$, we have $\mf{q} = R^*\mf{q}_\oplus$. Finally we establish the respective inequalities resulting from the imposition of Neumann and Dirichlet boundary conditions.

To see that $R$ is a bijective isometry, observe that if $u\in H$, then
\[ \|Ru\|^2 = \sum_i \|u\|_{L^2(\Gamma_i)}^2 = \sum_i \int_{\Gamma_i} u|_{\Gamma_i}^2\ dv = \int_M u^2\ dv = \|u\|^2. \]
For $v\in H_i$, define
\[ \bar{v}(x) = \begin{cases} v(x),\ &x\in\Gamma_i\\ 0,\ &x\notin\Gamma_i \end{cases} \]
If $u_1\oplus \cdots\oplus u_N\in \oplus_i H_i$, then $R(\sum_i\bar{u_i}) = \oplus_i u_i$. So $R$ is surjective and its inverse is defined by mapping $\oplus_iu_i\mapsto \sum_i\bar{u}_i$. If $Ru = Rv$, then $R(u-v) = 0$ and as $R$ is an isometry, we must have $u=v$ in $H$. Thus $R$ is a bijective isometry of Hilbert spaces. By Lemma \ref{inj_isom_subsp} $R^*\mf{q}_\oplus$ is isospectral to $\mf{q}_\oplus$, so $\mf{q}$ is subspectral to $\mf{q}_\oplus$ iff $\mf{q}$ is subspectral to $R^*\mf{q}_\oplus$, and likewise $\mf{q}$ is superspectral to $\mf{q}_\oplus$ iff $\mf{q}$ is superspectral to $R^*\mf{q}_\oplus$.

Now suppose $u,v\in W\cap V$. Then we have
\begin{align*}
R^*\mf{q}_{\oplus}(u,v) &= \mf{q}_\oplus(Ru, Rv)\\
	&= \sum_i \int_{\Gamma_i} \langle\nabla (u|_{\Gamma_i}), \nabla (v|_{\Gamma_i})\rangle\ dv\\
	&= \int_M \langle \nabla u, \nabla v\rangle\ dv\\
	&= \mf{q}(u, v)
\end{align*}
So $\mf{q} = R^*\mf{q}_\oplus$ on $W\cap V$.

Suppose we have imposed Dirichlet conditions on the interior boundaries of the $\Gamma_i$. Any $u\in\oplus_i D^{\nu_i}$ is supported away from the internal boundary of $\Gamma_i$, so the image of $\oplus_i D^{\nu_i}$ under $R^{-1}$ is contained within $D^\nu$. Because $\mf{q} = R^*\mf{q}_\oplus$ we have $W\subseteq V$ and the result for internal Dirichlet conditions follows from Theorem \ref{subspectr1}.

Suppose we have imposed Neumann conditions on the interior boundaries of the $\Gamma_i$. As the restriction of an element of $V$ to $\Gamma_i$ is an element of $V_i$ we have that $V\subseteq W$. The result for internal Neumann conditions follows from Theorem \ref{subspectr1}.
\end{proof}

For Neumann conditions on the interior boundaries of the partition sets, this is the well-known Neumann domain monotonicity theorem. The well-known Dirichlet domain monotonicity theorem follows from a short proof.

\begin{cor}[Dirichlet domain monotonicity]\label{dirichlet_monotonicity}
Let $M$ be a Riemannian manifold with piecewise smooth boundary. Let $\Gamma_i\subset M$ be finitely many pairwise disjoint subdomains. Denote by $\sigma_k$ the ordering with multiplicity of the union of the Dirichlet spectra of the $\Gamma_i$. Then $M$ is subspectral to the disjoint union of the $\Gamma_i$.
\end{cor}
\begin{proof}
Let
\[\Gamma_0 = M - \bigg(\bigcup_i \Gamma_i\bigg). \]
By Corollary \ref{partitionthm}, and then observing that counting functions are always positive, we have
\[ N_M \geq \sum_{i\geq 0}N_{\Gamma_i} \geq \sum_{i\geq 1}\Gamma_i \]
which is as desired.
\end{proof}

We have the following conjecture:
\begin{conj}\label{neumann_monotonicity}
Let $M$ be a normal manifold and let $\{\Gamma_i\}$ be a finite open cover of $M$ by normal manifolds. Then, with interior boundary components given Neumann conditions, $\sqcup_i\Gamma_i$ is subspectral to $M$.
\end{conj}

By Weyl's law, because $\sum_i |\Gamma_i| > |M|$, we must have that $\sqcup_i\Gamma_i$ is eventually subspectral to $M$. Thus any counterexamples must be low-order eigenvalues. By the quantitative Weyl law proven below, this conjecture can be verified with numerical computation for any given compact Euclidean domain and open cover.


We include Propositions \ref{neumann_general_1} and \ref{neumann_general_2} as work toward this conjecture.

\begin{defn}\label{l2part1}If $M$ is a compact Riemannian manifold with piecewise smooth boundary and $\{\Gamma_i\}$ is an open cover of $M$, we say that $\Phi = \{\phi_i\}$ is an $L^2$ partition of unity subordinate to $\{\Gamma_i\}$ provided the $\phi_i$ are positive real-valued functions on $M$, smooth in the interior of $M$, and the collection $\phi_i^2$ is a partition of unity subordinate to $\{\Gamma_i\}$.\end{defn}

The existence of these follows by taking the positive square roots of partitions of unity subordinate to the same cover. The existence of partitions of unity is established in standard differential topology introductions such as Warner \cite{Warner1961} 1.11.

\begin{prop}[Generalized Neumann monotonicity 1]\label{neumann_general_1}
Let $M$ be a compact Riemannian manifold with piecewise smooth possibly empty boundary. Let $\{\Gamma_i\}$ be a finite open cover of $M$. Denote by $\sigma_k$ the ordered collection of Neumann eigenvalues of the $\Gamma_i$, counted with multiplicity, and denote by $\lambda_k$ the Neumann eigenvalues of $M$.

For an $L^2$ partition of unity $\Phi = \{\phi_i\}$ subordinate to $\{\Gamma_i\}$ we set
\[\mc{D}_\Phi = \sup_M \sum_i |\nabla\phi_i|^2\]
and $\mc{D} = \inf_\Phi D_\Phi$.
Then for each $k$,
\[ \sigma_k \leq \mu_k + \mc{D} \]
\end{prop}
\begin{proof}
Denote by $V$ the form domain of the Neumann Laplacian on $M$ and by $V_i$ the form domain of the Neumann Laplacian on $\Gamma_i$. Denote by $H$ the space $L^2(M)$ and by $H_i$ the space $L^2(\Gamma_i)$. Denote by $\mf{q}$ the form of the Neumann Laplacian on $M$ and by $\mf{q}_\oplus$ the form of the orthogonal sum of the Neumann Laplacians of $\Gamma_i$ on its form domain $\oplus_i V_i$.

For convenience in computation, make the following notation. For two mult-indices $I = i_1\cdots i_k$, $J=j_{k+1}\cdots j_N$  such that $(i_1,\ldots,i_k,j_{k+1},\ldots,j_N)$ is a permutation of $(1,2,\ldots,N)$, we say 
\[\Gamma_{IJ} = (\Gamma_{i_1}\cap\cdots\cap\Gamma_{i_k})\cap(\Gamma_{j_{k+1}}\cup\cdots\cup\Gamma_{j_N})'\]
is a leaf set of $\{\Gamma_i\}$. Note that for any integrable function $f$ we have
\[\int_M f\ dx = \sum_{\mbox{leaf sets }\Gamma_{IJ}} \int_{\Gamma_{IJ}} f\ dx\]

The weighted restriction map $\Phi:H\to \oplus_i H_i$ defined by
\[\Phi: v\mapsto \oplus_i\phi_i\cdot v|_{\Gamma_i}\]
is an isometry. It maps $D^1$ into $\oplus_i D^1$, hence maps $V$ into $\oplus_i V_i$. The map $\Phi$ is not an isometry when restricted to $V$, as we compute. Let $u\in V$. Then:
\begin{align*}
\mf{q}_\oplus(\Phi u) &= \sum_i \int_{\Gamma_i} |\nabla (\phi_i u)|^2\ dv \\
	&= \sum_i \int_{\Gamma_i} (|\nabla \phi_i|^2 |u|^2 + 2u\phi_i\langle \nabla u, \nabla\phi_i\rangle + \phi_i^2|\nabla u|^2\ dv\\
	&= \sum_{\mbox{leaf sets}} \int_{\Gamma_{IJ}} \bigg( |\nabla u|^2 \sum_{i\in I} \phi_i^2 + \langle \nabla u, \sum_{i\in I}2\phi_i\nabla\phi_i\rangle\\ &\ \ \ \ + |u|^2\sum_{i\in I}|\nabla \phi_i|^2\ \bigg)\ dv
\end{align*}
On a leaf set $\Gamma_I$, we have $\sum_{i\in I}\phi_i^2 = \sum_i\phi_i^2 = 1$ because for $i\notin I$ we have $\phi_i = 0$. Thus the first term in each integral is $|\nabla u|^2$.

Taking the gradient of both sides of $\sum_{i\in I}\phi_i^2 = 1$ we have $\sum_{i\in I} 2\phi_i\nabla \phi_i = 0$, so the second term vanishes. 

By the observation that integration over $M$ is equal to a sum of integrals over the leaf sets, we have
\[ Q_\oplus(\Phi u) = Q_M(u) + \delta_\Phi(u) \]
where we define
\[ \delta_\Phi(u) = \int_M |u|^2\sum_i |\nabla \phi_i|^2\ dv. \]

Let $\hat{E}_k(\Gamma)$ denote the span in $\oplus_i V_i$ of the first $k$ eigenfunctions of $\mf{q}_\oplus$. Let $\hat{E}_k(M)$ denote the span in $V$ of the first $k$ eigenfunctions of $\mf{q}$. Because $\Phi$ is injective the image under $\Phi$ of $\hat{E}_k(M)$ is $k$-dimensional, hence there exists some vector $v_\phi\in \Phi(\hat{E}_k(M))$ which is perpendicular to $\hat{E}_{k-1}(\Gamma)$.

We therefore have
\[ \sigma_k\|v_\Phi\|^2 \leq \mf{q}_\oplus(v_\Phi) = \mf{q}(v_\Phi) + \delta_\Phi(v_\Phi) \leq \lambda_k\|v_\Phi\|^2 + \delta_\Phi(v_\Phi). \]

Applying H\"older's inequality, $\|fg\|_1 \leq \|f\|_p\|g\|_q$ for $p^{-1} + q^{-1} = 1$, to $\delta_\Phi(v_\Phi)$ with $p=1,q=\infty$, and dividing through by $\|v_\Phi\|^2$, gives
\[ \sigma_k \leq \lambda_k + \frac{\delta_\Phi(v_\Phi)}{\|v_\Phi\|^2} \leq \lambda_k + \sup_M \sum_i |\nabla\phi_i|^2. \]
Taking the infimum over all $L^2$ partitions of unity yields the claimed result.
\end{proof}

If additional information could be deduced about the function $v_\Phi$, this argument might be extended to prove Conjecture \ref{neumann_monotonicity}. For example, if for each $k$ one could find a sequence of $L^2$ partitions of unity $\phi_m$ such that $\delta_{\Phi_m}(v_{\Phi_m})\to 0$, the conjecture would be established.

We remark that we may not apply Theorem \ref{subspectr1} to the proof of Proposition \ref{neumann_general_1} because the map $\Phi$ is not an isometry on form domains.

In a similar vein we have:
\begin{prop}[Generalized Neumann monotonicity 2]\label{neumann_general_2}
Let $M$ be a normal manifold. Let $\{\Gamma_i\}$ be a finite collection of codimension zero normal submanifolds of $M$ such that $M \subset \cup_i\Gamma_i$. Denote by $\sigma_k$ the ordered collection of Neumann eigenvalues of the $\Gamma_i$, counted with multiplicity, and denote by $\lambda_k$ the Neumann eigenvalues of $M$.

Let $G = \sup_{x\in M}|\{\Gamma_i\ |\ x\in\Gamma_i\}|$. Then $\sigma_k \leq (1+G)\lambda_k$.
\end{prop}
\begin{proof}
Denote by $V$ the form domain of the Neumann Laplacian on $M$ and by $V_i$ the form domain of the Neumann Laplacian on $\Gamma_i$. Denote by $H$ the space $L^2(M)$ and by $H_i$ the space $L^2(\Gamma_i)$. Denote by $\mf{q}$ the form of the Neumann Laplacian on $M$ and by $\mf{q}_\oplus$ the form of the orthogonal sum of the Neumann Laplacians of $\Gamma_i$ on its form domain $\oplus_i V_i$. Define the leaf sets $\Gamma_{IJ}$ of the cover $\{\Gamma_i\}$ as in the proof of Proposition \ref{neumann_general_1}.

Define the restriction map $R:H\to\oplus_i H_i$ by $R(u) = u|_{\Gamma_1}\oplus\cdots\oplus u|_{\Gamma_n}$. Notice that $R$ is injective as the $\Gamma_i$ form an open cover and $R$ takes $V$ into $\oplus_i V_i$. Denote by $W$ the domain of the pullback $R^*\mf{q}_\oplus$. Because $R$ takes $V$ into $\oplus_i V_i$, we have that $V\subset W$.

Suppose $u\in W\cap V$. Then
\begin{align*}
\mf{q}_\oplus(R(u)) &= \sum_i \int_{\Gamma_i} |\nabla u|^2\ dv\\
	&= \sum_{\mbox{leaf sets }\Gamma_{IJ}} \sum_{i\in I}\int_{\Gamma_i}|\nabla u|^2\ dv\\
	&= \int_M|\nabla u|^2\ dv + \sum_{|I|>1}|I|\int_{\Gamma_i} |\nabla u|^2\ dv \\
	&\leq \int_M|\nabla u|^2\ dv + G\int_M|\nabla u|^2\ dv \\
	&= (1 + G)\mf{q}(u)
\end{align*}
As $R^*\mf{q}_\oplus\leq (1+G)\mf{q}$ and $V\subset W$, by Theorem \ref{subspectr1} we have that $R^*\mf{q}_\oplus$ is subspectral to the form $(1+G)\mf{q}$. 

Suppose $u\in H$. Then 
\[ \|Ru\|_{\oplus_i H_i}^2 = \sum_i \int_{\Gamma_i} |u|^2\ dv \geq \int_M |u|^2\ dv = \|u\|_H^2\]
so $R$ satisfies the conditions of Lemma \ref{pullback_subsp}. Therefore we have that $\mf{q}_\oplus$ is subspectral to $R^*\mf{q}_\oplus$.

Since the spectrum of $(1+G)\mf{q}$ is obtained by multiplying each Neumann eigenvalue of $M$ by $(1+G)$, we have established the desired result.
\end{proof}

We remark that Neumann monotonicity is a corollary to Proposition \ref{neumann_general_2}.


\section{Quantitative Weyl law}

In this section we prove a quantitative Weyl law to use in constructing examples of pairs of domains where one is subspectral to the other. We use the term quantitative Weyl law because instead of proving asymptotic isospectrality, these results give subspectrality and superspectrality of the Laplace operator to modified Weyl polynomials. Weyl's law for Euclidean domains follows as a corollary.

The following lemma is based on a lattice counting argument attributed to Gauss; the earliest reference the author could find is Rayleigh \cite{Rayleigh1900}. The argument bounds pointwise error terms.

\begin{lemma}[Quantitative Weyl law for rectangles]
Let $a_1,\ldots,a_n$ be positive numbers and let $R = \times_{j=1}^n [0,a_j]$ denote the rectangular prism (unique up to isometry) with side lengths $a_i$ in $\mb{R}^n$. Let $N_\nu$ denote the counting function of $\Delta_\nu$ on $R$ where $\nu=0$ represents Dirichlet conditions and $\nu=1$ represents Neumann conditions. Denote by $d$ the codiagonal of $R$, defined by $d^2 = \sum_i \frac{\pi^2}{a_i^2}$. Denote by $w$ the one-term Weyl polynomial of $R$.

Then:
\begin{itemize}
\item $R$ is Neumann-subspectral to $w$ and Neumann-superspectral to $w\cdot (\lambda\mapsto 1 + d/\sqrt{\lambda})^n$
\item $R$ is Dirichlet-superspectral to $w$ and Dirichlet-subspectral to $w\cdot(\lambda\mapsto 1 - d/\sqrt{\lambda})^n$
\end{itemize}
\end{lemma}
\begin{proof}First note that by Definition \ref{SubspFunc} it suffices to bound $N_\nu$ above and below.

Recall that the eigenvalues of $\Delta_\nu$ are of the form
\[ \lambda_{i_1i_2\cdots i_n} = \bigg(\frac{i_1\pi}{a_1}\bigg)^2 + \cdots + \bigg(\frac{i_n\pi}{a_n}\bigg)^2 \]
where each $i_j$ ranges over positive integers if $\nu=0$ and nonnegative integers if $\nu=1$. Denote by $Q_1$ the closed first quadrant and $Q_0$ the open first quadrant. Each eigenvalue of $\Delta_\nu$ corresponds to the squared length of exactly one element of $Q_\nu\cap (\oplus_j (\pi/ a_j\mb{Z})$. For ease of notation, set $\Lambda_a = \oplus_j (\pi/a_j)\mb{Z}$ and denote by $B(0,r)$ the ball of radius $r$ in $\mb{R}^n$. 

The eigenvalue counting function $N_\nu$ of $\Delta_\nu$ satisfies
\[ N_\nu(\lambda) = \big|Q_\nu\cap B(0,\sqrt{\lambda}) \cap \Lambda_a\big|.\]

We estimate $N_0$. Associate to each $v\in Q_\nu\cap B(0,\sqrt{\lambda}) \cap \Lambda_a$ the cell of $\Lambda_a$ whose vertices are $v, v+(\pi/a_1,0,\ldots,0), v+(0,\pi/a_2,\ldots,0),\ldots, v+(0,0,\ldots,\pi/a_n)$. (This is the cell of $\Lambda_a$ which is closest to the origin of those cells adjacent to $v$.) The volume of each cell of $\Lambda_a$ is $\pi^n/a_1a_2\cdots a_n = \pi^n/|R|$.

By comparing the areas of cells of $Q_0\cap \Lambda_a\cap B(0,\sqrt{\lambda})$ to the volume of the hemisphere $S(\sqrt{\lambda})\cap Q_0$ we have the following relation:

\[ |S(\sqrt{\lambda} - d)\cap Q_0|\leq N_0(\lambda)\frac{\pi^n}{|R|}\leq |S(\sqrt{\lambda})\cap Q_0| \]

As the volume of $B(0,1)$ is equal to $\omega_n$, we have $|S(\lambda)\cap Q_\nu| = \frac{\omega_n}{2^n}\lambda^n$ for $\nu=0,1$. Thus:
\[ \frac{\omega_n}{2^n}\bigg(\sqrt{\lambda} - d\bigg)^n \leq N_0(\lambda)\frac{\pi^n}{|R|} \leq  \frac{\omega_n}{2^n}\lambda^{n/2}\]

Simplifying, we have the desired result for the Dirichlet Laplacian on $R$:
\[\frac{\omega_n|R|}{(2\pi)^n}\bigg(1 - \frac{d}{\sqrt{\lambda}}\bigg)^n\lambda^{n/2} \leq N_0(\lambda) \leq  \frac{\omega_n|R|}{(2\pi)^n}\lambda^{n/2}\]

We estimate $N_1$. To each lattice point $v\in Q_1\cap \Lambda_a$ associate the cell furthest from the origin of those adjacent to $v$. The hemisphere $S(\sqrt{\lambda})$ is contained in the union of these cells, so we have the area estimate
\[ \frac{\omega_n}{2^n}\lambda^{n/2} \leq N_1(\lambda)\frac{\pi^n}{|R|} \leq  \frac{\omega_n}{2^n}(\sqrt{\lambda} + d)^n \]
which gives
\[ \frac{\omega_n|R|}{(2\pi)^n}\lambda^{n/2} \leq N_1(\lambda) \leq  \frac{\omega_n|R|}{(2\pi)^n}\bigg(1 + \frac{d}{\sqrt{\lambda}}\bigg)^n\lambda^{n/2} \]
as desired.


\end{proof}

\begin{figure}[h]
	\centering
	\includegraphics[scale=0.75]{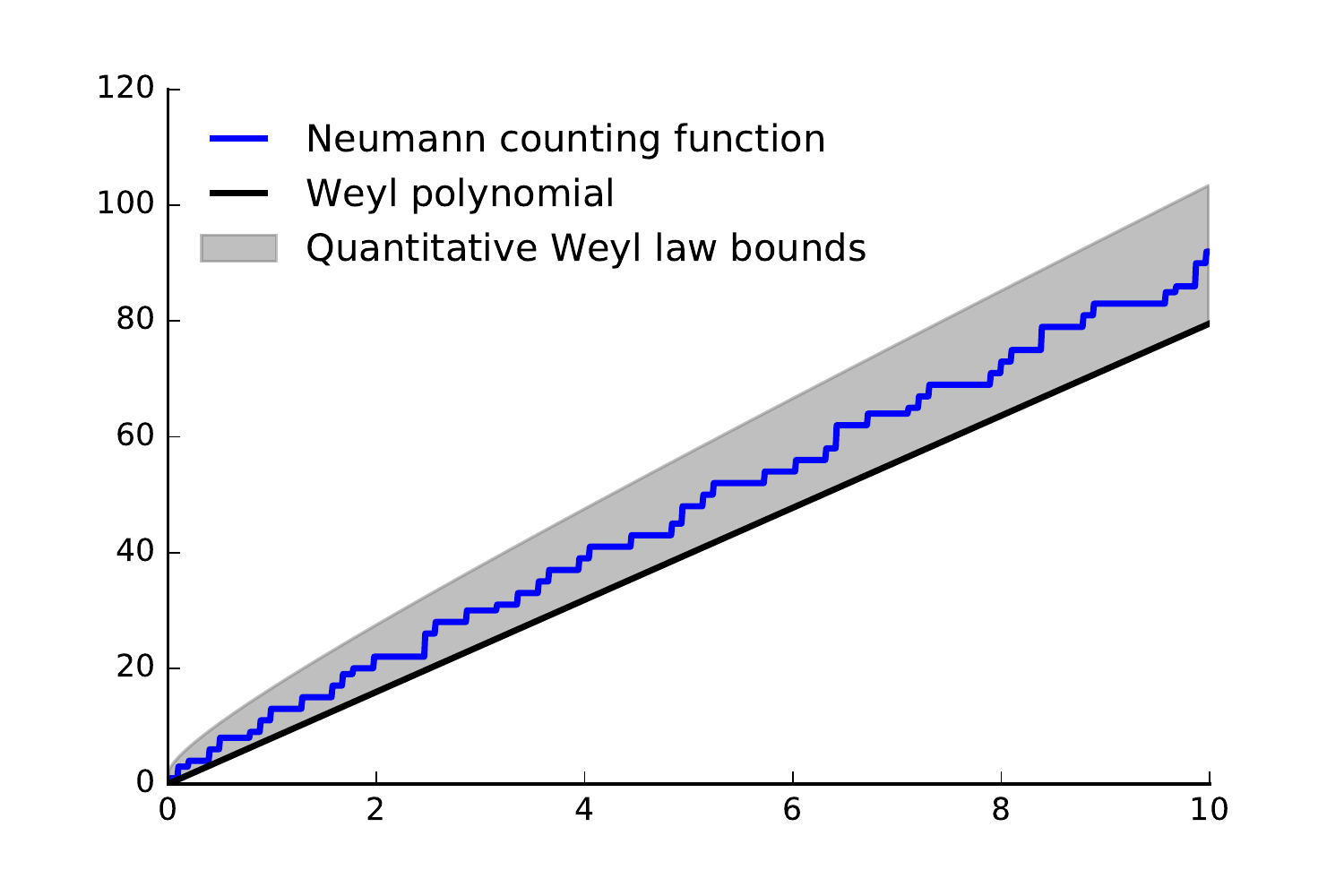}
	\caption[Bounds for quantitative Weyl law]{The shaded area represents the bounds given by the quantitative Weyl law for rectangles.}
	\label{quantitative_weyl_rect}
\end{figure}

We illustrate with the graphic in Figure \ref{quantitative_weyl_rect}, depicting the Neumann counting function of the square $[0,10]\times[0,10]$. 

We have the following statement as an immediate corollary.
\begin{cor}[Quantitative Weyl law for squares]
If $R$ is a square of side length $\eps>0$ in $\mb{R}^n$ then the Neumann Laplacian is superspectral to 
\[x\mapsto w^R(\lambda)\bigg(1 + \frac{\pi\sqrt{n}}{\eps\sqrt{\lambda}}\bigg)^n \]
and subspectral to $w^R$.
The Dirichlet Laplacian is subspectral to 
\[ \lambda\mapsto w^R\bigg(1 - \frac{\pi\sqrt{n}}{\eps \sqrt{\lambda}}\bigg)^n \]
and superspectral to $w^R$
\end{cor}
\begin{proof}
By substituting $a_i=\eps$ in the previous lemma, we have $|R| = \eps^n$ and $d =\frac{\pi}{\eps}\sqrt{n}$. The Dirichlet counting function $N_0$ then satisfies
\[ \frac{\omega_n\eps^n}{(2\pi)^n}\bigg(1 - \frac{\pi\sqrt{n}}{\eps\sqrt{\lambda}}\bigg)^n\lambda^{n/2} \leq N_0(\lambda) \leq  \frac{\omega_n\eps^n}{(2\pi)^n}\lambda^{n/2} \]
and the Neumann counting function $N_1$ satisfies
\[ \frac{\omega_n\eps^n}{(2\pi)^n}\lambda^{n/2} \leq N_1(\lambda) \leq  \frac{\omega_n\eps^n}{(2\pi)^n}\bigg(1 + \frac{\pi\sqrt{n}}{\eps\sqrt{\lambda}}\bigg)^n\lambda^{n/2} \]
as desired.
\end{proof}

In general we establish a quantitative Weyl law for the Dirichlet Laplacian on Euclidean domains. The argument is an adaptation of the proof of Weyl's law due to Weyl \cite{Weyl1911} as described in Courant-Hilbert \cite{Courant1953}. We modify it to track the error bounds on the counting function of each square.

\begin{thm}[Quantitative Weyl law for Euclidean domains]\label{quantitative_weyl}
Let $\Omega$ be a normal domain in $\mb{R}^n$. For any $\eta > 0$ denote by $\Omega^{-\eta}$ the set of points in $\Omega$ of distance greater than $\eta$ from $\partial\Omega$ and denote by $\Omega^{\eta}$ the set of points in $\mb{R}^n$ of distance no more than $\eta$ from any point of $\Omega$.

Define the function $E^\eps_\pm$ by
\begin{align*}
E^\eps_\pm(\lambda) &= w^{\Omega^{\pm\eps\sqrt{n}}}(\lambda)\bigg(1 \pm \frac{\pi\sqrt{n}}{\eps\sqrt{\lambda}}\bigg)^n\\
	&= \frac{\omega_n|\Omega^{\pm\eps\sqrt{n}}|}{(2\pi)^n}\bigg(1 \pm \frac{\pi\sqrt{n}}{\eps\sqrt{\lambda}}\bigg)^n\lambda^{\frac{n}{2}}
\end{align*}

Then for any $\eps > 0$ we have $\Omega$ is Dirichlet-superspectral to $E^\eps_+$ and Dirichlet-subspectral to $E^\eps_-$.
\end{thm}
\begin{proof}
Denote by $N_0^\Omega$ the Dirichlet Laplace eigenvalue counting function.

Let $\eps > 0$ be fixed. Consider the lattice generated by $(\eps\mb{Z})^n$; the cells of this lattice are cubes of side length $\eps$. Denote by $\Omega_{out}$ the union of all cells in the lattice which intersect $\Omega$ and denote by $\Omega_{in}$ the union of all cells in the lattice which are contained in $\Omega$. Denote by $\#_{out}^\eps$ the number of cells comprising $\Omega_{out}$ and by $\#_{in}^\eps$ the number of cells comprising $\#_{in}$. Denote by $N_\nu^\eps$ the counting function of the Laplacian on the $n$-cube of side length $\eps$, where $\nu=0$ is the Dirichlet counting function and $\nu=1$ is the Neumann counting function.

By Dirichlet domain monotonicity, we have for arbitrary $\lambda \geq 0$:
\[ \#_{in}^\eps N_0^\eps(\lambda) \leq N_0^{\Omega_{in}}(\lambda) \leq N_0^\Omega(\lambda) \]
and
\[ N_0^\Omega(\lambda) \leq N_0^{\Omega_{out}}(\lambda) \leq N_1^{\Omega_{out}}(\lambda) \leq \#_{out}^\eps N_1^\eps(\lambda) \]

Now we have by the previous lemma
\[ \#^\eps_{in}\frac{\omega_n\eps^n}{(2\pi)^n}\bigg(1 - \frac{\pi\sqrt{n}}{\eps\lambda}\bigg)^n\lambda^{n/2} \leq N_0^\Omega(\lambda) \leq \#^\eps_{out}\frac{\omega_n\eps^n}{(2\pi)^n}\bigg(1 + \frac{\pi\sqrt{n}}{\eps\sqrt{\lambda}}\bigg)^n\lambda^{n/2}. \]

As the diagonal of each cell has length $\eps \sqrt{n}$ note that
\[ \Omega^{-\eps \sqrt{n}} \subset \Omega_{in} \subset \Omega \subset \Omega_{out} \subset \Omega^{\eps \sqrt{n}}. \]
The inequality $|\Omega^{-\eps \sqrt{n}}| \leq |\Omega| \leq |\Omega^{\eps \sqrt{n}}|$ yields the result.
\end{proof}

We remark that Weyl's law for the Dirichlet Laplacian on Euclidean domains follows as a corollary by setting $\epsilon = (\log lambda)^{-1}$ and noting that $\epsilon\sqrt{\lambda}\to\infty$ and $\epsilon\to 0$ as $\lambda\to\infty$.

Similar proof of a quantitative Weyl law for the Neumann Laplacian would follow from Conjecture \ref{neumann_monotonicity}.

Known proofs of asymptotic isospectrality to a Weyl polynomial with more terms involve a Tauberian theorem applied to the analysis of either the heat kernel or the wave propagation operator. The author is not aware of a Tauberian theorem that provides pointwise rather than asymptotic estimates; such a result would be of interest. (The closest the author has found is a paper of Brownell \cite{Brownell1955} providing log-Gaussian error bounds to the Laplace transform of the heat trace.)

\section{Subspectral Riemannian metrics}

We analyze subspectrality between two Riemannian metrics on a fixed manifold.

Let us fix a normal $n$-dimensional manifold $M$. If $M$ has boundary, we fix boundary conditions $\nu$. Suppose $g$ and $h$ are two Riemannian metrics on $M$.

To study Neumann eigenvalues we make the following definition.
\begin{defn}[Boundary-conformal]\label{bdry_conformal}
Say two Riemannian metrics on a compact smooth manifold $M$ are boundary-conformal provided the outward unit normal fields on $\partial M$ with respect to $g$ and $h$ differ by multiplication by a nowhere-zero smooth function.
\end{defn}
Boundary-conformal is an equivalence relation. Any two Riemannian metrics on a closed manifold are boundary conformal.

Use $g$ as a reference metric. In every tangent space, we have that $h$ is a symmetric bilinear operator which has a discrete, strictly positive spectrum. Denote by $0<\delta_1\leq\delta_2\leq\cdots\leq\delta_n$ the functions where $\delta_i(p)$ is the $i^{th}$ eigenvalue of $h$ with respect to $g$ in $T_pM$, for $i=1,\ldots,n$. An eigenvector $v_{ip}\in T_pM$ is a solution of the generalized eigenvalue problem $h(u,v_{ip}) = \delta_i g(u,v_{ip})$ for all $u\in T_pM$. Let $\delta_+ = \sup_M(\delta_n)$ and let $\delta_- = \inf_M(\delta_1)$. As an example, if $g$ and $h$ are conformally related, then $\delta_- = \delta_1(p) = \delta_n(p) = \delta_+$ for all $p$.

This gives bounds relating the volume forms and gradients of $g$ and $h$. Let $e_i$ be a local orthonormal frame field for $g$. Let $h_{ij} = h(e_i,e_j)$ be the coordinate matrix for $h$ with respect to $e_i$. In the $e_i$ coordinates, the volume form $dv_h = \sqrt{\det{h_{ij}}}e_1\wedge\cdots\wedge e_n = \sqrt{\det{h_{ij}}}dv_g$. As $\det(h_{ij}) = \delta_1\cdots\delta_n$, we have the bound $\delta_-^{n/2}dv_g \leq dv_h \leq \delta_+^{n/2} dv_g$. Likewise, in the $e_i$ coordinates, the differential $du$ of a smooth function $u$ can be written as a linear combination of the dual basis of $e_i$. The gradient with respect to $h$ in coordinates is $\nabla_h u = (h_{ij})^{-1}du$, and its pointwise norm can be written
\[|\nabla_h u|_h^2 = ((h_{ij}^{-1} du)^T h_{ij} (h_{ij}^{-1} du) = du^T h_{ij}^{-1} du\]
Recalling that in these orthonormal coordinates the inner product of $du$ with respect to $g$ is given by $du^T du$ we have the bounds
\[ \delta_+^{-1} |\nabla_g u|_g^2 \leq |\nabla_h u|_h^2 \leq \delta_-^{-1}|\nabla_g u|_g^2. \]

We now show that the domains of the quadratic forms of the $g$-Laplacian and $h$-Laplacian are canonically equal.
\begin{lemma}\label{metr_subsp_lemma1}
Let $M$ be compact with (possibly empty) boundary. If $\partial M$ is nonempty, set a boundary condition $\nu$. Let $g$ and $h$ be two smooth Riemannian metrics on $M$. If we set Neumann conditions, let $g$ and $h$ be boundary-conformal. Denote by $\Delta_g^\nu$ and $\Delta_h^\nu$ the $\nu$-Laplacians of $g$ and $h$, respectively, acting on $D^\nu$. Let $n_g$ and $n_h$ be the $L^2$ inner products acting on smooth functions with respect to $g,h$, resp. Then $\Delta_g^\nu$ and $\Delta_h^\nu$ are comparable operators in the sense of Definition \ref{def_comparable}.
\end{lemma}

\begin{proof}
Recall we have set $\mf{q}^\nu_g$ and $\mf{q}^\nu_h$ to be the energy forms corresponding to the $\nu$-Laplacians $\Delta_g^\nu$ and $\Delta_h^\nu$. Because $\nu$ is fixed we drop it for the duration of this proof.

We first show that the spaces $L^2(M,g)$ and $L^2(M,h)$ are equal. For $*\in\{g,h\}$, define the norm $\|u\|_*^2 = \int_M u^2\ dv_*$. Because $\delta_-^{n/2} dv_g\leq dv_h\leq \delta_+^{n/2}dv_g$, we have that for all $u\in C^\infty(M)$, 
\[\delta_-^{n/2}\|u\|_g^2\leq \|u\|_h^2\leq \delta_+^{n/2}\|u\|_g^2.\]
Any sequence in $C^\infty(M)$ is Cauchy with respect to one norm if and only if it is Cauchy with respect to the second. Therefore $L^2(M,g) = L^2(M,h)$. Let us denote the identified space $H$.

We have the norms on $D_\nu$ used in the construction of the Friedrichs extension: 
\[ \|u\|_{V,*}^2 = \|u\|_*^2 + \mf{q}_*(u,u) \]
for $*\in \{g,h\}$. Call the norms the $g$-norm and the $h$-norm.

Let $V_g$ (resp $V_h$) be the form domain of $\Delta_g$ (resp $\Delta_h$). To show that $V_g = V_h$, we show that a sequence $u_1,u_2,\cdots\in D_\nu$ is Cauchy with respect to the $g$-norm if and only if it is Cauchy with respect to the $h$-norm. We therefore show that the norms are equivalent. By the inequalities above, we have for any $u\in D_\nu$
\[ \int_M |\nabla_h u|_h^2\ dv_h + \int_M u^2\ dv_h \leq \delta_+^{n/2}\on{max}\{\delta_-^{-1},1\}\bigg( \int_M |\nabla_g u|_g^2\ dv_g + \int_M u^2\ dv_g \bigg) \]
and likewise
\[ \delta_-^{n/2}\on{min}\{\delta_+^{-1},1\}\bigg( \int_M |\nabla_g u|_g^2\ dv_g + \int_M u^2\ dv_g \bigg) \leq \int_M |\nabla_h u|_h^2\ dv_h + \int_M u^2\ dv_h \]
so that we have
\[ c\|u\|_{V,g}^2 \leq \|u\|_{V,h}^2 \leq C\|u\|_{V,g}^2 \]
for any $u\in D_\nu$, with $c = \delta_-^{n/2}\on{min}\{\delta_+^{-1},1\}$ and $C = \delta_+^{n/2}\on{max}\{\delta_-^{-1},1\}$, so the $g$-norm and the $h$-norm are equivalent.

Therefore $\Delta_g$ and $\Delta_h$ are comparable.
\end{proof}

We now apply Theorem \ref{metr_subsp} of Chapter 1.

\begin{prop}[Metric subspectrality]\label{metric_subspectrality}
Let $M$ be a normal manifold with specified boundary conditions $\nu$. Let $g$ and $h$ be two Riemannian metrics. If we have chosen Neumann conditions, let $g$ and $h$ be boundary-conformal. Denote by $\lambda_k(*)$ the $k^{th}$ eigenvalue of the Laplacian $\Delta_*^\nu$. Then
\[ \frac{1}{\delta_+}\bigg(\frac{\delta_-}{\delta_+}\bigg)^{n/2}  \leq \frac{\lambda_k(h)}{\lambda_k(g)} \leq \frac{1}{\delta_-}\bigg(\frac{\delta_+}{\delta_-}\bigg)^{n/2}. \]
\end{prop}
\begin{proof}
Recall the inequalities 
\[ \delta_+^{-1}\delta_-^{n/2} \mf{q}_g \leq \mf{q}_h \leq \delta_-^{-1}\delta_+^{n/2} \mf{q}_g \]
and
\[ \delta_-^{n/2} \|\cdot\|_g \leq \|\cdot\|_h \leq \delta_+^{n/2}\|\cdot\|_g \]
In the notation of Theorem \ref{metr_subsp}, we have 
\[c_H = \delta_-^{n/2}\] 
\[C_H = \delta_+^{n/2}\] 
\[c_V = \delta_+^{-1}\delta_-^{n/2}\]
 and 
 \[C_V = \delta_-^{-1}\delta_+^{n/2}.\] 
Then by Theorem \ref{metr_subsp}
\[ \frac{1}{\delta_+}\bigg(\frac{\delta_-}{\delta_+}\bigg)^{n/2}  = \frac{c_V}{C_H} \leq \frac{\lambda_k(h)}{\lambda_k(g)} \leq \frac{C_V}{c_H} = \frac{1}{\delta_-}\bigg(\frac{\delta_+}{\delta_-}\bigg)^{n/2} \]
\end{proof}

\begin{figure}
	\centering
	\includegraphics[width=0.7\textwidth]{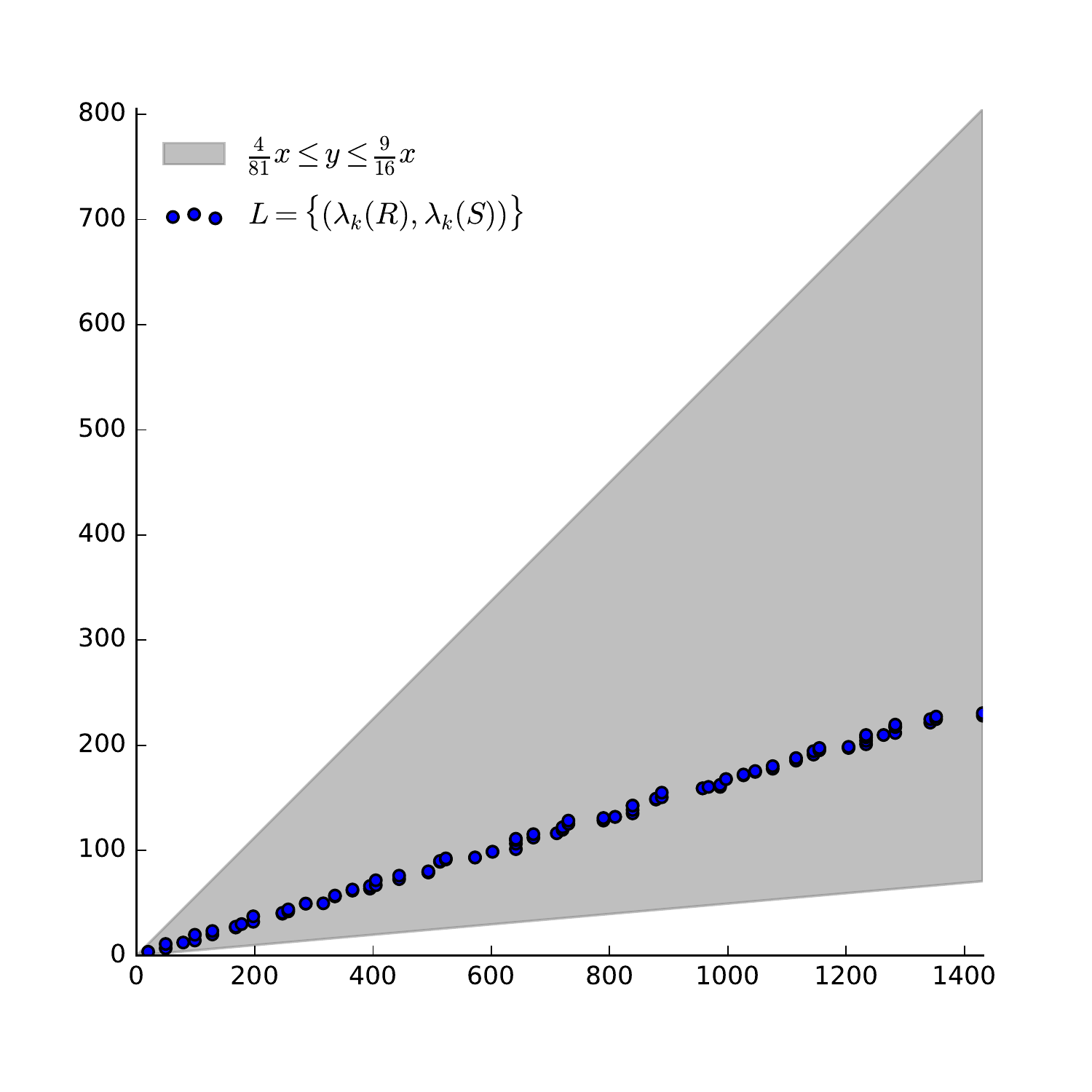}
	\caption[Illustration of bounds in Theorem \ref{metric_subspectrality}]{The set $L = \{(\lambda_k(R),\lambda_k(S))\ |\ k\in\mb{N}\}$ where $S = [0,1]\times[0,1]$ and $R = [0,2]\times[0,3]$. Computing the pullback of the metric on $R$ to $S$ by the linear map, we have $\delta_- = 4$, $\delta_+=9$. The theorem states $L$ lies within the shaded gray cone.}
	\label{metr_subsp_figure}
\end{figure}

We illustrate the proposition in Figure \ref{metr_subsp_figure}. 

We may now apply Theorem \ref{metr_subsp}. If $f:(M,g)\to (N,h)$ is a diffeomorphism of compact Riemannian manifolds (which extends to a diffeomorphism of their boundaries), we can obtain two metrics on $M$, $g$ and $f^*h$. Thence we have $\delta_+$ and $\delta_-$, defined as above.

\begin{cor}[Sufficient conditions for subspectrality]\label{subspectr_manifolds} Suppose $M$ and $N$ are manifolds with boundary and $f:M\to N$ is a diffeomorphism that extends to a diffeomorphism between $\partial M$ and $\partial N$. Fix the same boundary condition on $M$ and $N$. If considering Neumann boundary conditions, suppose that $f$ maps normal vectors to normal vectors. If $\delta_-\geq \delta_+^{1+2/n}$, then $M$ is subspectral to $N$. If $\delta_- \geq \delta_+^{\frac{n}{n+2}}$ then $M$ is superspectral to $N$.
\end{cor}
\begin{proof}
Recall that we have established, in the notation of Chapter 1, the following facts:
\[c_V = \delta_+^{-1}\delta_-^{n/2} \leq  \delta_-^{-1}\delta_+^{n/2} = C_V\]
and
\[c_H = \delta_-^{n/2} \leq \delta_+^{n/2} = C_H\]

Applying Lemma \ref{metr_subsp_lemma1} and Proposition \ref{metric_subspectrality} yields the desired result.
\end{proof}

We have the following illustration of this fact.
\begin{cor}
Suppose $M$ is a smooth closed manifold and $g$ is a Riemannian metric on $M$. If $h = e^fg$ is a conformal deformation of $g$ with $f$ a smooth function everywhere nonnegative, then $(M,h)$ is subspectral to $(M,g)$.
\end{cor}
\begin{proof}
As $\delta_+ \geq \delta_- = e^{\inf_M f} \geq 1$, we satisfy the conditions of the previous theorem.
\end{proof}

\section{Examples of subspectral rectangles}

Since the simplest domains of dimension greater than one with computable Laplace spectra are rectangles, we study subspectrality in rectangles. We construct pairs of rectangles which do not embed in each other where one is subspectral to the other.

\begin{prop}[Non-embedding subspectral pairs of rectangles]
Given a rectangle $R$, there exists a rectangle $R'$ such that $R'$ is Neumann subspectral to the $R$ but $R$ does not embed in $R'$.

Given a rectangle $R$, there exists a rectangle $R'$ such that $R$ is Dirichlet subspectral to $R$ but $R'$ does not embed in $R$.
\end{prop}

This is a corollary to the following lemmata.
\begin{lemma}\label{rect_event_subspect}
Suppose $R$ and $R'$ are rectangular prisms with codiagonals $d, d'$ and volumes $|R|$, $|R'|$ respectively. (Recall that the codiagonal of a rectangular prism is defined in Lemma \ref{quantitative_weyl_rect}.) Let
\[ \lambda_0 = (d')^2 \bigg(1 - \bigg[ \frac{|R|}{|R'|}\bigg]^{\frac{1}{n}}\bigg)^{-2} \tag{$*$} \]
and
\[ \lambda_1 = d^2\bigg(\bigg[\frac{|R'|}{|R|}\bigg]^{\frac{1}{n}} - 1\bigg)^{-2} \tag{$**$}\]
If $|R|< |R'|$ then:
\begin{itemize}
\item $R'$ is Dirichlet subspectral to $R$ beyond $\lambda_0$
\item $R'$ is Nemann subspectral to $R$ beyond $\lambda_1$
\end{itemize}
\end{lemma}

\begin{proof}
Notice that condition $(*)$ is equivalent to $|R| = (1 - d'/\sqrt{\lambda_0})^n|R'|$ and the condition $(**)$ is equivalent to $|R'| = (1 + d/\sqrt{\lambda_1})^n|R|$.

We apply Theorem \ref{quantitative_weyl_rect}, the quantitative Weyl law for rectangles. Recall for any rectangular prism $\Omega$ with codiagonal $\delta$ we have for all $\lambda > 0$ the Dirichlet eigenvalue counting function bounds
\[\frac{|\Omega|}{\pi^{n/2}\Gamma(n/2 + 1)}\bigg(1 - \frac{\delta}{\sqrt{\lambda}}\bigg)^n\lambda^{n/2} \leq N_0(\lambda) \leq  \frac{|\Omega|}{\pi^{n/2}\Gamma(n/2 + 1)}\lambda^{n/2} \tag{$D$}\]
and the Neumann eigenvalue counting function bounds
\[ \frac{|\Omega|}{\pi^{n/2}\Gamma(n/2 + 1)}\lambda^{n/2} \leq N_1(\lambda) \leq  \frac{|\Omega|}{\pi^{n/2}\Gamma(n/2 + 1)}\bigg(1 + \frac{\delta}{\sqrt{\lambda}}\bigg)^n\lambda^{n/2}. \tag{$N$}\]

Suppose $\lambda > \lambda_0$. Then applying $(D)$ we have
\begin{align*}
N_0(\lambda) &\leq \frac{\omega_n}{(2\pi)^n}|R|\lambda^{n/2} \\
	&= \frac{\omega_n}{(2\pi)^n}|R'|\bigg(1 - \frac{d'}{\sqrt{\lambda_0}}\bigg)^n\lambda^{n/2}\\
	&\leq \frac{\omega_n}{(2\pi)^n}|R'|\bigg(1 - \frac{d'}{\sqrt{\lambda}}\bigg)^n\lambda^{n/2}\\
	&\leq N_0'(\lambda)
\end{align*}
yielding that $R'$ is subspectral to $R$ beyond $\lambda_0$ as claimed.

Now suppose $\lambda > \lambda_1$. Then applying $(N)$ we have
\begin{align*}
N_0(\lambda) &\leq \frac{\omega_n}{(2\pi)^n}|R|\bigg( 1 + \frac{d}{\sqrt{\lambda}} \bigg)^n\lambda^{n/2}\\
	&\leq \frac{\omega_n}{(2\pi)^n}|R|\bigg( 1 + \frac{d}{\sqrt{\lambda_1}} \bigg)^n\lambda^{n/2} \\
	&= \frac{\omega_n}{(2\pi)^n}|R'|\lambda^{n/2} \\
	&\leq N_0'(\lambda)
\end{align*}
\end{proof}
Note that $(*)$ and $(**)$ are derived from $(D)$ and $(N)$, respectively.

\begin{lemma}
Given a rectangle $R = [0,L]\times[0,W]$, there exist $r>0$ and $1 > \eps > 0$ such that the map $(x,y)\mapsto (\eps x, ry)$ carries $R$ to a rectangle $R'$ which is Dirichlet subspectral to $R$.	
\end{lemma}
\begin{proof}
Let $d$ denote the codiagonal of $R$. Recall that $d^2$ is equal to the smallest Dirichlet eigenvalue of $R$. It will suffice to find some $R'$ subspectral to $R$ beyond $d^2$, because for $\lambda < d^2$ the counting function $N_R(\lambda) = 0$.

Choose $\xi\in(0,d^2)$. Choose $\epsilon > \sqrt{\frac{1}{1 + (L/W)^2 - \xi L^2/\pi^2}}$. Then $\frac{\pi^2}{\eps^2 L^2} < d^2 - \xi$. Let $r > d^2/(\eps \xi^2)$. Denote by $R'$ the rectangle $[0,\eps L]\times [0, rW]$.

By Lemma \ref{rect_event_subspect} we have that $R'$ is subspectral to $R$ beyond 
\[ \lambda_0 = (d')^2 \bigg(1 - \bigg[ \frac{|R|}{|R'|}\bigg]^{\frac{1}{2}}\bigg)^{-2}.\]
By choice of $\eps$ we have that $(d')^2 = \pi^2/(\eps L)^2 + \pi^2/(rW)^2 < d^2 - \xi$. By choice of $r$ we have the following implications:
\begin{align*}
r > \frac{d^2}{\eps\xi^2} &\Rightarrow \frac{1}{\eps r} < \frac{\xi^2}{d^2} \\
	&\Rightarrow \sqrt{\frac{1}{\eps r}} < \frac{\xi}{d} \\
	&\Rightarrow 1 - \sqrt{\frac{1}{\eps r}} > 1 - \frac{\xi}{d} \\
	&\Rightarrow \frac{1}{1 - \sqrt{\frac{1}{\eps r}}} < \frac{1}{1 - \frac{\xi}{d}}
\end{align*}

Thus
\begin{align*}
\lambda_0 &= (d')^2 \bigg(1 - \bigg[ \frac{|R|}{|R'|}\bigg]^{\frac{1}{2}}\bigg)^{-2} \\
	&< (d-\xi)^2\bigg(1 - \sqrt{\frac{1}{\eps r}}\bigg)^{-2} \\
	&= \bigg( \frac{d - \xi}{1 - \sqrt{\frac{1}{\eps r}}}\bigg)^2 \\
	&< \bigg( \frac{d - \xi}{1 - \frac{\xi}{d}} \bigg)^2 \\
	&= d^2
\end{align*}
so $R'$ is indeed subspectral to $R$.
\end{proof}

\begin{lemma}
Given a rectangle $R = [0,L]\times [0,W]$, for any $\eps > 0$ there exists an $r_0 > 0$ such that for all $r > r_0$ the map $(x,y)\mapsto (\eps x, ry)$ carries $R$ to a rectangle $R'$ which is Neumann subspectral to $R$.
\end{lemma}
\begin{proof}
Let $\eps$ be given. For any $r > 1/\eps$, let $R' = [0,\eps L]\times [0, rW]$. By Lemma \ref{rect_event_subspect}, we have $R'$ is subspectral to $R$ beyond 
\[\bigg( \frac{d}{\sqrt{\eps r} - 1} \bigg).\]
Denote by $\xi$ the first nonzero Neumann eigenvalues of $R$. Choose $r_0 = \frac{1}{\eps}( + d\sqrt{\xi})^2$. As this quantity is greater than $1/\eps$, the condition of Lemma \ref{rect_event_subspect} are satisfied for any $r > r_0$ and
\[ \bigg( \frac{d}{\sqrt{\eps r} - 1} \bigg) < \xi. \]

Because $R$ is connected, we have $N_R(\lambda) = 1$ for all positive $\lambda < \xi$ and because $N_{R'} \geq 1$ we have that $R'$ is subspectral to $R$.
\end{proof}

\section{Converses to domain monotonicity}

We ask to what extent subspectrality between domains with the same boundary conditions implies that one domain can be embedded in the other. We therefore investigate the relationship between subspectrality and embedding and construct examples where Dirichlet subspectrality holds but domain containment does not.

\begin{prop}[Counterexample to Dirichlet domain monotonicity converse]
For any regular Euclidean domain $\Omega$ there exists a domain $\Omega'$ such that $\Omega$ is Dirichlet subspectral to $\Omega'$ but $\Omega'$ does not embed in $\Omega$.
\end{prop}
\begin{proof}
Suppose $\Omega\subset\mb{R}^n$. We produce the domain $\Omega'$ as a rectangular prism $R$ isometric to $[0,s]^{n-1}\times[0,L]$ for some large $L$ and small $s > 0$ to be chosen.

We will choose $L$ and $s$ such that $|R| = Ls^{n-1} < |\Omega|$ so we are guaranteed that $\Omega$ is subspectral to $R$ beyond some $\lambda_0$; by the quantitative Weyl laws for rectangles and domains, Proposition \ref{quantitative_weyl_rect} and Theorem \ref{quantitative_weyl}, we will choose $L$ and $s$ such that $\lambda_0$ is the infimum of the set of $\lambda$ such that
\[ N_0^R(\lambda) \leq C|R|\lambda^{n/2} \leq C\bigg|\Omega^{-\eps\sqrt{n}}\bigg|\bigg(1 - \frac{\pi}{\eps}\sqrt{\frac{n}{\lambda}}\bigg)^n\lambda^{n/2} \leq N_0^\Omega(\lambda) \]
where $C = (\pi^{n/2}\Gamma(n/2 + 1))^{-1}$.

Canceling $C$ and $\lambda^{n/2}$ and rearranging terms, we see this holds exactly when
\[ \frac{\pi\sqrt{n}}{\eps} \frac{|\Omega^{-\eps\sqrt{n}}|^{1/n}}{|\Omega^{-\eps\sqrt{n}}|^{1/n} - |R|^{1/n}} \leq \sqrt{\lambda}. \]
As the square root function is monotone we have 
\[\lambda_0 = \frac{\pi^2n}{\eps^2} \bigg(\frac{|\Omega^{-\eps\sqrt{n}}|^{1/n}}{|\Omega^{-\eps\sqrt{n}}|^{1/n} - |R|^{1/n}} \bigg)^2\]

Note that $\lambda_0$ depends on $\eps$, $L$, and $s$. Now let $\eps$ be small enough so that $|\Omega^{-\eps\sqrt{n}}| > 0$ and let $L > \diam\Omega$. Having chosen $\eps$ and $L$, note that $\lambda_0$ is a continuous function of $s$.

Notice that given choices of $L$ and $s$ we have the first Dirichlet eigenvalue of $R$ is $(n-1)\frac{\pi^2}{s^2} + \frac{\pi^2}{L^2}$. Given choices of $L$ and $\eps$ as in the previous paragraph, if there exists an $s$ with $(n-1)\frac{\pi^2}{s^2} + \frac{\pi^2}{L^2} > \lambda_0(s)$ then $\Omega$ is Dirichlet-subspectral to $R(L,s)$. But in the inequality
\[ (n-1)\frac{\pi^2}{s^2} + \frac{\pi^2}{L^2} > \frac{\pi^2n}{\eps^2} \bigg(\frac{|\Omega^{-\eps\sqrt{n}}|^{1/n}}{|\Omega^{-\eps\sqrt{n}}|^{1/n} - |R|^{1/n}}\bigg)^2 \]
as $s\to 0$ the right side tends to $\frac{\pi^2n}{\eps^2}$ while the left side tends to infinity. This establishes the existence of such an $s$, proving the theorem.
\end{proof}

\section{Equal-area subspectral domains}

Notice that the use of Theorem \ref{quantitative_weyl} and Proposition \ref{subspectr_manifolds} requires the two domains to have different areas. As reported in Ivrii's survey \cite{Ivrii2016}, Weyl \cite{Weyl1913} conjectured that any regular domain $\Omega$ has
\[N_\nu(x) = \frac{|\Omega|}{4\pi}x + (-1)^\nu \frac{|\partial \Omega|}{8\sqrt{\pi}}\sqrt{x} + o(1) \]
where we use $\nu$ to denote Dirichlet or Neumann boundary conditions.

Suppose we have two domains $\Omega_1$, $\Omega_2$ with equal area but $|\partial\Omega_1| > |\partial\Omega_2|$. Then we might expect that $\Omega_1$ is Dirichlet subspectral and Neumann superspectral to $\Omega_2$.

As we a priori do not have uniform control of the error terms in the counting function, we cannot use a quantitative Weyl law to explore this intuition. However, we can conduct numerical exploration. A computer count of eigenvalues for equal-area rectangles in $\mb{R}^2$ results in several observations.

If the rectangles' perimeters are close, then one rectangle is eventually subspectral to the other, with violations to subspectrality appearing only for small eigenvalues. If one rectangle's perimeter is much greater than another's, then subspectrality appears to hold.

\begin{figure}[p]
	\centering
	\includegraphics[scale=0.8]{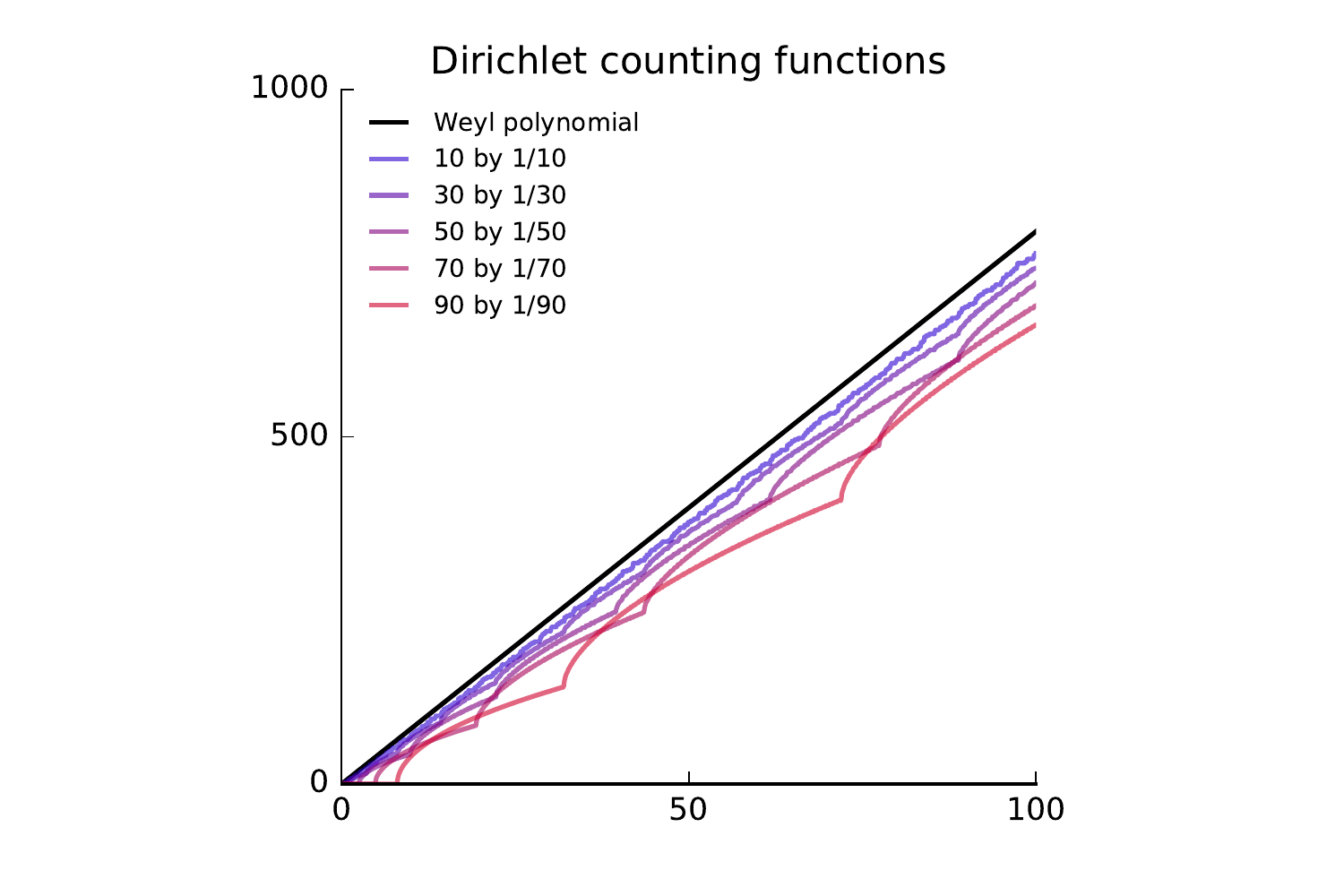}
	\caption[Weyl functions and Dirichlet counting functions for five rectangles]{Dirichlet eigenvalue counting functions and Weyl functions for five rectangles of area $100$}
	\label{multiple_rects_dirichlet}
\end{figure}
\begin{figure}[p]
	\centering
	\includegraphics[scale=0.8]{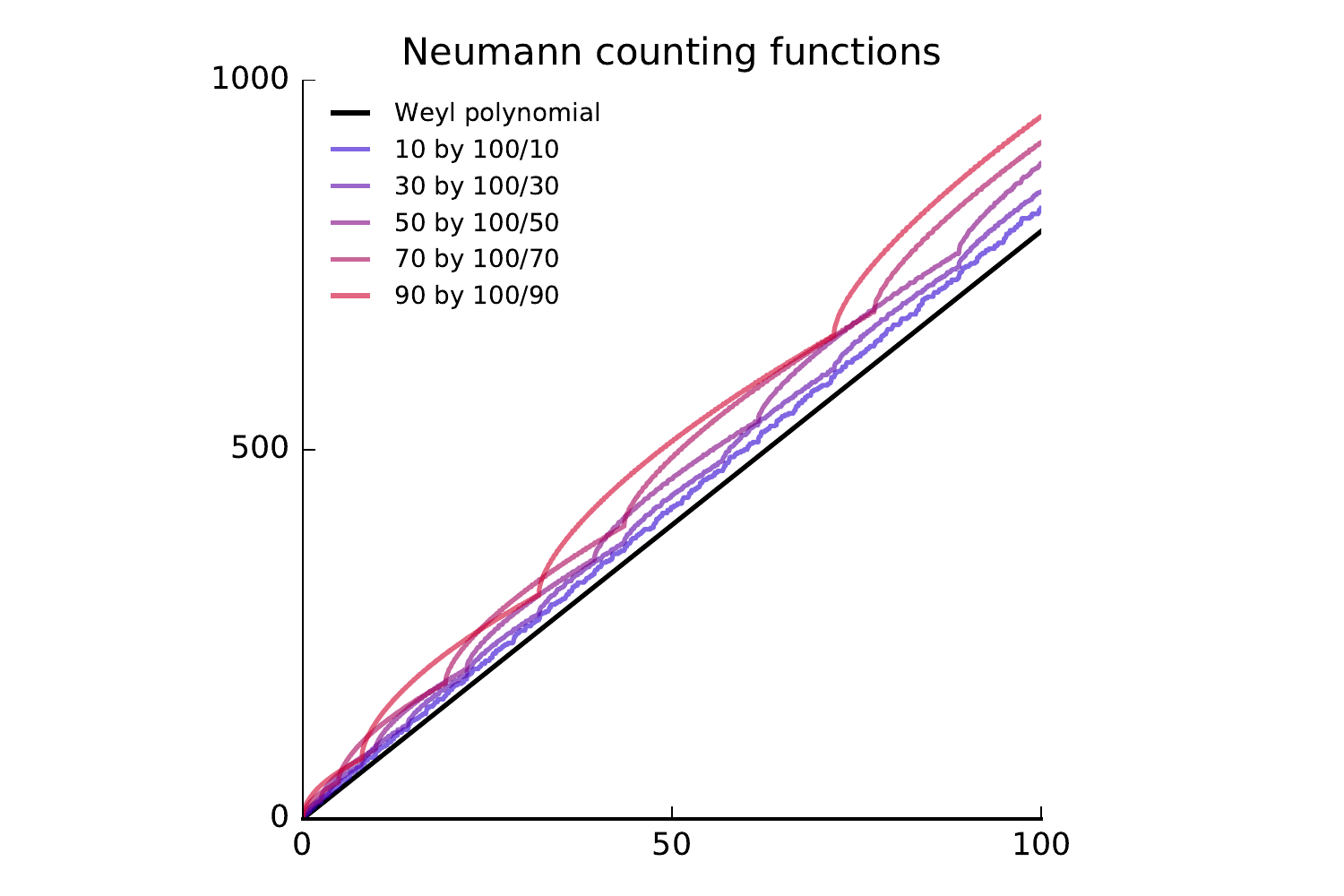}
	\caption[Weyl functions and Neumann counting functions for five rectangles]{Neumann eigenvalue counting functions and Weyl polynomial for five rectangles of area $100$}
	\label{multiple_rects_neumann}
\end{figure}

Plotting a variety of rectangle eigenvalue counting functions in Figures \ref{multiple_rects_dirichlet} and \ref{multiple_rects_neumann} clarifies this intuition. As the rectangle's perimeter increases, the counting function moves further from the one-term Weyl polynomial. We might interpret this to mean that the error term of the two-term Weyl law is relatively small even for small eigenvalues.


This numerical exploration leads to the following conjecture:

\begin{conj}[Equal-area subspectral rectangles]
For any rectangle $R$, there exists an $s>0$ such that for all $0<\epsilon<s$ the rectangle $R_\epsilon$ obtained from $R$ by applying the map $(x,y)\mapsto (\eps x, y/\eps)$ is Dirichlet-superspectral to and Neumann-subspectral to $R$.
\end{conj}

\begin{figure}[p]
	\centering
	\includegraphics[scale=0.5]{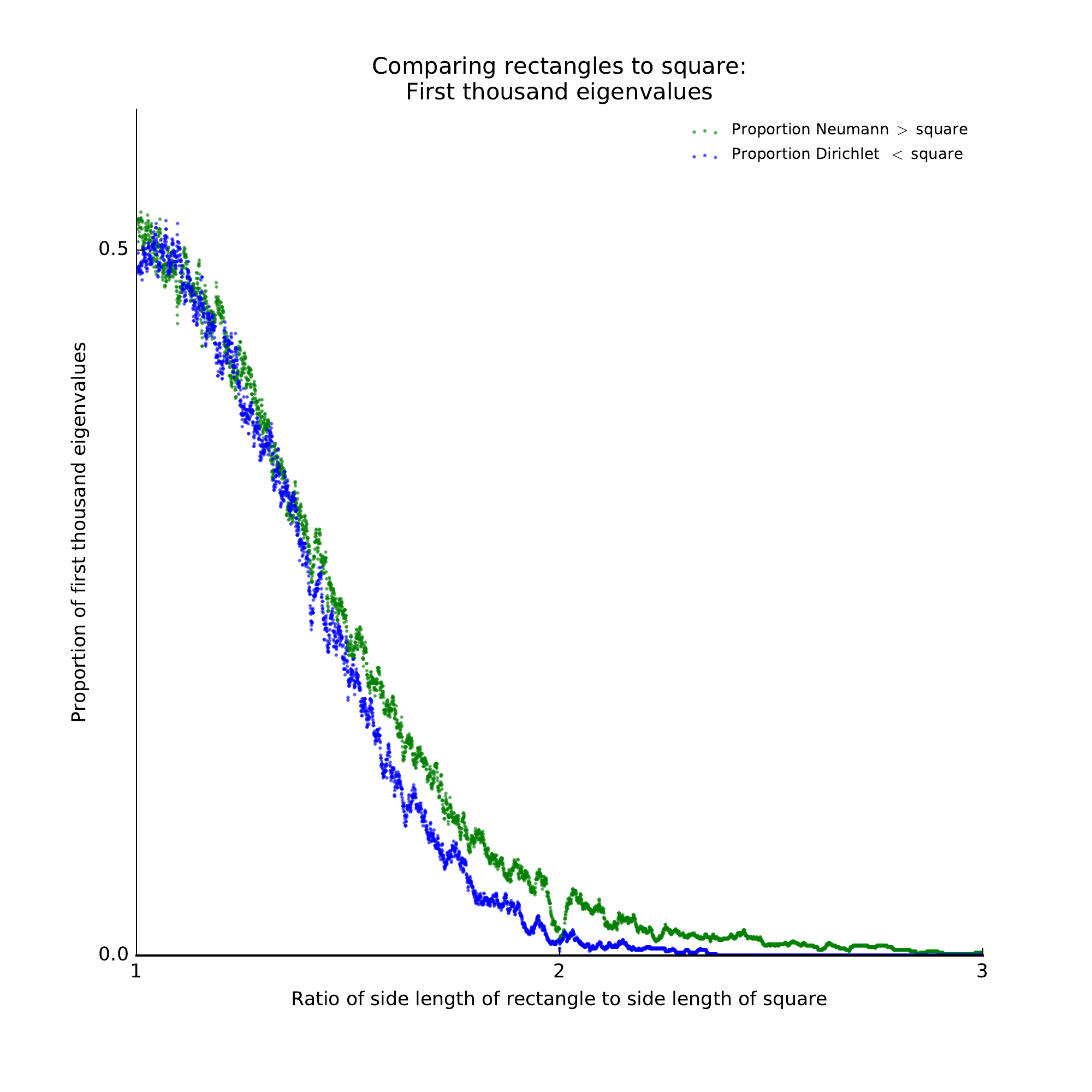}
	\caption[Proportion of low-frequency eigenvalues violating expected subspectrality as a function of ratio of side lengths]{A rectangle is expected to be Dirichlet-superspectral and Neumann-subspectral to a square with the same area. The horizontal axis in this figure measures the ratio of the rectangle's longer side to the side length of a square of equal area. The two series plotted here are the proportion of Dirichlet (resp. Neumann) eigenvalues which violate the expected super-(resp. sub-)spectrality.}
	\label{rect_square_1000}
\end{figure}

We now briefly explore measuring how quickly eventual subspectrality takes hold for equal-area rectangles where one is not subspectral to the other. One measurement is the proportion of low frequency eigenvalues which violate the subspectrality indicated by the two-term Weyl law. We numerically examine this measure by computing the first thousand Neumann and Dirichlet eigenvalues for a variety of rectangles of area $100$ and varying perimeter, and plotting the proportion of eigenvalues violating expected sub- and superspectrality to the square as a function of side length. This is illustrated in Figure \ref{rect_square_1000}.

A second measurement is measuring the proportion of the first $n$ eigenvalues whcih violate expected subspectrality, as a function of $n$. We call this the subspectral mean and define it precisely:
\begin{defn}
Let $\Omega_1,\Omega_2$ be two domains in $\mb{R}^n$. Fix the same boundary conditions for the Laplacian on both. Denote by $\lambda_k(\Omega_i)$, $k\geq \nu$, the Laplace spectrum of $\Omega_i$. For each $k\in\mb{N}$ define
\[ 1_{\lambda_k(\Omega_1) > \lambda_j(\Omega_2)} = \begin{cases} 1, & \lambda_k(\Omega_1) > \lambda_k(\Omega_2) \\
																 0, & \lambda_k(\Omega_1) \leq \lambda_k(\Omega_2)
																 \end{cases}\]
Define the subspectral mean of $\Omega_1$ and $\Omega_2$ at $n$ to be the quantity
\[ \frac{1}{n}\sum_{k=\nu}^n 1_{\lambda_k(\Omega_1) > \lambda_k(\Omega_2)} \]
\end{defn}

\begin{figure}[p]
	\centering
	\includegraphics[scale=0.7]{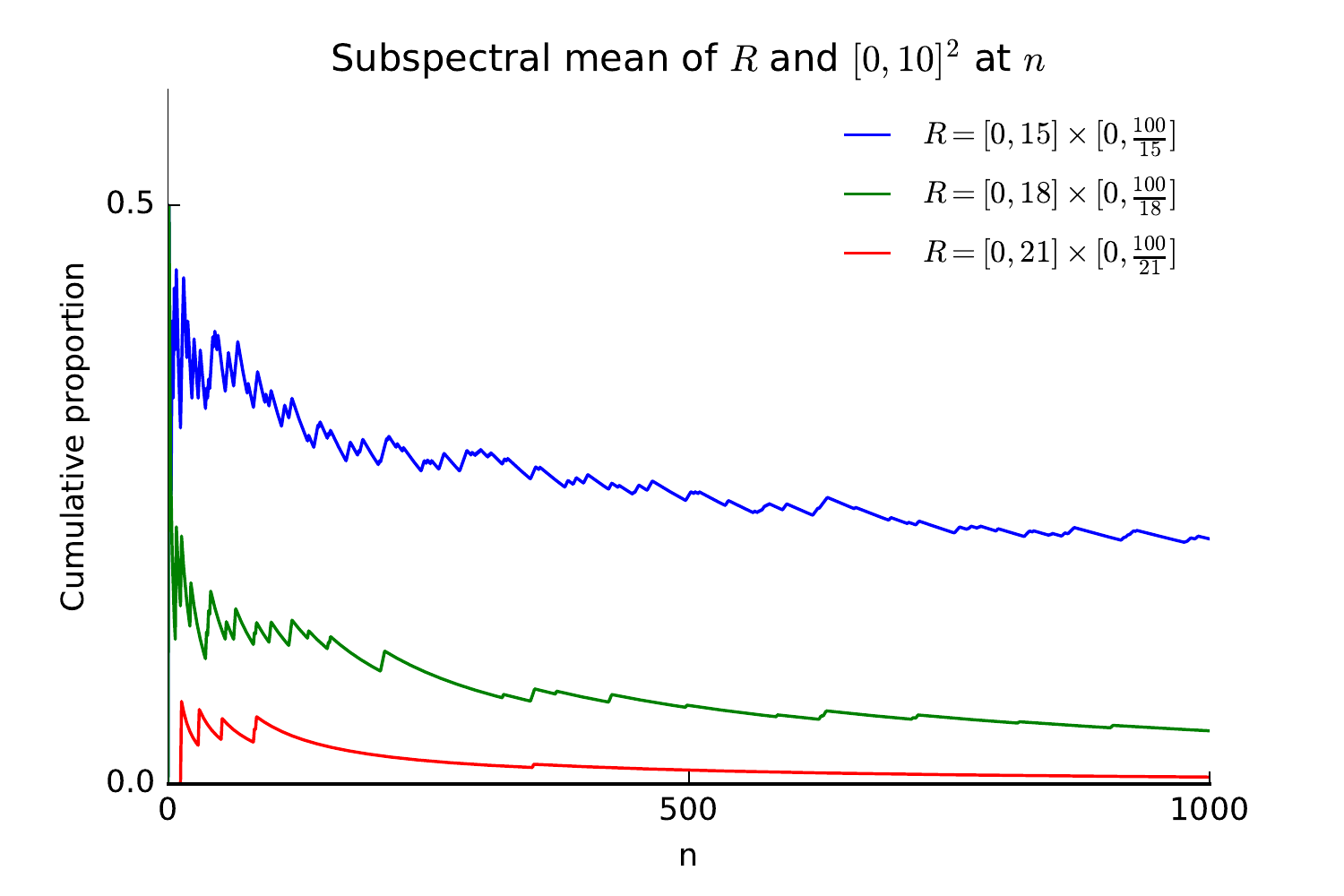}
	\caption[Dirichlet subspectral mean of three rectangles to the square of the same area]{The subspectral mean of $\Omega_1$ and $\Omega_2$ at $n$ is the proportion of the first $n$ eigenvalues of $\Omega_1$ that are less than the same-index eigenvalues of $\Omega_2$. Here we compute the Dirichlet subspectral mean of $R$ to the square of the same area, for the three values of $R$ indicated in the figure.}
	\label{subspec_mean}
\end{figure}

The subspectral mean of $\Omega_1$ and $\Omega_2$ at $n$ measures the proportion of eigenvalues of index no greater than $n$ which violate the assertion that $\Omega_1$ is subspectral to $\Omega_2$. If $\Omega_1$ is subspectral to $\Omega_2$ then the subspectral mean of $\Omega_1$ and $\Omega_2$ is equal to zero for each $n$. If $|\Omega_1|>|\Omega_2|$, then the limit of the subspectral mean of $\Omega_1$ and $\Omega_2$ as $n\to\infty$ is equal to zero. We illustrate this function in Figure \ref{subspec_mean}.

Some natural questions: How does the asymptotic behavior of the subspectral mean of $\Omega_1$ and $\Omega_2$ as $n\to\infty$ relate to the geometry of $\Omega_1$ and $\Omega_2$? Does there exist a pair of equal-area domains $\Omega_1,\Omega_2$ such that the subspectral mean does not converge to $0$ or $1$?



\chapter{Necessary conditions for subspectrality}%
We study implications of subspectrality for closed manifolds.

\section{Weyl law}

The simplest implications of subspectrality are derived from Weyl's law. The sharpest known form of Weyl's law is due to Ivrii and Melrose independently \cite{Ivrii1980}:
\begin{thm}\label{weyl_ivrii}[Weyl's law (Ivrii, Melrose)]
Suppose $M$ is a manifold with boundary such that the set of periodic points for the geodesic billiard flow has measure zero. Then 
\[ N_{\Delta_\nu} = \frac{\omega_n}{(2\pi)^{n}}|M|x^{n/2} + (-1)^{1-\nu} \frac{|\partial M|\omega_{n-1}}{4(2\pi)^{n-1}}x^{(n-1)/2} + o(x^{(n-2)/2} \]
\end{thm}
For other results in this direction, see a recent survey also by Ivrii \cite{Ivrii2016}.

We have the following by combining Lemma \ref{subspectr_sum} with Theorem \ref{weyl_ivrii}.
\begin{cor}
Suppose $M$ and $N$ are, respectively, $m$- and $n$-dimensional manifolds with smooth boundary. Suppose $M$ is subspectral to $N$ and both manifolds have periodic geodesic billiards. Then $m\geq n$. If $m=n$, we have $|M|\geq |N|$. If $|M| = |N|$, we have $(-1)^{1-\nu}|\partial M| \geq (-1)^{1-\nu}|\partial N|$.
\end{cor}
\begin{proof}
By Theorem \ref{weyl_ivrii}, both $M$ and $N$ are asymptotically isospectral to their two-term Weyl polynomials. We apply Theorem \ref{subspectr_sum} with $T = \Delta^\nu(M)$, $T' = \Delta^\nu(N)$, and $p$ and $p'$ the two-term Weyl polynomials for $M$ and $N$, respectively.
\end{proof}

\section{Heat trace}


We define the heat trace of $M$ to be the sum
\[ Z_M(t) = \sum e^{-\lambda_k(M)t}. \]
This is related to the heat kernel of the Laplacian on $M$; see Chavel \cite{Chavel1984} VII.3(31-2) and VI.1.


By work of Minakshisundaram and Pliejel (\cite{Minakshisundaram1949a}, \cite{Minakshisundaram1949}, \cite{Minakshisundaram1953}) the heat trace obeys the asymptotic expression
\[ \lim_{t\to 0} \frac{Z(t)}{(4\pi t)^{-n/2}\sum_{j=0}^J t^ja_j} = 1 \]
where $J > \dim(M)$. For derivation, see Chavel \cite{Chavel1984}, VI; further study for domains can be found in Kac \cite{Kac1966} and for polygons in van den Berg-Srisatkunarajah \cite{VandenBerg1988}.

\begin{prop}
Suppose $M$ and $N$ are closed manifolds of dimensions $m$ and $n$, respectively. Suppose that $M$ is subspectral to $N$. Denote by $a_j(M)$ (resp. $a_j(N)$) the $j^{th}$ Minakshisundaram-Pleijel coefficient of $M$ (resp. $N$). If $m=n$ and for all $j<k$, $a_j(M) = a_j(N)$. Then $a_k(M) > a_k(N)$.
\end{prop}
\begin{proof}
This follows immediately by applying Lemma \ref{subspectr_lapl} to the asymptotic expansion of the heat trace.
\end{proof}

We apply this to a result of McKean and Singer \cite{Mckean1967} to produce two corollaries.
\begin{cor}[Global curvatura integra comparisons]
Let $M$ and $N$ be two closed manifolds with $M$ subspectral to $N$. Let $K$ denote Gaussian curvature. Then if $|M| = |N|$, we have 
\[ \int_M K \geq \int_N K. \]
\end{cor}
\begin{proof}
McKean and Singer \cite{Mckean1967} prove that for a closed manifold $M$ the second coefficient $a_2(M)$ of the heat trace is proportional to the integral over the manifold of the Gaussian curvature with a proportionality constant that depends only on dimension.
\end{proof}

\begin{cor}
Suppose $M$ and $N$ are two closed surfaces with $M$ subspectral to $N$. If $|M| = |N|$, then the genus of $M$ is at most the genus of $N$.
\end{cor}
\begin{proof}
Recall that the genus $g$ of a surface is related to its Euler characteristic $\chi$ by $\chi = 2 - 2g$. By the Gauss-Bonnet theorem (see Hubbard \cite{Hubbard2006} 2.4.15 or Taylor \cite{Taylor2011b} Appendix C 5.31), the Euler characteristic of a surface is proportional to the integral of its Gaussian curvature, so we have $\chi(M) \geq \chi(N)$, hence $g(N) \geq g(M)$.
\end{proof}

As a partial converse, if we have an inequality of the heat trace between two manifolds, then we may conclude an inequality for the first unequal eigenvalues of their respective Laplacians. 
\begin{prop}
Suppose $M$ and $N$ are regular manifolds with boundary of the same dimension. Apply the same boundary conditions to their Laplacians. If $Z_M(t) \geq Z_N(t)$ for all $t$ and $\lambda_j(M) = \lambda_j(N)$ for all $j < k$, then $\lambda_k(M) \leq \lambda_k(N)$. 
\end{prop}
\begin{proof}
By assumption,
\[(4\pi t)^{-m/2}Z_M(t) \geq (4\pi t)^{-n/2}Z_N(t)\]
that is
\[ (4\pi t)^{-m/2}\sum_{j=1}^\infty e^{-\lambda_j(M)t} \geq (4\pi t)^{-n/2}\sum_{j=1}^\infty e^{-\lambda_j(N)t} \]
Let $k = \min\{k_1,k_2\}$. By the assumption on dimension, cancel $(4\pi t)^{n/2} = (4\pi t)^{m/2}$.

We now apply Corollary \ref{heat_ineq_cor} with 
\[F(t) = \sum_{j=1}^\infty e^{-\lambda_j(M)t}\]
and
\[G(t) = \sum_{j=1}^\infty e^{-\lambda_j(N)t}\]
to conclude that $\lambda_k(M) \leq \lambda_k(N)$.
\end{proof}

Note that we must have the same dimension, for otherwise we could not cancel the power of $t$. 

\section{Poisson summation formula}

Let us restrict our attention to $n$-dimensional flat tori, that is, manifolds $T_\Lambda = \mb{R}^n/\Lambda$ where $\Lambda$ is a lattice acting cocompactly by translations on $\mb{R}^n$. Let $\Lambda^*$ be the dual lattice. The Laplace spectrum of $T_\lambda$ is equal to the set of squared lengths of elements of $\Lambda^*$, so that $N(\lambda) = \#(\Lambda^*\cap B(0,\sqrt{\lambda}))$.

The Poisson summation formula is
\[ \sum_{k\in\Lambda^*} f(|k|) = |\Lambda| \sum_{\ell\in\Lambda} \widetilde{f}(|\ell|) \]
where $\widetilde{\cdot}$ is the integral transform
\[ \widetilde{f}(s) = \frac{1}{(2\pi)^n} \int_0^\infty \int_{S^{n-1}} e^{isr\eta} f(r)r^{n-1}d\eta dr, \]
the function $f$ is a test function on $\mb{R}$, and $|\Lambda|$ is the volume of a Dirichlet domain of $\Lambda$.

In particular, for $f = e^{-w^2t}$ a Gaussian depending on the parameter $t$, we have $\widetilde{f}(u) = \frac{|\Lambda|}{(4\pi t)^{n/2}}e^{-u^2/4t}$, and the Jacobi identity: for all $t>0$,
\[ \sum_{k\in\Lambda^*} e^{-|k|^2t} = \frac{|\Lambda|}{(4\pi t)^{n/2}} \sum_{\ell\in\Lambda} e^{-|\ell|^2/4t}. \]

We now prove:
\begin{prop}
Let $T_1$ and $T_2$ be two $n$-dimensional flat tori. If $T_2$ is subspectral to $T_1$ and $|T_1| = |T_2|$, then the systole of $T_1$ is at most the systole of $T_2$.
\end{prop}

If we drop the restriction that the tori have the same volume, the statement is false: consider the square torus $T_2 = \mb{R}^2/\mb{Z}^2$, and the smaller torus $T_1 = \mb{R}^2/(\epsilon\mb{Z}^2)$, for $\epsilon<1$. Then the systole of $T_1$ is shorter than the systole of $T_2$, but $T_2$ is certainly subspectral to $T_1$, as $T_1$ is an $\epsilon$ dilation of $T_2$.

\begin{proof}
Suppose we have two lattices, $\Lambda_1$ and $\Lambda_2$, defining two tori $T_i = \mb{R}^n/\Lambda_i$, $i=1,2$, such that for each $k>0$, $\lambda_k(T_1) > \lambda_k(T_2)$. The left hand sides of the Jacobi identities for the two lattices $\Lambda_1$ and $\Lambda_2$ combined with subspectrality give the following inequality:
\[ \frac{|T_1|}{(4\pi t)^{n/2}}\sum_{\ell\in\Lambda_1} e^{-\ell^2/4t} = \sum_{k\in\Lambda_1^*} e^{-|k|^2t} < \sum_{k\in\Lambda_2^*}e^{-|k|^2t} = \frac{|T_2|}{(4\pi t)^{n/2}}\sum_{\ell\in\Lambda_2} e^{-\ell^2/4t} \]

By the hypothesis that $T_1$ and $T_2$ are isovolumetric, we may subtract the constant terms, yielding:
\[ \sum_{\ell\in\Lambda_1-\{0\}} e^{-\ell^2/4t} < \sum_{\ell\in\Lambda_2-\{0\}} e^{-\ell^2/4t}. \]

Denote by $\ell_1$ and $\ell_2$ the systoles of $T_1$ and $T_2$, respectively. Substituting $\tau$ for $1/t$ we may apply Lemma \ref{heat_ineq_lemma} with $F(\tau) = \sum_{\ell\in \Lambda_1} e^{-(\ell^2/4)\tau}$ and $G(\tau) = \sum_{\ell\in \Lambda_2} e^{-(\ell^2/4)\tau}$ to conclude that $\ell_1^2/4 \leq \ell_2^2/4$; the desired result follows from multiplying both sides by $4$ and taking positive square roots.
\end{proof}

It is not clear that there exists a pair of tori to which this theorem applies. 



\section{Hyperbolic surfaces}

We consider closed surfaces of genus $g$ with constant sectional curvature $-1$. By Gauss-Bonnet we must have $g>1$. Such surfaces can be realized as the quotient of the Poincare disk, the unit disk in $\mb{R}^2$ with distance element $ds^2 = \frac{dx^2+dy^2}{(1-x^2-y^2)^2}$, by the isometric cocompact action of a discrete group of isometries. Let $M$ be such a surface. As a closed Riemannian manifold, $M$ has Laplacian acting on smooth functions and the closed eigenvalue problem has a discrete set of eigenvalues $0=\lambda_0(M)<\lambda_1(M)\leq\cdots\to \infty$. For an introduction, see Buser \cite{Buser1992}, Ch 7.

We concern ourselves with the relationship between the eigenvalue spectrum and the geodesic length spectrum of $M$. For discussion of the following concepts see Buser \cite{Buser1992} Ch 1 and 9. Each free homotopy class of noncontractible curves has a unique geodesic representative (Buser \cite{Buser1992} 1.6.6). Therefore to each free homotopy class we may associate the length of its geodesic representative.

\begin{defn}[Length spectrum; simple length spectrum; systole]
The sequence of lengths of closed geodesics on $M$, ordered by magnitude, is called the length spectrum of $M$. A primitive geodesic is one which does not retraverse itself. The sequence of lengths of primitive geodesics, ordered by magnitude, is the primitive length spectrum of $M$ and is denoted $P(M)$. The length of the shortest geodesic on $M$ is called the systole of $M$.
\end{defn}
These can be found in Definitions 9.2.5 and 9.2.8 and Lemma 9.2.6 of Buser \cite{Buser1992}.

The following analogue of the Jacobi identity, due to McKean and derivable from the Selberg trace formula, holds:
\begin{align*}
\sum e^{-\lambda_k t} &= \frac{|M|e^{-t/4}}{(4\pi t)^{3/2}}\int_0^\infty \frac{re^{-r^2/4t}}{\sinh{r/2}}\ dr \\
	&+ \frac{e^{-t/4}}{2(4\pi t)^{1/2}}\sum_{n=1}^\infty \sum_{\gamma\in P(M)} \frac{|\gamma|}{\sinh(|\gamma^n|/2)} e^{-|\gamma^n|^2/4t}.
\end{align*}
This is 9.2.11 in Buser \cite{Buser1992}.

We express this more succinctly. Set $C_g =  \frac{|M|e^{-t/4}}{(4\pi t)^{3/2}}\int_0^\infty \frac{re^{-r^2/4t}}{\sinh{r/2}}\ dr$. Note that the left side is the trace of the heat kernel, denoted $Z_M(t)$. Let us order the lengths of the primitive geodesics and denote them $l_1, l_2, \ldots$. Then we may rewrite the formula:
\[ Z_M(t) = C_g + \frac{e^{-t/4}}{2(4\pi t)^{1/2}}\sum_{n,k=1}^\infty \frac{l_k}{\sinh(nl_k/2)}e^{-(nl_k)^2/4t}. \]

\subsection{Laplace subspectrality implies systolic inequality}

We prove the following proposition:
\begin{prop}Let $M_1$ and $M_2$ be two hyperbolic surfaces of the same genus. If $M_2$ is subspectral to $M_1$, the systole of $M_2$ is no greater than the systole of $M_1$.\end{prop}

\begin{proof}
Since $M_1$ and $M_2$ have the same genus, by Gauss-Bonnet their areas are equal. Let $\ell_i$ be the systole of $M_i$. Enumerate the primitive length spectrum of $M_1$ as $l_k$ and the primitive length spectrum of $M_2$ as $m_k$. The proof then proceeds as in the case of a torus. By Lemma \ref{subspectr_lapl}, because $M_2$ is subspectral to $M_1$, then for all $t$, we have $Z_1(t) \leq Z_2(t)$. Cancelling constant terms, we have
\[ \sum_{n,k=1}^\infty \frac{l_k}{\sinh(nl_k/2)}e^{-(nl_k)^2/4t} \leq \sum_{n,k=1}^\infty \frac{m_k}{\sinh(nm_k/2)}e^{-(nm_k)^2/4t}. \]
Order the sets $\{ nl_k\ |\ n,k\in\mb{N}\}$ and $\{nm_k\ |\ n,k\in\mb{N}\}$. (This is possible because for each $L$ there are only finitely many geodesics with length less than $L$.) For each $i$, if $r_i = nl_k$ then set $a_i = l_k/\sinh(nl_k/2)$ and if $s_i = nm_k$ then set $b_i = m_k/\sinh(nm_k/2)$. Substitute $\tau = 1/t$ and apply Lemma \ref{heat_ineq_lemma} to the inequality to conclude that $l_1^2/4 \geq m_1^2/4$ so that we have $l_1\geq m_1$ as claimed.
\end{proof}

\subsection{Length subspectrality implies principal eigenvalue\\ inequality}

Additionally, we can use an inequality of the length spectrum to show an inequality in the lowest nonzero eigenvalue. (Thanks to Dylan Thurston for discussion and some ideas in the proof.) If $M_1$ and $M_2$ are manifolds such that $l_k(M_1) \geq l_k(M_2)$ for all $k$, then say $M_2$ is length-subspectral to $M_1$.
\begin{prop}Let $M_1$ and $M_2$ be two hyperbolic surfaces of the same genus. If $M_2$ is length-subspectral to $M_1$, then the first nonzero eigenvalue of $M_2$ is greater than the first nonzero eigenvalue of $M_1$.\end{prop}

\begin{proof}
The proof proceeds similarly to the previous two: we exploit the Jacobi identity and then apply Lemma \ref{heat_ineq_lemma}.

In order to use the Jacobi identity, we need to show that length subspectrality implies inequality of length trace. Set $F_n(x) = \frac{x}{\sinh(nx/2)}.$ With this, the trace formula for $M$ becomes
\[ Z_M(t) = C_g + \sum_{n,k} F_n(l_k)e^{-(nl_k)^2/4t}. \]
As the exponential is decreasing in $l_k$, to show that the right hand side of the trace formula is decreasing in $l_k$, it suffices to show that $F_n(x)$ is decreasing in $x$ for $x>0$. To see this, differentiate: $F_n'(x) = \sinh^{-2}(nx/2)(\sinh(nx/2) - (nx/2)\cosh(nx/2))$. Since $\sinh(x) > 0$ for $x>0$ to establish the sign of $F_n'$ we examine the second factor. Substituting $t = nx/2$ gives
\begin{align*}
\sinh(t) - t\cosh(t) &= \sum_{k\mbox{ odd}} \frac{t^k}{k!} - \sum_{k\mbox{ even}} \frac{t^{k+1}}{k!}\\
	&= \sum_{k = 3, 5, 7,\ldots} \bigg( \frac{1}{k!} - \frac{1}{(k-1)!}\bigg)t^k \\
	&\leq 0
\end{align*}
Since the series converge absolutely, indeed $F_n$ is decreasing for all $n$.

Hence if $M_2$ is length-subspectral to $M_1$, for each $k$ and $n$ we have $F_n(l_k) < F_n(m_k)$, the length trace series of $M_2$ is greater than the length trace series of $M_1$, and so $Z_2(t) \geq Z_1(t)$. The result follows from Lemma \ref{heat_ineq_lemma}.
\end{proof}

It is not clear the hypotheses of the propositions of this section are ever satisfied. It would be interesting to know whether there exists a pair of hyperbolic surfaces such that one is subspectral to the other, or whether there exists a pair of hyperbolic surfaces such that one is length-subspectral to the other.

\chapter{Subspectrality and Polya's eigenvalue theorem}%
Recall Definition \ref{SubspFunc}. We generalize a result of Polya comparing the eigenvalue counting function to the one-term Weyl function. Recall from Definition \ref{weylfunc} that we denote the one-term Weyl function by $w^\Omega$, or $w$ when context is clear.

\section{Background}

If $\Omega$ is a subset of $\mb{R}^n$, we say that $\Omega$ tiles $\mb{R}^n$ if $\mb{R}^n$ can be expressed as the union of domains congruent to $\Omega$ such that no two of the domains share interior points. Polya \cite{Polya1961} (2.1, 2.2) proved:
\begin{thm}[Polya's theorem]\label{Polya}Let $\Omega\subset\mb{R}^n$ be a normal domain that tiles $\mb{R}^n$. Then $\Omega$ is Dirichlet-superspectral to $w$ and Neumann-subspectral to $w$.\end{thm}

Polya conjectured that this holds for any planar domain. This conjecture is still open.

\subsection{Packings}
A packing is a generalization of the notion of tiling measuring of how badly the domain fills $\mb{R}^n$. Using packing allows us to generalize Theorem \ref{Polya}. We therefore make the following definitions. For a more thorough summary of basic concepts in packing, we refer the reader to Fejes-Toth \cite{Toth1993}, section 2, and to Groemer \cite{Groemer1986}. 

\begin{defn}[Packing constant]
A packing $\mc{P}$ of $\mb{R}^n$ by $\Omega$ is a collection of pairwise disjoint congruent copies of $\Omega$. If $G$ is a domain, define the inner density $d_{\on{inn}}$ and outer density $d_{\on{out}}$ of $\Omega$ with respect to $G$ to be
\[ d_{\on{inn}}(\mc{P}|G) = \frac{1}{|G|}\sum_{A\in\mc{P},A\subset G}|A| \]
\[ d_{\on{out}}(\mc{P}|G) = \frac{1}{|G|}\sum_{A\in\mc{P},A\subset G}|A| \]
We call $(G, o)$ a gauge.

We define the inner (resp. outer) densities of $\mc{P}$ with respect to the gauge $(G,o)$ as
\[ d_-(\mc{P},G,o) = \liminf_{\lambda\to\infty} d_{\on{inn}}(\mc{P}|\lambda G) \]
\[ d_+(\mc{P},G,o) = \limsup_{\lambda\to\infty} d_{\on{out}}(\mc{P}|\lambda G) \]
where $\lambda G$ is the image of $G$ under a homothety of scale $\lambda$ fixing $o$.

The packing density $\delta$ of $\Omega$ is defined to be the supremum of the outer densities $d_+(\mc{P},G,o)$ over all packings of $\Omega$ and all choices of $(G,o)$.
\end{defn}

\subsection{Results related to Polya's theorem}

Urakawa \cite{Urakawa1984} proved a result related to Polya's theorem, bounding Dirichlet eigenvalues in terms of the \emph{lattice packing constant} of a domain. The lattice packing constant is defined by restricting packings to only those which are obtained by translating a domain by a discrete subgroup of isometries.

Li and Yau \cite{Li1983} proved that
\begin{thm}[Li-Yau]\label{LiYau}
Suppose $\Omega$ is a domain in $\mb{R}^n$. Then $\Omega$ is Dirichlet-superspectral to $\big(\frac{n+2}{n}\big)^{n/2}w$.
\end{thm}

They stated this as an inequality of eigenvalues, which we have translated into an equivalent statement about subspectrality using Lemma \ref{SubspFunc}. Li-Yau proved this as a corollary of an inequality involving the sequence of partial sums of the eigenvalue sequence, $\sum_{j=1}^k \lambda_j$. Kroger \cite{Kroger1994} proved a corresponding inequality for Neumann eigenvalues. 

Recent improvements to Theorem \eqref{LiYau} have proceeded by studying the partial sums of the eigenvalue sequence in greater detail and extending the inequality to more general settings. For details, we refer the reader to Laptev \cite{Laptev1997}, Melas \cite{Melas2003}, Wei \cite{Wei2010}, Geisinger \cite{Geisinger2011}, Hatzinikitas \cite{Hatzinikitas2013}, Yolcu and Yolcu-Yolcu \cite{Yolcu2013a} \cite{Yolcu2013} \cite{Yolcu2014}, and Kovarik \cite{Kovarik2015}.

\section{Generalization of Polya's theorem}

Following Polya's original argument and generalizing Urakawa's argument, we use packings of a Euclidean domain $\Omega$. By doing so, the packing constant $\delta$ of $\Omega$ enters the inequality. In fact, we are able to replace $\frac{n+2}{n}$ in Li-Yau's estimate with $\delta$. This replacement sacrifices universality, but strengthens the inequality for domains with high packing constant.

We begin with a lemma describing the behavior of packings.

\begin{lemma}\label{packing_limit}
Let $\Omega$ be a bounded domain in $\mathbb{R}^n$. Let $\mc{P}$ be a maximum-density packing of $\Omega$. For $\sigma > 0$, let $\nu(\sigma)$ be the number of components of $\mc{P}$ entirely contained in $[-\sigma/2,\sigma/2]^n$. Then
\[ \lim_{\sigma\to\infty} \frac{\sigma^n}{\nu(\sigma)} = \frac{|\Omega|}{\delta}. \]
\end{lemma}

\begin{proof}
According to a theorem of Groemer \cite{Groemer1986}, c.f. also section 2 of \cite{Toth1993}, for every compact domain $\Omega$ in $\mb{R}^n$, there exists a packing $\mc{P}$ by congruent copies of $\Omega$ such that 
\[ d_+(\mc{P},G,o) = d_-(\mc{P},G,o) = \delta \]
for every gauge $(G,o)$. We call such a packing a maximum-density packing.

Therefore we may choose a suitable gauge pair: $([-1/2,1/2]^n,0)$. Then
\[ \delta = \lim_{\lambda\to\infty} d_{\on{inn}}(\mc{P},\lambda G) = \lim_{\sigma\to\infty} \frac{1}{\sigma^n |\Omega|}\sum_{A\in\mc{P},A\subset\lambda\Omega}|A|. \]

In view of the fact that every $A$ is congruent to $\Omega$, we have
\begin{align*}
\delta &= \lim_{\sigma\to\infty} \sigma^{-n}\sum_{A\in\mc{P},A\subset\lambda\Omega} |\Omega|\\
	&= \lim_{\sigma\to\infty} \frac{|\Omega| \nu(\sigma)}{\sigma^n}.
\end{align*}

Therefore, as desired,
\[ \lim_{\sigma\to\infty} \frac{\sigma^n}{\nu(\sigma)} = \frac{|\Omega|}{\delta}. \]
\end{proof}

We now prove the following theorem:
\begin{thm}[Generalization of Polya's theorem]\label{PolyaGeneral}
Let $\Omega$ be a bounded domain in $\mathbb{R}^n$. Let $|\Omega|$ be the volume of $\Omega$ and $\delta$ be its packing constant. Then $\Omega$ is Dirichlet-superspectral to $w/\delta$.
\end{thm}

The proof proceeds in two steps. First, we apply Dirichlet domain monotonicity and use Polya's theorem. Second, we equate the limit $\lim_{\eps\to 0}\nu(\eps)\eps^n$ with the packing constant of $\Omega$, according to Lemma \ref{packing_limit} above.

\begin{proof}Let $\sigma > 0$ be given. Set $G = [-\sigma/2,\sigma/2]^n$. Note that $|G| = \sigma^n$. Denote by $N_G$ and $N_\Omega$ the eigenvalue counting functions of $G$ and $\Omega$, resp.

Let $\mc{P}$ be a maximal packing of $\Omega$. Let $\nu(\sigma)$ denote the number of components of $\mc{P}$ contained within $G$. Let $x$ be an arbitrary positive real number. By Dirichlet domain monotonicity,
\[ \sum_{A\in\mc{P}, A\subset G} N_A(x) \leq N_G(x). \]
Since $G$ tiles $\mb{R}^n$, by Polya's theorem, Theorem \ref{Polya}, we have 
\[ N_\Omega(x) \leq \frac{N_G(x)}{\nu(\sigma)} \leq \frac{\omega_n}{(2\pi)^n} \frac{\sigma^n}{\nu(\sigma)} x^{n/2}. \]
This is true for every $\sigma > 0$. Letting $\sigma\to\infty$ and using Lemma \ref{packing_limit}, 
\[ N_\Omega(x) \leq \frac{\omega_n}{(2\pi)^n} \frac{|\Omega|}{\delta} x^{n/2} = \delta^{-1}w(x). \]

As $x$ was chosen arbitrarily, this inequality holds for all $x\in[0,\infty)$. This completes the proof. 
\end{proof}

Equivalently, by Lemma \ref{subsp_func_lemma},
\[ \lambda_k(\Omega) \geq \frac{4\pi^2}{\omega_n^{2/n}} \bigg( \delta \frac{k}{V} \bigg)^{2/n} \]
for all $k$.

Theorem \ref{PolyaGeneral} is a generalization of Polya's theorem, as it makes no assumptions about lattice tiling and the packing constant of a tiling domain is equal $1$. The theorem also permits us to replace the factor $[(n+2)/n]^{n/2}$ in the Li-Yau inequality \ref{LiYau} with the factor $\delta^{-1}$. In particular, if $\delta > [n/(n+2)]^{n/2}$, then the inequality in Theorem \ref{PolyaGeneral} is stronger than the inequality in Theorem \ref{LiYau}. Such domains are not difficult to construct; for instance, Theorem \ref{PolyaGeneral} is stronger than the inequality in Li and Yau's Theorem \ref{LiYau} for any domain in dimension $n\geq 2$ which has a bounding parallelopiped with less than twice the volume of the domain.

General lower bounds on packing constants for various classes of domain in all dimensions tend to be weak. For example, a theorem of Minkowski-Hlawka guarantees that the packing constant for a convex, centrally symmetric domain in $\mb{R}^n$ is no less than $\zeta(n)/2^{n-1}$ where $\zeta$ is the Riemann zeta function. Schmidt proved that there is a constant $c$ such that every convex domain in $\mb{R}^n$ has $\delta \geq cn^{3/2}/4^n$. Compare the discussion in Toth \cite{Toth1993}.

Dimensions $2$ and $3$ are better-studied. We reproduce a portion of a table from a survey by Bezdek \cite{Bezdek2013}, modifying the last row with information from section 8.4 of the same paper:
\bigbreak
\begin{tabular}{|c|c|}
\hline
Body & Lower bound for packing density \\ 
\hline
Unit ball & $\frac{\pi}{\sqrt{18}} = 0.7408\ldots$ \\
Regular octahedron & $\frac{18}{19} = 0.9473\ldots$ \\
Cylinder over a plane domain $K$ & $\delta(K)$ \\
Doubled cone & $\pi\sqrt{6}/9 = 0.855\ldots$ \\
Tetrahedron & $ 0.856\ldots $ \\
\hline
\end{tabular}
\bigbreak
Here a cylinder over a plane domain $K$ is the Minkowski sum of $K\times\{0\}$ with a line segment $s$ (which is assumed non-parallel to $K$). Observe that the packing constant of cylinders implies the three-dimensional case of Laptev's  \cite{Laptev1997} proof of Polya's conjecture for products of tiling domains and arbitrary domains.

In \cite{Torquato2009}, the authors provide a survey of known lattice packing constants for Platonic and Archimedean solids. All the Platonic and Archimedean solids have packing densities in excess of $0.5$.

Specializing to $n=2$, Kuperberg-Kuperberg \cite{Kuperberg1990}, later improved by Doheny \cite{Doheny1995}, found lower bounds for packing constants of convex planar domains: 
\begin{thm}[Kuperberg-Kuperberg, Doheny]
If $\Omega$ is a convex planar domain, then its packing constant is at least $\sqrt{3}/2$.
\end{thm}
This gives the following corollary to Theorem \ref{PolyaGeneral} improving Li-Yau's inequality inequality Theorem \ref{LiYau} by a factor of $\sqrt{3}$.
\begin{cor}\label{PolyaConvexDim2}
Let $\Omega$ be a convex planar domain. Then $\Omega$ is Dirichlet-superspectral to $x\mapsto\frac{|\Omega|x}{2\pi\sqrt{3}}$ and, for all $k$, by Lemma \ref{subsp_func_lemma} its Dirichlet eigenvalues satisfy
\[ \lambda_k > 2\sqrt{3}\pi\frac{k}{|\Omega|}. \]
\end{cor} 
After the author posted this result on the arxiv \cite{Coleman2015}, Iosif Polterovich informed the author that this result was known to him and Olivier Mercier.

Note that in general, packing constants are greater than lattice packing constants. For instance, the regular tetrahedron has a low lattice packing constant and admits non-lattice packings with higher density; see section 8 of Bezdek \cite{Bezdek2013} for more information. In fact, many domains have high packing constants but low lattice packing constants. This is the case even in $\mathbb{R}^2$; for instance, there are triangles tile the plane, but have lattice packing constants are strictly less than one. Replacing lattice packing constant with general packing constant corrects this oversight in \cite{Urakawa1984}.

\section{Polya's inequality for spheres}\label{polya_sphere}

We recall the following facts about round spheres $\mb{S}^n = \{x\in\mb{R}^{n+1}\ |\ |x|=1\}$ and their Laplace spectra. The Riemannian metric on $\mb{S}^n$ is induced by the Euclidean metric on $\mb{R}^{n+1}$. The eigenfunctions of the Laplacian on $\mb{S}^N$ are linear combinations of restrictions to the sphere of homogeneous harmonic polynomials on $\mb{R}^{n+1}$. The eigenvalues of the Laplacian are $\{ k(k+n-1)\ |\ k = 0,1,2,\ldots\}$ and the multiplicity of the $k^{th}$ eigenvalue is equal to $\binom{n+k}{k} - \binom{n+k-2}{k-2}$ for $k\geq 2$, equal to $n+1$ for $k=1$, and equal to $1$ for $k=0$.

Recall that the value of the counting function at a point $x$ is equal to the sum of the multiplicities of the eigenvalues less than or equal to $x$. The sum telescopes. If we establish the convention that $\binom{m}{j} = 0$ for $j<0$ and denote the counting function for $\mb{S}^n$ by $N_n$, we have:
\[ N_n(k(k+n-1)) = \binom{n+k}{k} + \binom{n+k-1}{k-1} \]

\begin{lemma}[Sphere Weyl function]
The one-term Weyl function of $\mb{S}^n$ is
\[ w^n(x) = \frac{2}{n!}x^\frac{n}{2} \]
\end{lemma}
\begin{proof}
Denote by $\omega_n$ the volume of the unit disk $\mb{D}^n = \{x\in\mb{R}^n\ |\ |x|\leq 1\}$. Denote by $\mf{s}_n$ the Riemannian volume of $\mb{S}^n$. Note that $\omega_0 = 1$ and $\mf{s}_0 = 2$. The following recursive relation is established by integrating in spherical and toroidal coordinates, respectively:
\[ \begin{cases}  
\omega_n &= \frac{1}{n}\omega_{n-1} \\
\mf{s}_n &= 2\pi\omega_{n-1} \\
\end{cases} \]
The Weyl function of $\mb{S}^n$ is $w^n(x) = \frac{\omega_n\mf{s}_n}{(2\pi)^n}x^{n/2}$. By the recursive relation we observe that $\omega_n\mf{s}_n = \frac{2\pi}{n}\omega_{n-1}\mf{s}_{n-1}$. Inductively we have 
\[\omega_n\mf{s}_n = \frac{(2\pi)^n}{n!}\omega_0\mf{s}_0 = 2\frac{(2\pi)^n}{n!}. \]
Substituting into the numerator of the coefficient in the Weyl function yields the claim.
\end{proof}

\begin{lemma}\label{sphere_counting}
For all $n$, and all $k>0$, the counting function and Weyl function of the sphere $\mb{S}^n$ satisfy
\[ w^n\big(k(k+n-1)\big) < N_n\big(k(k+n-1)\big). \]
\end{lemma}
\begin{proof}
This assertion is equivalent to:
\[\frac{2}{n!}\big(k(k+n-1)\big)^{n/2} < \binom{k+n}{k} + \binom{k+n-1}{k-1}.\]
We verify by direct computation.

Expand the right side:
\begin{align*}
\binom{k+n}{k} + \binom{k+n-1}{k-1} &= \frac{1}{n!}(k+n)(k+n-1)\cdots(k+1)\\ &\ \ \ \ + \frac{1}{n!}(k+n-1)(k+n-2)\cdots k \\
	&= \frac{1}{n!}(k+n-1)\cdots (k+1)(k+n+k) \\
	&= \frac{2}{n!}\bigg(k+\frac{n}{2}\bigg)(k+n-1)\cdots(k+1).
\end{align*}
Canceling the factor of $2/n!$ and squaring both sides, the assertion holds if and only if we have
\[ \bigg[k(k+n-1)\bigg]^n < \bigg(k+\frac{n}{2}\bigg)^2\big(k+n-1\big)^2\cdots (k+1)^2 \]
which is true if and only if
\[ \frac{k(k+n-1)}{(k+n/2)^2}\prod_{j=1}^{n-1} \frac{k(k+n-1)}{(k+j)(k+n-j)} < 1. \]
First notice $k(k+n-1)/(k+n/2)^2 < 1$, as
\[ k^2 + (n-1)k < k^2 + nk + n^2/4. \]
Now, for each $j=1,\ldots,n-1$ we have
\[ k^2 + (n-1)k < k^2 + nk + j(n-j) \]
and so each factor in the product is less than one. Therefore the entire product is less than one. This establishes the inequality.
\end{proof}

\begin{thm}\label{polya_sphere_thm}
For all $n\geq 1$, the sphere $\mb{S}^n$ is not superspectral to its Weyl function.
\end{thm}
\begin{proof}
By Lemma \ref{sphere_counting}, for each distinct positive eigenvalue $\lambda$ of $\mb{S}^n$, we have $w^n(\lambda) <  N_n(\lambda)$. By Lemma \ref{subsp_func_lemma}, if $s = N( k(k+n-1))$ then $k(k+n-1) = \lambda_s < (w^n)^{-1}(s)$ and so $\mb{S}^n$ cannot be superspectral to $w^n$.
\end{proof}

\begin{figure}[h]
	\includegraphics[scale=0.75]{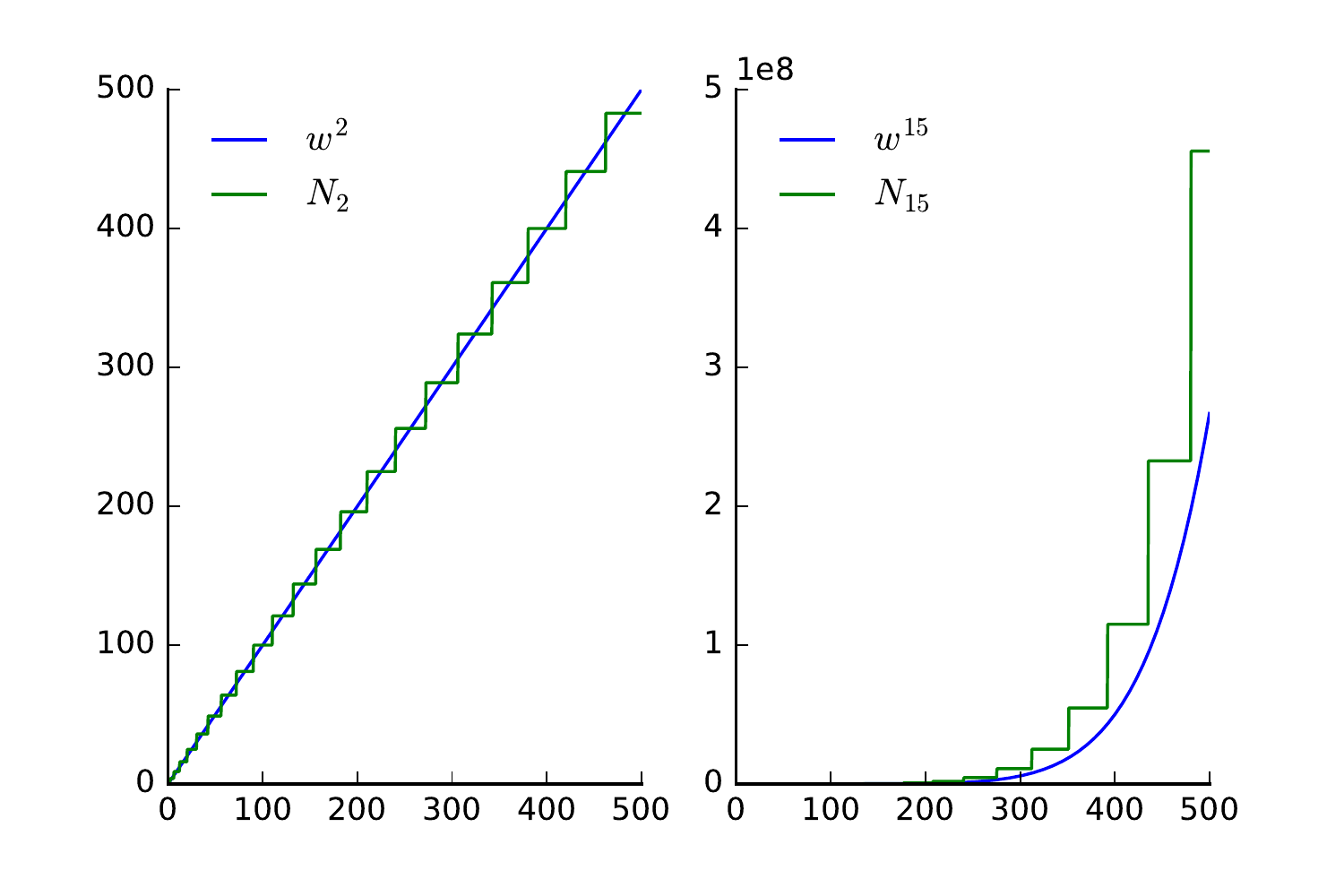}
	\caption[Weyl function and counting functions for $\mb{S}^2$ and $\mb{S}^{15}$]{The left plot graphs the counting function and Weyl function for $\mb{S}^2$. The right plot graphs the counting function and Weyl function for $\mb{S}^{15}$. (Note the vertical axis of the right plot is in units of $10^8$.)}\label{sphere_fig}
\end{figure}

In fact, for the sphere $\mb{S}^2$, we do not have $w^2 \leq N_2$: for sufficiently small $\eps>0$, we have $w^2(k(k+n-1) - \eps) \geq N_2( (k-1)(k+n-2) ) = N_n( k(k+n-1) - \eps)$. For illustration, see Figure \ref{sphere_fig}.

\section{Polya-type inequality for heat traces}

We apply a version of the domain packing argument used to prove Theorem \ref{PolyaGeneral} to produce a bound on the Dirichlet heat trace of a domain. First we prove an inequality between the heat trace on a subset of a domain and the heat trace of the domain.

\begin{lemma}\label{heatineqlemma}
Let $\Omega$ be a normal manifold and let $\Omega_1,\ldots,\Omega_m$ be pairwise disjoint submanifolds of $\Omega$. Then
\[ \sum_{k=1}^m Z_{\Omega_k} \leq Z_\Omega. \]
\end{lemma}
\begin{proof}
Denote by $\sigma_k$ the $k^{th}$ element of the union of the Dirichlet eigenvalues of the $\Omega_i$, when they have been listed in increasing order. Then for each $k$, we have $\sigma_k \geq \lambda_k$. Therefore,
\[ \sum_{k=1}^\infty e^{-t\sigma_k} \leq \sum_{k=1}^\infty e^{-t\lambda_k}. \]
As the sums are absolutely convergent for all $t>0$, we have
\[ \sum_{k=1}^m Z_{\Omega_k}(t) \leq Z_\Omega(t) \]
as claimed.
\end{proof}

We now prove
\begin{prop}
Let $\Omega$ be a compact domain in $\mb{R}^n$ with Dirichlet heat trace $Z(t)$. Then for all $t>0$,
\[ Z(t)\leq \frac{1}{\delta}\frac{|\Omega|}{(4\pi t)^{n/2}}. \]
where $\delta$ is the packing constant of $\Omega$ in $\mb{R}^n$.
\end{prop}

\begin{proof}
Recall that if we scale $\Omega$ by a factor $s>0$, the eigenvalues of $\Omega$ scale by a factor of $s^{-2}$, so
\[Z_{s\Omega}(t) = \sum_k e^{-\lambda_k(s\Omega)t} = \sum_k e^{-\lambda_k(\Omega)s^{-2}s} = Z_\Omega(s^{-2}t).\]

Now take some packing of $\Omega$ into $\mb{R}^n$ with density $\delta$. Consider the square $G_s = [0,s]^n$, let $G=G_1$. Denote the heat trace of $G_s$ by $Z_s$, and the heat trace of $G$ by $Z_G$. For each $s$, let $\nu_s$ denote the number of copies of $\Omega$ contained within $G_s$. By Lemma \ref{heatineqlemma}, for all $s>0$ we have $\nu_s Z_\Omega \leq Z_s$.

Let $t>0$ be arbitrary. Then we have
\[ Z_\Omega(t) \leq \frac{1}{\nu_s} Z_s(t) = \frac{s^n}{\nu_s} \frac{1}{s^n}Z_G(s^{-2}t). \]

We now let $s\to\infty$. The factor $s^n/\nu_s$ tends to $|\Omega|/\delta$ by Lemma \ref{packing_limit}.

We evaluate the expression $\lim_{s\to\infty}\frac{1}{s^n}Z_G(s^{-2}t)$ by substituting $\tau = s^{-2}t$, so that $s^n = (\tau/t)^{-n/2},$ and taking $\tau\to 0$. This yields 
\[s^{-n}Z_G(s^{-2}t) = \frac{1}{t^{n/2}}\tau^{n/2}Z_G(\tau).\]
By the asymptotic expansion of the heat trace, 
\[\lim_{\tau\to 0}(4\pi \tau)^{n/2}Z_G(\tau) =|G|.\]
Therefore, taking the limit $\tau\to 0$ yields
\[\frac{1}{t^{n/2}} \tau^{n/2}Z_G(\tau) \xrightarrow{\tau\to 0} \frac{|G|}{(4\pi t)^{n/2}}.\]

Recalling that $|G| = 1$ and substituting back, we have 
\begin{align*}
Z_\Omega(t) &\leq \lim_{s\to\infty}\frac{s^n}{\nu_s} \frac{1}{s^n}Z_G(s^{-2}t)\\
	&= \frac{1}{\delta}\frac{|\Omega|}{(4\pi t)^{n/2}}
\end{align*}
as claimed.
\end{proof}

\section{Bounds on Laplace spectra of sequences of domains}

In the previous sections, we demonstrated subspectrality by studying the asymptotic behavior of the spectrum of a sequence of domains. We prove the following result on the eigenvalues of a sequence of Euclidean domains whose boundaries satisfy a certain property.
\begin{prop}
Let $\Omega_k\subset\mb{R}^n$ be a sequence of compact Euclidean domains with Dirichlet counting functions $N_k$. Suppose there exists a sequence $c_k\in(0,1)$ such that $c_k\to 0$, $c_k\diam\Omega_k\to \infty$, and 
\[\frac{|\partial\Omega_k^{c_k}|}{|\Omega_k|}\to 0\]
where we define
\[\partial\Omega_k^{c_k} = \{p\in\mb{R}^n\ |\ d(p,\partial\Omega_k)<c_k\diam\Omega_k\}.\]
Then for every $x>0$ we have
\[ \lim_{k\to \infty} \frac{N_k(x)}{|\Omega_k|} = \frac{\omega_n}{(2\pi)^n}x^{n/2}. \]
\end{prop}
\begin{proof}
The proposition is a consequence of the quantitative Weyl law in Theorem \ref{quantitative_weyl}. For every $k>0$, every $x>0$, and every $\epsilon>0$ we have the following inequalities:
\[|\Omega_k^{-\eps\sqrt{n}}|\frac{\omega_n}{(2\pi)^n}x^{n/2}\bigg(1 - \pi\sqrt{\frac{n}{\epsilon^2 x}}\bigg)^n 
\leq N_k(x) 
\leq |\Omega_k^{\eps\sqrt{n}}|\frac{\omega_n}{(2\pi)^n}x^{n/2}\bigg(1 + \pi\sqrt{\frac{n}{\epsilon^2 x}}\bigg)^n \]

Dividing through by the volume of $\Omega_k$, we have:
\[ \frac{|\Omega_k^{-\eps\sqrt{n}}|}{|\Omega_k|}\frac{\omega_n}{(2\pi)^n}x^{n/2}\bigg(1 - \pi\sqrt{\frac{n}{\epsilon^2 x}}\bigg)^n 
\leq \frac{N_k(x)}{|\Omega_k|}
\leq \frac{|\Omega_k^{\eps\sqrt{n}}|}{|\Omega_k|}\frac{\omega_n}{(2\pi)^n}x^{n/2}\bigg(1 + \pi\sqrt{\frac{n}{\epsilon^2 x}}\bigg)^n \]
We leave $x$ fixed. Recall that in the proof of Theorem \ref{quantitative_weyl}, we tiled $\Omega_k$ with cubes congruent to $[0,\eps]^n$. The diameter of each cube is $\epsilon\sqrt{n}$.

Let $\epsilon_k = \frac{1}{2\sqrt{n}}c_k\diam\Omega_k$. Then the set $(\Omega_k - \partial\Omega_k^{c_k})$ is a subset of the union of the cubes that are entirely contained in $\Omega_k$, and the union $\Omega_k\cup \partial\Omega_k^{c_k}$ contains the union of the cubes which have nonempty intersection with $\Omega_k$.

Therefore $|\Omega_k| - |\partial\Omega_k^{c_k}| \leq |\Omega^{-\eps_k\sqrt{n}}|$ and $|\Omega^{\eps_k\sqrt{n}}| \leq |\Omega_k| + |\partial\Omega_k^{c_k}|$.

Then we have for all $k$:
\begin{align*}
\frac{\omega_n}{(2\pi)^n}x^{n/2}\bigg(1 - \frac{|\partial\Omega_k^{c_k}|}{|\Omega_k|}\bigg)\bigg(1 - \pi\sqrt{\frac{n}{\epsilon^2x}}\bigg)
	&\leq \frac{N_k(x)}{|\Omega_k|}\\
	&\leq \frac{\omega_n}{(2\pi)^n}x^{n/2}\bigg(1 + \frac{|\partial\Omega_k^{c_k}|}{|\Omega_k|}\bigg)\bigg(1 + \pi\sqrt{\frac{n}{\epsilon^2x}}\bigg)
\end{align*}
By the hypothesis that $|\partial\Omega_k^{c_k}|/|\Omega_k|\to 0$, we have
\[ \bigg(1 \pm \frac{|\partial\Omega_k^{c_k}|}{|\Omega_k|}\bigg) \xrightarrow{k\to\infty} 1, \]
Because $c_k\diam\Omega_k\to\infty$ and $\eps_k = \frac{1}{2\sqrt{n}}c_k\diam\Omega_k\to\infty$, we have that
\[ \sqrt{\frac{n}{\eps_k^2x}}\to 0 \]
thus
\[ \bigg(1 \pm \sqrt{\frac{n}{\epsilon^2x}}\bigg) \xrightarrow{k\to\infty} 1 \]
Combining these limits yields
\[ \lim_{k\to\infty} \frac{N_k(x)}{|\Omega_k|} = \frac{\omega_n}{(2\pi)^n}x^\frac{n}{2} \]
as claimed.
\end{proof}
The set $\partial\Omega_k^{c_k}$ is a ``thickened boundary.'' The number $c_k$ measures how thick $\partial\Omega_k^{c_k}$ is relative to the diameter of $\Omega_k$. The idea in the proof is that we can let actual thickness, $c_k\diam\Omega_k$, tend to infinity, while the relative thickness tends to zero.

If the condition does not hold, then there are examples where the proposition need not hold. Consider the family of two-dimensional Euclidean rectangles $R_k = [0,k^2]\times [0,1/k]$. Any sequence of $c_k$ satisfying the hypotheses of the theorem must have $c_k\to 0$ and $c_k\diam R_k \to \infty$. For sufficiently large $k$, we have $R_k\subset \partial R_k^{c_k}$, but no point in $R_k$ is ever more than $1/k$ from $\partial R_k$, and so the ratio $|\partial R_k^{c_k}|/|R_k|$ does not tend to $0$.

\section{Generalization of Polya's theorem to Riemannian manifolds}

We conjecture the following generalization of Polya's theorem to Riemannian manifolds.
\begin{conj}
Suppose that $M$ is a complete, contractible Riemannian manifold and that a group $\Gamma$ acts properly, effectively, and cocompactly by isometries on $M$. Every Dirichlet domain of $\Gamma$ satisfies Polya's conjecture for both Neumann and Dirichlet eigenvalues.
\end{conj}

The motivation for this conjecture is the observation that the proof of Polya's theorem on packing domains relies on the interaction between domain monotonicity and a packing of the ambient space by isometric copies of a domain. Thus, it may be possible to extend the argument to domains in manifolds which pack the manifold. To the author's knowledge, there is not yet a quantitative Weyl law for Riemannian manifolds.


\section{Numerical evidence for Polya conjecture}

We numerically approximate the low-frequency Neumann spectra of random pentagons. We also numerically approximate the low-frequency Neumann spectra of annuli. Polya's conjecture holds in each case.

\begin{figure}[p]
	\centering
	\includegraphics[scale=0.9]{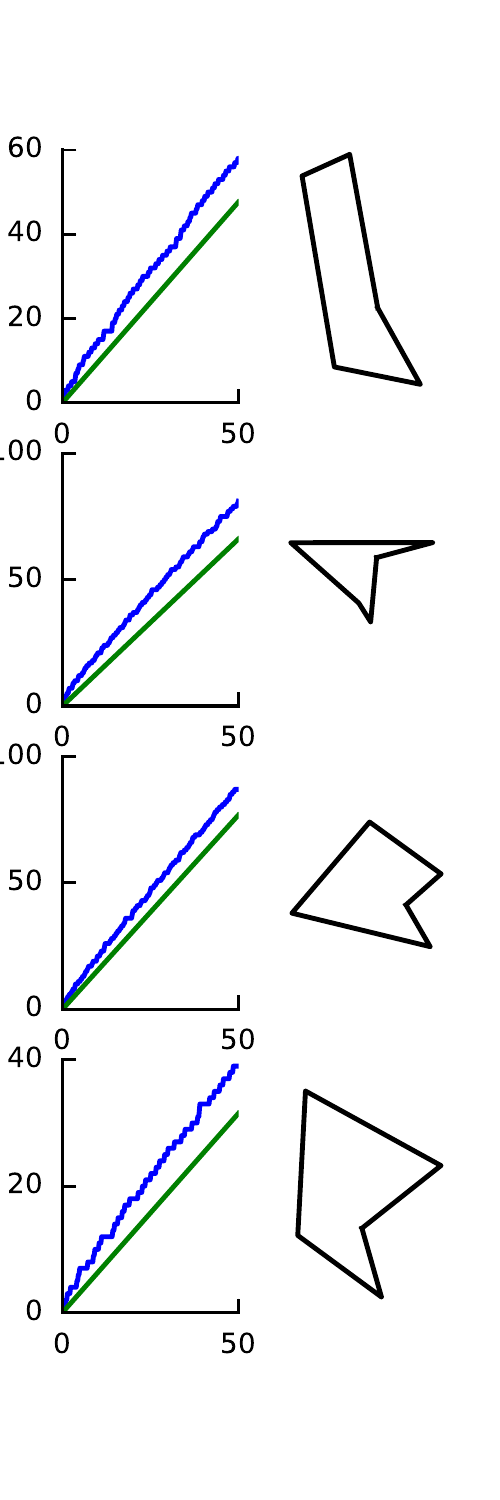}
	\caption[Numerical estimates of low-frequency Neumann spectrum for pentagons]{This figure depicts numerical estimates of the low-frequency Neumann Laplace spectrum of four randomly generated pentagons. The pentagons are depicted in the right column. The Neumann counting function of each pentagon is depicted in its respective plot in blue, while its Weyl function is depicted in green. Polya's conjecture holds in each case for the eigenvalues examined. The pentagons were generated procedurally by taking the origin as one vertex and randomly choosing four points as the other vertices, one in each quadrant, each with norm between $0.5$ and $5$.}\label{random_pentagons}
\end{figure}

\begin{figure}[p]
	\centering
	\includegraphics[scale=0.9]{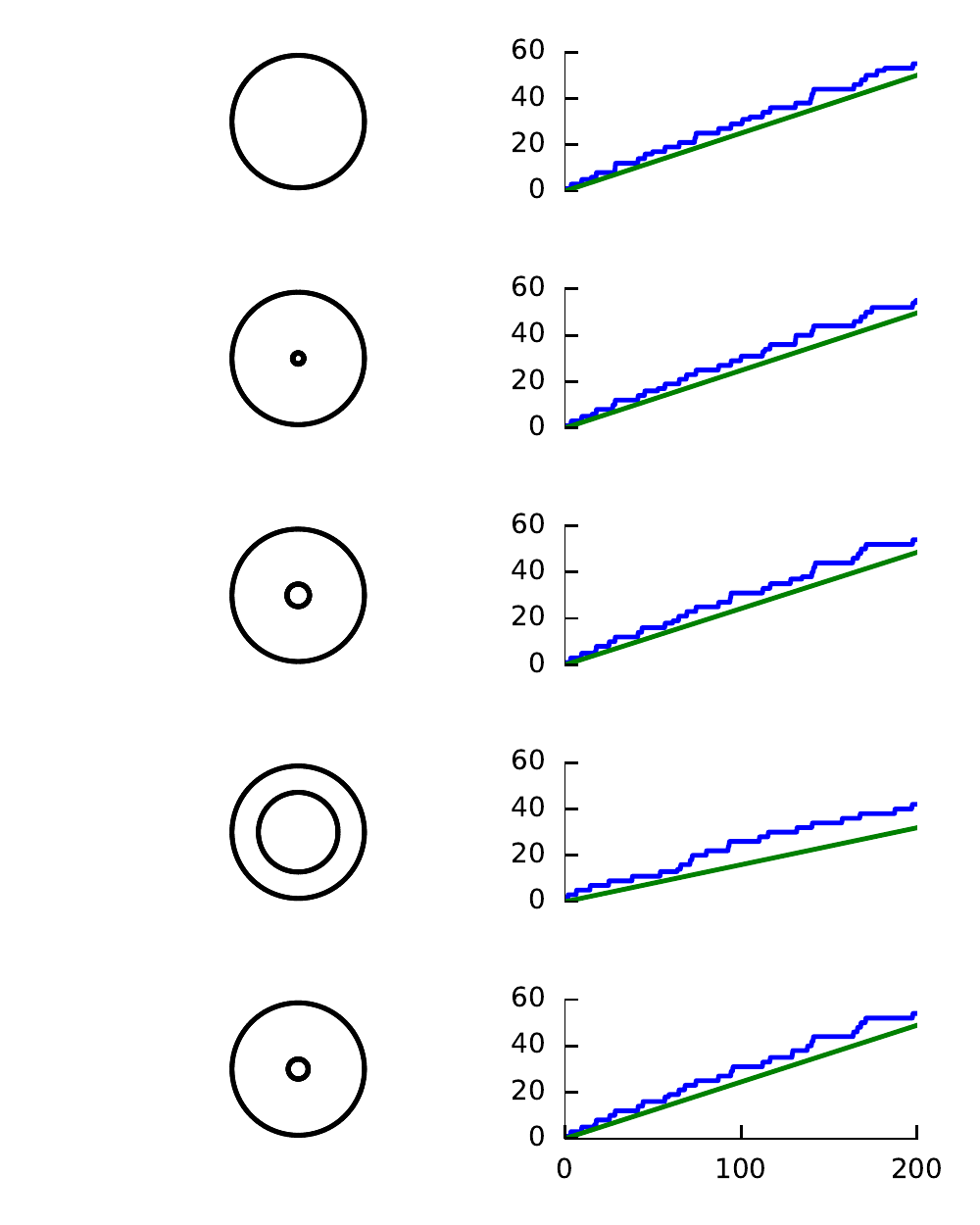}
	\caption[Numerical estimates of low-frequency Neumann spectrum for annuli]{This figure depicts numerical estimates of the low-frequency Neumann Laplace spectrum for the unit disk and four annuli  with major radius equal to $1$ and randomly chosen inner radius less than $1$. Each row holds data for one annulus. The left column contains an image of each annulus, while the right column contains numerically estimated counting Neumann functions and the Weyl function for each annulus. Polya's conjecture holds in each case for the eigenvalues examined.}\label{random_annuli}
\end{figure}

In Figure \ref{random_pentagons} we have a panel of four randomly-generated pentagons and the plot of their corresponding Neumann eigenvalue counting functions for low frequencies. The pentagons were generated by choosing four points $x_i\in\mb{R}^2$, with $x_i$ in the $i^{th}$ quadrant and $0.5\leq \|x_i\|\leq 5$ for each $i$, then defining the pentagon as the shape bounded by the line segments connecting $(1,0)$ to $x_1$, $x_1$ to $x_2$, $x_2$ to $x_3$, $x_3$ to $x_4$, and $x_4$ to $(1,0)$. Meshes with triangle area bounded above by $0.001$ were generated, and Neumann eigenvalues were approximated by order-$1$ finite elements. Each pentagon satisfies Polya's conjecture in the low-frequency spectrum numerically estimated. To illustrate, the Neumann counting function is plotted in blue and the one-term Weyl function is plotted in green; this provides evidence that each pentagon is Neumann-subspectral to the one-term Weyl function over the low frequencies estimated.

%

We also study annuli. In Figure \ref{random_annuli} are four annuli with randomly chosen inradius, and a disk for comparison. All have outradius equal to $1$. The inradius is chosen randomly between $0$ and $0.8$. The annuli are approximated by a mesh with triangle area bounded above by $0.001$ and Neumann eigenfunctions are approximated by order-$1$ finite elements. Each annulus generated satisfies Polya's conjecture. Likewise to the pentagons, the Neumann counting function is plotted in blue and the Weyl function is plotted in green, and we see evidence that for the low frequency eigenvalues estimated, the annuli and disk are Neumann-subspectral to the Weyl function.

Note that, without explicit error bounds on finite elements, these computations do not comprise proof. Of interest would be a study of a posteriori error bounds on the finite element method able to give explicit bounds for individual computations, rather than simply asymptotic bounds on convergence. Such a theory would be analogous to the quantitative Weyl law proven in Chapter 1. In combination with a quantitative Weyl law bounding error away from the two-term Weyl function, Polya's conjecture could be proven for arbitrary planar domains.

\appendix
\chapter{Spectral theorem}%
\section{Friedrichs extension and the spectral theorem}

We include a proof of the Friedrichs extension theorem and the spectral theorem for compactly resolved, self-adjoint operators.

\subsection{Friedrichs extension}

Suppose $H$ is a Hilbert space with inner product $(\cdot,\cdot)$. Let $T$ be a densely defined unbounded operator mapping from a domain $D$ to $H$ bounded below by some constant $k$.

We construct the Friedrichs extension of $T$.

\begin{thm}[Friedrichs extension]\label{friedrichs}
Suppose $T$ is a symmetric operator defined in a Hilbert space $H$. The following exist:
\begin{itemize}
	\item A Hilbert space $V$ which has a natural bounded embedding $\iota$ into $H$;
	\item A self-adjoint operator $\wt{T}:\iota(V)\to H$ extending $T$.
\end{itemize}
\end{thm}

The proof follows that in Riesz-Nagy \cite{Riesz1955} section 124. 
\begin{proof}
Define $\mf{t}(u,v) = (Tu,v)$. For $u,v\in D$, we have $\mf{t}(u,v) = \ol{\mf{t}(v,u)}.$

Choose $\sigma > \max\{1,-k\}$. Then the shifted operator $T+\sigma$ is positive and so we may define the sesquilinear form $\mf{t}_\sigma:D\times D\to \mb{C}$ by
\[ \mf{t}_\sigma(u,v) = \sigma(u,v) + (Tu,v). \]
This defines an inner product $(u,v)_V = \mf{t}_\sigma(u,v)$ on $D$. Let $V$ be the completion of $D$ with respect to $\mf{t}_\sigma$. Note that if a sequence $u_k$ is $V$-Cauchy, then it is $H$-Cauchy, and so converges in $H$. Thus the inclusion $D\hookrightarrow H$ extends to an injective map $\iota:V\hookrightarrow H$.

We have chosen $\sigma$ so for all $u\in V$ we have $\|u\|_V^2 \geq \|\iota u\|^2$, hence the injection $\iota$ is bounded.

\begin{defn}
We say the space $V$ is the form domain of $T$.
\end{defn}

Define by $\rho_X$ the Riesz operator mapping $X\to X^*$ for $X=V,H$. Define the map
\[ K: H\xrightarrow{\rho_H} H^*\xrightarrow{\iota^*} V^*\xrightarrow{\rho_V} V \]
Note that $K$ has the property that for all $u\in H,v\in V$,
\[ (Ku,v)_V = (u,v). \]
We use $K$ to construct a self-adjoint extension of $T$.

We claim that $K$ has the following properties:
\begin{itemize}
\item $K$ is bounded
\item $K$ is symmetric
\item $K$ is injective
\item $K$ is positive
\item $\on{ran} K$ is dense in $V$
\end{itemize}
To see that $K$ is bounded, we compute:
\[ \|Kh\|_V^2 = (Kh, Kh)_V = (h,Kh) \leq \|h\|\|Kh\| \leq \|h\|\|Kh\|_V \]
so $\|Kh\|_V \leq \|h\|.$

To see $K$ is symmetric on $V$, compute:
\[ (Ku,v)_V = (u,v) = \ol{(v,u)} = \ol{(Kv,u)_V} = (u,Kv)_V \]
To see $K$ is symmetric on $H$, compute:
\[ (Ku,v) = (K^2u,v)_V = (Ku,Kv)_V = (u,Kv). \]

To see $K$ is injective, suppose $u,v\in H$ have $Ku=Kv$. Then for arbitrary $w\in H$ we have
\[ (u,w) = (Ku,w)_V = (Kv,w)_V = (v,w) \]
and as $w$ is arbitrary we must have $u=v$.

To see $K$ is positive, observe
\[ (v,Kv) = (Kv,Kv)_V \geq 0. \]

To see $K$ has dense image, suppose $w\in V$ has $(Kh, w) = 0$ for all $h\in H$. Then 
\[ 0 = (Kh,w)_V = (h,w) \]
and so $w=0$. 

Because $K$ is injective, $K$ has an inverse defined on its range $\mc{D}$, call it $L:\mc{D}\to H$. We claim that $L$ is self-adjoint and extends $T+\sigma$.

To see that $L$ is self-adjoint, define the maps $S$ and $U$ on $H\oplus H$ in the following way: Let $S(v\oplus w) = w\oplus v$, and let $U(v\oplus w) = -w\oplus v$. Note that $SU = -US$. Note also that for a subspace $X\subset H\oplus H$, we have $S(X^\perp) = (SX)^\perp$.
We have $S$ is an inversion operator in the sense that
\[ \on{graph} L = S(\on{graph} K). \]

Recall that the $\on{graph} L^* = (U\on{graph} L)^\perp$. Then
\begin{align*}
\on{graph}L^* &= (U\on{graph}L)^\perp \\
	&= (US\on{graph}K)^\perp \\
	&= (-SU\on{graph}K)^\perp \\
	&= -S(U\on{graph}K)^\perp \\
	&= -\on{graph} L \\
	&= \on{graph} L
\end{align*}
and so we have that $L$ is self-adjoint.

We now show that $\dom L$ is dense in $V$. Suppose $v\in V$ is orthogonal in $V$ to $\dom L$. Then for any $w\in \dom A$, 
\[ 0 = (v, w)_V = (v, Lw) \]
and as $L$ surjects onto $H$, we must have $v=0$.

Now we show that $L$ extends $T+\sigma$. Let $v,w\in D$ be arbitrary. Then note, by the defining property of $K$, that
\[ (v,w)_V = (v, (T+\sigma)w) = (v, K(T+\sigma)w)_V \]
and as $\dom L = \on{ran} K$ is dense in $V$ we have for all $w\in D$ that $w = K(T+\sigma)w$.

Thus $D\subset \on{ran} K$ and for $w\in D$ we have 
\[ Lw = L(K(T+\sigma))w = (T+\sigma)w. \]

We have therefore constructed a self-adjoint extension $L$ of $T+\sigma$, defined in a dense subset of the form domain of $T$.
\end{proof}

\subsection{Spectral theorem for compactly resolved\\self-adjoint operators}

For completeness we include a proof of the spectral theorem for compactly resolved self-adjoint operators. References here include Kato \cite{Kato1976} VI.5, Gilbarg-Trudinger \cite{Gilbarg1998} chapters 5 and 8. For an alternate proof of the spectral theorem using the functional calculus, see Rudin \cite{Rudin1987} 13.21.

We begin our study with the 
\begin{thm}[Fredholm alternative]
Suppose $H$ is a Hilbert space and $K:H\to H$ is a compact mapping. Then the operator $I-K$ either has a nontrivial kernel, or has a bounded inverse. That is, either $x-Tx=0$ has a nontrivial solution in $H$ or for every $y$ there is a unique $x$ satisfying $x-Tx=y$, and $(I-K)^{-1}$ is bounded.
\end{thm}

From this, we have the
\begin{thm}[Spectral theorem for compact operators]
Suppose $H$ is a Hilbert space and $K:H\to H$ is compact. Then the spectrum of $K$ consists of a countable set of eigenvalues, each with finite multipliciy, and no limit points except possibly $0$.
\end{thm}

We will use Schauder's theorem:
\begin{thm}[Adjoint of compact is compact]
Suppose $X,Y$ are Banach spaces and $T:X\to Y$ is compact. Then $T^*:Y^*\to X^*$ is compact.
\end{thm}

The following lemma regarding compactly resolved maps will prove useful.

\begin{lemma}[Compactly resolved operators have discrete spectrum]
Let $T$ be a densely defined linear operator in $H$. For any $\lambda$ in the resolvent set of $T$, define the resolvent
\[ R_\lambda = (T-\lambda)^{-1}. \]
If there exists $\lambda_0$ in the resolvent set of $T$ such that $R_{\lambda_0}$ is compact, then the spectrum of $T$ is discrete and comprises eigenvalues.
\end{lemma}
\begin{proof}
Let $K = R_\lambda: H\to\dom T$ and suppose $K$ is compact. Let $L = K^{-1}$. We claim that the spectrum of $L$ is pure point and discrete. Because $K$ is invertible, $0$ is not one of its eigenvalues. By the spectral theorem, the spectrum of $K$ is discrete away from zero and its nonzero elements are eigenvalues.

As for any eigenvalue $u$ of $K$ has $Lu = \mu^{-1}u$, we have
\[\{\mu^{-1}\ |\ \mu \mbox{ an eigenvalue of $K$} \}\subset \mbox{spectrum of } L.\]

We now show that $L$ has no additional spectrum. 

Claim: 
\[ \{\zeta^{-1}\ |\ \zeta\in\mbox{ resolvent set of $K$}\}\subset \mbox{resolvent set of } L. \]

Suppose $\zeta$ is not in the spectrum of $K$. Then $K-\zeta$ is bounded and invertible, so
\[ K-\zeta = \zeta K(\zeta^{-1} - K^{-1}) = -\zeta K (L-\zeta^{-1}). \]
Solving for $L-\zeta^{-1}$ gives
\[ L-\zeta^{-1} = -\zeta L(K-\zeta).\]
Because $K-\zeta$ is invertible and $L$ is (by definition!) invertible, we may invert both sides which gives
\[ (L-\zeta)^{-1})^{-1} = -\zeta^{-1}K(K-\zeta)^{-1}. \]
As $K$ is compact and $(K-\zeta)^{-1}$ is bounded by assumption, we have that $(L-\zeta^{-1})^{-1}$ is bounded, thus $\zeta^{-1}$ is an element of the resolvent set of $L$.

Thus the spectrum of $T$ comprises eigenvalues, accumulates only at $\infty$, and is given by 
\[ \bigg\{ \mu^{-1} - \sigma \bigg|\ \mu\in\mbox{ spectrum of $K$} \bigg\} \]
\end{proof}

We now state and prove the spectral theorem for self-adjoint, bounded-below, compactly resolved operators.

\begin{thm}[Spectral theorem for self-adjoint, bounded-below, compactly resolved operators]\label{spectral_theorem}
Suppose an operator $T$ is self-adjoint and bounded below with domain $D$ in a Hilbert space $H$. Denote by $\mf{t}$ the quadratic form associated to $T$. If $T$ is compactly resolved, then the spectrum of $T$ is discrete and accumulates only at infinity.
\end{thm}
\begin{proof}
Because $\iota: V\to H$ is compact, by Schauder's theorem we have that $\iota^*:H^*\to V^*$ is compact. Then the map $K$ constructed in the proof of the Friedrichs extension is, as a composition of bounded and compact maps, a compact map.

The Friedrichs extension of $T+\sigma$ is the map $L:\mc{D}\to H$ and it is inverted by $K$. Thus we may write
\[ K = (T+\sigma)^{-1} \]
As the domain of $K$ is all of $H$ and $K$ is injective, we have that $-\sigma$ is an element of the resolvent set of $T$, and $K = R_{-\sigma}$.

Thus by the lemma, the spectrum of $T$ is discrete and comprises eigenvalues.
\end{proof}

\section{The form domain of the Laplacian and Sobolev spaces}

In this section, we prove the following proposition:
\begin{prop}[Spectral theorem for the Laplace operator]\label{spectral_thm_laplacian}
Suppose $M$ is a normal Riemannian manifold. Then:
\begin{itemize}
\item If $M$ is closed, the Friedrichs extension of the Laplacian has discrete spectrum comprising eigenvalues
\item If $M$ has boundary, the Friedrichs extensions of both the Dirichlet and Neumann Laplacians have discrete spectrum comprising eigenvalues
\end{itemize}
\end{prop}

In this subsection we suppose $M$ is a normal Riemannian manifold. This requires defining Sobolev spaces on $M$.

\begin{defn}[Differential operator]
A differential operator in $M$ is a linear map $P:C^\infty(M)\to C^\infty(M)$ such that for all elements $p$ of the interior of $M$ and for all local coordinate systems $\phi: U \to M$, there exist $f^\alpha\in C^\infty(U)$ such that for each $u\in C^\infty(M)$ we have
\[ \phi^*(Pu)(x) = \sum_{|\alpha|\leq k} f^\alpha(x)\partial_\alpha \phi^*u(x) \]
for any $u\in C^\infty(M)$. The set of all $k$-order differential operators is written $\on{Diff}^k(M)$.
\end{defn}
Here we use multi-index notation: a multi-index $\alpha$ of order $j$ on $n$ indices is an ordered $j$-tuple $(\alpha_1,\alpha_2,\ldots,\alpha_j)\in \{1, 2, \ldots, n\}^j$. The order of $\alpha$ is denoted $|\alpha| = j$. We write the operation $\partial_\alpha u = \frac{\partial^j u}{\partial_{\alpha_1}\cdots\partial_{\alpha_j}}$ and the monomial $x^\alpha = x^{\alpha_1}\cdots x^{\alpha_j}$.

\begin{lemma}The space $\on{Diff}^k(M)$ is a finitely-generated $C^\infty(M)$ module.\end{lemma}

\begin{proof}
Let $\mc{U}$ be a finite cover of $M$ by coordinate neighborhoods $U_I$ with coordinate charts $\phi_I$, and let $\psi_I$ be an $L^2$ partition of unity subordinate to $\mc{U}$. (See Definition \ref{l2part1}.) For each $\alpha$ define $P_{\alpha, I}$ by $P_{\alpha, I} u = \psi_I \partial_\alpha (\phi_I^*u)$.

Let $P$ be a differential operator. Let $u\in C^\infty(M)$ be arbitrary. For each $I$, we have the existence of $f^\alpha_I$ so that $P = \sum_{|\alpha \leq k|} f^\alpha_I \partial_\alpha$. Claim: $P = \sum_{\alpha,I} \psi_If^\alpha_I P_{\alpha,I}$. Let $u\in C^\infty(M)$ and $p\in M$ be arbitrarily chosen. Then for each $I$ such that $p\in U_I$ we have $\phi^*(Pu)(\phi^{-1}(p)) = \sum_\alpha f^\alpha_I(p)\partial_\alpha \phi^* u(\phi^{-1}(p))$, so that 
\begin{align*}
\sum_{\alpha,I}\psi_I f^\alpha_I P_{\alpha,I}u &= \sum_I \bigg(\sum_\alpha f^\alpha_I \psi_I\partial_\alpha u\bigg) \\ 
	&= \sum_I\psi_I^2 \bigg( \sum_\alpha f^\alpha_I \partial_\alpha u\bigg) \\
	&= \sum_I\psi_I^2 Pu\\
	&= Pu
\end{align*}
As $u$ is arbitrary we have established this identity for $P$ and as $P$ was arbitrary we have shown that each element of $\on{Diff}^k(M)$ can be written as a finite $C^\infty(M)$-linear combination of differential operators.
\end{proof}

\begin{defn}[Distribution; distributional derivative]
A distribution on $M$ is a linear functional mapping $C^\infty(M)\to\mb{C}$. If $P\in\on{Diff}^k(M)$ then for any distribution $T$ we define $PT$ by $(PT)(u) = T(P^*u)$ where $P^*$ is the $L^2$-adjoint of $P$.
\end{defn}


\begin{lemma}[Identifying distributions with elements of $L^2$]
As $C^\infty_0(M)$ is dense in $L^2(M)$, all distributions are linear functionals densely defined in $L^2(M)$. If a distribution is bounded on $C^\infty_0(M)$, then by the Riesz representation theorem it can be given by the action of an element of $L^2(M)$ on itself.
\end{lemma}

\begin{defn}[Sobolev spaces]
We define the spaces
\[ H^k(M) = \{ u\in L^2(M)\ |\ Pu\in L^2(M)\mbox{ for all } P\in\on{Diff}^k(M)\} \]
and
\[ H^k_0(M) = \mbox{ closure of } C^\infty_0(M) \mbox{ in } H^k(M)  \]
\end{defn}
These are Definitions 4.7.1 and 4.7.2 in Taylor \cite{Taylor2011a}.

\begin{lemma}
The space $H^k(M)$ is a Hilbert space when given the $k$-norm
\[ \|u\|_k^2 =  \sum_i \|P_iu\|_{L^2}^2 \]
where $P_i$ is a finite generating set of $\on{Diff}^k(M)$. The topology of $H^k(M)$ does not depend on the choice of generating set.
\end{lemma}
\begin{proof}
The norm can be obtained from an inner product: for any $u,v\in H^k$
\[ \|u\|_k^2 + \|v\|_k^2 = \sum_i \bigg( \|P_iu\|^2 + \|P_iv\|^2\bigg) \]
and the parallelogram law follows from the parallelogram law for the $L^2$ norm applied to each summand.

We now show that $H^k(M)$ is complete under the $k$-norm. First we note that a sequence that is Cauchy in the $k$-norm is Cauchy in the $L^2$ norm. As multiplication by $1$ is a $k$-order differential operator, we have that there exist $f^i\in C^\infty$ such that $1 = \sum f^iP_i$.
\begin{align*}
\| u_j - u_l \|^2 &= \|1(u_j-u_l)\|^2 \\
	&= \bigg\| \sum f^i P_i(u_j - u_l)\bigg\|^2 \\
	&\leq \sum_i \| f^i P_i(u_j - u_l)\|^2 \\
	&\leq \sum_i \sup_M |f^i|^2 \|P_i(u_j - u_l)\|^2 \\
	&\leq \sup_{M,i} |f^i|^2 \sum_i \|P_i(u_j - u_l)\|^2 \\
	&= \sup_{M,i} |f^i|^2 \|u_j - u_l\|_k^2
\end{align*}
where we have used H\"older's inequality and the triangle inequality.

Now suppose $u_j$ is a $\|\cdot\|_k$-Cauchy sequence. It also converges in $L^2$ to some $u\in L^2$. We now claim that $u\in H^k(M)$.

Let $P\in\on{Diff}^k(M)$ be arbitrary. Consider the sequence $Pu_j$. Note that $P = \sum f^iP_i$. Then by the same computation we have
\[ \|Pu_j - Pu_l\|^2 \leq \sup_{M,i} |f^i|^2 \sum \|P_i (u_j - u_l)\|^2 = \sup_{M,i} |f^i|^2 \|u_j - u_l\|_k^2 \]
and because $u_j$ is Cauchy with respect to the $k$-norm, the sequence $Pu_j$ is Cauchy in $L^2$ hence converges in $L^2$ to some $v$.

We now show that $v = Pu$. Suppose $g\in C^\infty_0(M)$ is arbitrary. Then as the natural embedding of a Hilbert space in its dual is continuous by Cauchy-Schwarz, we have
\[ \int_M v\bar{g}\ dV = \lim_j \int_M (Pu_j)\bar{g}\ dV = \lim_j \int_M u_j \overline{P^*g}\ dV = \int_M u\overline{P^*g}\ dV.\]
This holds for arbitrary smooth $g$, so we have that $v$ is equal to $Pu$, hence $u\in H^k(M)$. This shows that $H^k(M)$ is complete with respect to the $k$-norm.

To see that the topology does not depend on choice of coordinates, choose a different finite spanning set $Q_j$ of $\on{Diff}^k(M)$ and note that each $Q_j = \sum S_{ij}P_i$, while each $P_i = \sum T_{ji}Q_j$, where the $S_{ij}$ and $T_{ij}$ are smooth functions. By the triangle inequality and H\"older's inequality, for arbitrary $u\in H^k(M)$, we have:
\begin{align*}
\sum_i \|P_iu\|^2 &= \sum_i \bigg\|\sum_j T_{ij}Q_j u\bigg\|^2 \\
	&\leq \sum_j \bigg(\sum_i \sup_M |T_{ij}|^2\bigg) \|Q_ju\|^2 \\
	&\leq \sup_{i,j,M}|T_{ij}|^2 \sum_j \|Q_ju\|^2
\end{align*}
and likewise $\sum_j \|Q_ju\|^2 \leq \sup_{i,j,M} |S_{ij}|^2 \sum_i \|P_iu\|^2$. As the $S_{ij}, T_{ij}\in C^\infty(M)$ the suprema are finite quantities and so we have that the norms are equivalent. Thus the topology of $H^k$ does not depend on choice of spanning set.
\end{proof}

\begin{lemma}\label{gradient_spans}
Let $U_I$ be a cover of $M$ by coordinate neighborhoods and let $\psi_I$ be an $L^2$ partition of unity subordinate to $U_I$. Let $P_{I,i} = \psi_I \sum_j g^{ij}\partial_j$ where $\partial_j$ denotes the $j^{th}$ partial derivative in $U_I$. Then the set $\{1\}\cup \{P_{I,i}\}$ forms a basis of $\on{Diff}^1(M)$.
\end{lemma}
\begin{proof}
Let $P\in\on{Diff}^1(M)$. In the coordinate system $U_I$ we have $P = \sum_i h^i_I\partial_i$. We have the following smooth functions:  $h:U_I\to\mb{R}^n$ the coefficients of $P$; the metric $g:U_I\to O(n)$; and the inverse metric $g^{-1}:U_I\to O(n)$. The function $f = gh$ has components $f^{i,I}$ that satisfy $P = \sum_{i,j} f^{i,I}g^{ij}\partial_j$. We have
\begin{align*}
P &= \sum_I \psi_I^2 \sum_j h^j_I \partial_j \\
	&= \sum_I \psi_I^2 \sum_{i,j} f^{i,I}g^{ij}\partial_j \\
	&= \sum_{i,I} \psi_I f^{i,I} \bigg( \psi_I \sum_j g^{ij}\partial_j\bigg) \\
	&= \sum_{i,I} (\psi_I f^{i,I}) P_{i,I}  
\end{align*}
showing that indeed $P$ is a $C^\infty$ linear combination of the $P_{i,I}$.
\end{proof}

\begin{prop}\label{form_domain_cts_embeds}
If $M$ is closed, then the form domain of the Laplacian continuously embeds in $H^1(M)$. If $M$ has boundary, then the form domain of the Dirichlet Laplacian is equal to $H^1_0(M)$ and the form domain of the Neumann Laplacian continuously embeds in $H^1(M)$.
\end{prop}
\begin{proof}
The form domain of the Laplacian is defined as the completion of $C^\infty(M)$ with respect to the norm
\[ \| u \|^2_{V_\nu} = \|u\|^2_{L^2} + \int_M \langle \nabla u, \nabla u\rangle\ dx. \]
By Lemma \ref{gradient_spans} the operators $\{1, P_1,\ldots, P_n\}$ span $\on{Diff}^1(M)$. Hence the topology on $H^1(M)$ is defined by the norm
\[ \|u\|_1^2 = \|u\|^2 + \sum_{i=1}^n \|P_iu\|^2. \]
Recall that $\psi_I$ is an $L^2$ partition of unity. Compute:
\begin{align*}
\int_M |\nabla u|^2\ dV &= \sum_I \int_{U_I} \psi_I^2|\nabla u|^2\ dV\\
	&= \sum_I \int \psi_I \sum_j \big| g^{ij}\partial_j u\big|^2\ dV\\
	&= \sum_{I,i} \int \bigg| \psi_I g^{ij}\partial_j u\bigg|^2\ dV\\
	&= \sum \|P_{I,i} u\|^2.
\end{align*}

If $M$ is closed, or has boundary with Neumann boundary conditions imposed, then because $C^\infty(M)\subset L^2(M)$ and convergence under the norm on the form domain coincides with convergence in the $H^1$ norm, the form domain of the Laplacian (or Neumann Laplacian) continuously embeds into $H^1(M)$.

As the form domain of the Dirichlet Laplacian is equal to the closure of its domain with respect to the $V_0$ norm, the $V_0$ norm coincides with the norm on $H^1$, and the domain of the Dirichlet Laplacian is $C^\infty_0(M)$, we have that the form domain of the Dirichlet Laplacian coincides with $H^1_0(M)$.
\end{proof}

We now prove Proposition \ref{spectral_thm_laplacian}.

\begin{proof}
By the Rellich-Kondrachov theorem, c.f. Taylor \cite{Taylor2011a} Chapter 4 7.13, we have that $H^1$ compactly embeds into $L^2$. By Proposition \ref{form_domain_cts_embeds}, the form domain of the Laplacian compactly embeds into $L^2$. By Theorem \ref{spectral_theorem} we have that the Laplacian is compactly resolved and its spectrum is discrete comprising eigenvalues.
\end{proof}
\chapter{Code}%
In this appendix, we record the code used for numerical investigation and generating figures. The code is written in Python 3.5 and makes use of the following open-source libraries:
\begin{itemize}
\item Scipy
\item Triangle
\item Matplotlib
\item Pandas
\end{itemize}
They can be obtained from the python package repository with the following unix commands:
\begin{minted}[fontsize=\footnotesize]{bash}
pip install numpy
pip install triangle
pip install matplotlib
pip install pandas
\end{minted}

\section{Eigenvalues of rectangles}

The following code computes the Laplace spectrum of a rectangle.

\begin{quotation}
\begin{minted}[fontsize=\footnotesize]{python}
import functools
import numpy as np


def cartesian_product(arrays):
    """Compute the cartesian product of the list of arrays
    
    Parameters
    ----------
    arrays : a list of arrays whose product the method outputs
    """
    broadcastable = np.ix_(*arrays)
    broadcasted = np.broadcast_arrays(*broadcastable)
    rows, cols = functools.reduce(np.multiply,
                                  broadcasted[0].shape),
                                  len(broadcasted)
    out = np.empty(rows * cols, dtype=broadcasted[0].dtype)
    start, end = 0, rows
    for a in broadcasted:
        out[start:end] = a.reshape(-1)
        start, end = end, end + rows
    return out.reshape(cols, rows).T

def is_in_ellipse(L, W, lam, pt):
    """Test whether a point x,y is contained in the ellipse
    centered at 0 with semimajor and semiminor axes L\sqrt{lam}/\pi
    and W\sqrt{lam}/\pi
    
    Parameters
    ----------
    L : length
    W : width
    lam : size of ellipse
    pt : [x, y]
    """
    x,y = pt
    return int((x*np.pi/L)**2 + (y*np.pi/W)**2 < lam)

def count(L, W, lam, dirichlet):
    """Compute the eigenvalue counting function of a rectangle
    
    Parameters
    ----------
    L : length of the rectangle
    W : width of the rectangle
    lam : ceiling of the rectangle
    dirichlet : True - Dirichlet
                False - Neumann
    """
    # if dirichlet, start all counting ranges at 1.
    # otherwise, start at 0.
    start = int(dirichlet)
    max_x = np.floor(L*np.sqrt(lam)/np.pi)
    max_y = np.floor(W*np.sqrt(lam)/np.pi)
    x = np.arange(start, max_x+1, 1)
    y = np.arange(start, max_y+1, 1)
    grid = cartesian_product([x,y])
    def test(pt):
        return is_in_ellipse(L, W, lam, pt)
    if not list(grid):
        return 0
    else:
        return np.sum(np.apply_along_axis(test, 1, grid))
\end{minted}
\end{quotation}

\section{Code to generate Chapter 2 figures}

This is the code to generate Figure \ref{quantitative_weyl_rect}.

\begin{quotation}
\begin{minted}[fontsize=\footnotesize]{python}
x = np.arange(0.01, 10, 0.01)
y_err = (1. + np.pi*np.sqrt(0.1**2 + 0.1**2)/np.sqrt(x))**2
fig, ax = plt.subplots(1, 1)
square_eigs = np.array([count(10, 10, lam, dirichlet) for lam in x])
ax.plot(x, square_eigs, linewidth=2, label="Neumann counting function")
weyl_func = area*x/(4*np.pi)
ax.plot(x, weyl_func, color='k', linewidth=2, label="Weyl polynomial")
ax.fill_between(x, weyl_func , weyl_func*y_err , color='gray', alpha=0.5,
                label="Quantitative Weyl law bounds")
ax.spines["top"].set_visible(False)
ax.spines["right"].set_visible(False)
ax.xaxis.tick_bottom()
ax.yaxis.tick_left()
ax.legend(frameon=False, loc="upper left")
fig.savefig("weyl_square.pdf", extension="pdf")
plt.show()
\end{minted}
\end{quotation}

This is the code to generate Figure \ref{metr_subsp_figure}.

\begin{quotation}
\begin{minted}[fontsize=\footnotesize]{python}
N=100
upper_val = sq_eigvals[N]
x = np.arange(0, upper_val, 0.1)
lower_bd = (4./81.)*x
upper_bd = (9./16.)*x

fig, ax = plt.subplots(1,1)
ax.fill_between(x, lower_bd, upper_bd, color='gray', alpha=0.5,
			    label=r'$\frac{4}{81}x\leq y\leq\frac{9}{16}x$')
ax.scatter(sq_eigvals[:N], re_eigvals[:N],
		   label=r'$L = \{(\lambda_k(R),\lambda_k(S))\}$')
ax.spines["top"].set_visible(False)
ax.spines["right"].set_visible(False)
ax.xaxis.tick_bottom()
ax.yaxis.tick_left()
ax.set_xlim((0, upper_val))
ax.set_ylim((0, 9.*upper_val/16.))
ax.legend(frameon=False, loc="upper left")
fig.savefig("rectangle_compare.pdf", extension="pdf")
plt.show()
\end{minted}
\end{quotation}

This is the code to generate Figures \ref{multiple_rects_dirichlet} and \ref{multiple_rects_neumann}.

\begin{quotation}
\begin{minted}[fontsize=\footnotesize]{python}
rect_dict = {}
x = np.arange(0., 100., 0.01)
square_dirichlet = np.array(sq_dir_2)
square_neumann = np.array(sq_neu_2)
recs = np.arange(10., 100., 10)
pct_dir_subsp = []
pct_neu_subsp = []
for L in recs:
    if L not in rect_dict:
        print(str(L)+" ", end="")
        dir_eigs = np.array([count(L, 100./L, lam, True) for lam in x])
        neu_eigs = np.array([count(L, 100./L, lam, False) for lam in x])
        rect_dict[L] = {"D": dir_eigs, "N": neu_eigs}
    pct_dir_subsp.append(pct_subsp(square_dirichlet, dir_eigs))
    pct_neu_subsp.append(pct_subsp(neu_eigs, square_neumann))

x_fix = np.arange(0., 100., 0.01)
fig1, ax1 = plt.subplots()
fig2, ax2 = plt.subplots()
weyl = 100*x_fix/(4*np.pi)
ax1.title.set_text("Dirichlet counting functions")
ax2.title.set_text("Neumann counting functions")
ax1.plot(x_fix, weyl, label="Weyl polynomial", color="k", lw=1.5)
ax2.plot(x_fix, weyl, label="Weyl polynomial", color="k", lw=1.5)
for j in range(1, 10, 2):
    ax1.plot(x_fix, rect_dict[float(10*j)]["D"],
    		 label=(str(10*j)+' by 1/' + str(10*j)),
             color=(0.1 + 0.08*j, 0.0, 0.9 - 0.08*j, 0.6), lw=1.5)
    ax2.plot(x_fix, rect_dict[float(10*j)]["N"],
    		 label=(str(10*j)+' by 100/' + str(10*j)),
             color=(0.1 + 0.08*j, 0.0, 0.9 - 0.08*j, 0.6), lw=1.5)
for ax in (ax1, ax2):
    ax.spines["top"].set_visible(False)
    ax.spines["right"].set_visible(False)
    ax.spines["left"].set_linewidth(0.5)
    ax.spines["right"].set_linewidth(0.5)
    ax.xaxis.tick_bottom()
    ax.yaxis.tick_left()
    ax.locator_params(axis="x", nbins="3")
    ax.locator_params(axis="y", nbins="3")
    ax.set_aspect(ax.get_xlim()[1]/ax.get_ylim()[1])
    ax.legend(frameon=False, loc="best", fontsize=8)
plt.tight_layout()
fig1.savefig("multiple_rects_dirichlet.pdf", extension="pdf")
fig2.savefig("multiple_rects_neumann.pdf", extension="pdf")
plt.show()
\end{minted}
\end{quotation}

This is the code to generate Figure \ref{rect_square_1000}.

\begin{quotation}
\begin{minted}[fontsize=\footnotesize]{python}
eig_dict = {}
sq_eigvals = sorted([(m*np.pi)**2 + (n*np.pi)**2
                     for m in range(1,1000) for n in range(1,1000)])
re_eigvals = sorted([(m*np.pi/2.)**2 + (n*np.pi/3.)**2
                     for m in range(1,1000) for n in range(1,1000)])
for L in np.arange(0.025, 30., 0.005):
    if L not in eig_dict:
        eig_dict[L] = {}
        W = 100./L
        M = 30*W
        N = 30*L
        eig_dict[L]["D"] = np.array(
                               sorted([(m*np.pi/L)**2
                                       + (n*np.pi*L/100.)**2 
                                       for m in np.arange(1,M)
                                       for n in np.arange(1,N)])[:1000]
                                   )
        eig_dict[L]["N"] = np.array(
                               sorted([(m*np.pi/L)**2
                                       + (n*np.pi*L/100.)**2 
                                       for m in np.arange(M)
                                       for n in np.arange(N)])[:1000]
                                   )
        

def count_smaller(arr1, arr2):
    """Count the number of elements in arr1 which are less than
    the corresponding elements of arr2"""
    try:
        return len(arr1[arr1 < arr2])
    except ValueError:
        return -1

def count_larger(arr1, arr2):
    """Count the number of elements in arr1 which are greater than
    the corresponding elements in arr2"""
    try:
        return len(arr1[arr1 > arr2])
    except ValueError:
        return -1


Ls = [l for l in sorted(eig_dict.keys()) if l > 10.025]
num_dir_subsp = np.array([count_smaller(eig_dict[L]["D"],
                                        eig_dict[10.000000000000002]["D"])
                          for L in Ls])/1000.
num_neu_supsp = np.array([count_larger(eig_dict[L]["N"],
                                       eig_dict[10.000000000000002]["N"])
                          for L in Ls])/1000.
fig, ax = plt.subplots(1, 1)

ax.scatter((np.array(Ls)/10.), (num_neu_supsp),
           label="Proportion Neumann $>$ square", c="g", s=1,
           edgecolor="g", alpha=0.5)
ax.scatter((np.array(Ls)/10.), (num_dir_subsp),
           label="Proportion Dirichlet  $<$ square", c="b", s=1,
           edgecolor="b", alpha=0.5)

ax.spines["top"].set_visible(False)
ax.spines["right"].set_visible(False)
ax.spines["left"].set_linewidth(0.5)
ax.spines["right"].set_linewidth(0.5)
ax.set_xlim((1.0,3.0))
ax.set_ylim((0,0.6))
ax.xaxis.tick_bottom()
ax.yaxis.tick_left()
ax.set_xlabel("Ratio of side length of rectangle "
			  + "to side length of square")
ax.set_ylabel("Proportion of first thousand eigenvalues")
ax.title.set_text("Comparing rectangles to square:\nFirst thousand "
				  + "eigenvalues")
ax.locator_params(axis="x", nbins="3")
ax.locator_params(axis="y", nbins="2")
ax.legend(frameon=False, loc="best", fontsize=8)
fig.set_size_inches(8,8)
fig.savefig("rectangles_first1000.pdf", extension="pdf")
plt.show()
\end{minted}
\end{quotation}

This is the code to generate Figure \ref{subspec_mean}.

\begin{quotation}
\begin{minted}[fontsize=\footnotesize]{python}
import pandas as pd

eig_dict = {}
sides = [10, 15, 18, 21]
for L in sides:
    if L not in eig_dict:
        #print(str(L), end=' ')
        eig_dict[L] = {}
        W = 100./L
        M = 40*W
        N = 40*L
        eig_dict[L]['D'] = np.array(
                               sorted([(m*np.pi/L)**2
                                       + (n*np.pi*L/100.)**2 
                                       for m in np.arange(1,M)
                                       for n in np.arange(1,N)])[:1000])
        eig_dict[L]['N'] = np.array(
                               sorted([(m*np.pi/L)**2
                                       + (n*np.pi*L/100.)**2 
                                       for m in np.arange(M)
                                       for n in np.arange(N)])[:1000])


eigs = pd.DataFrame()
for L in sides:
    eigs[str(L)] = eig_dict[L]["D"]

eigs["10>15"] = (eigs["10"] >= eigs["15"]).apply(lambda x: int(x))
eigs["10>18"] = (eigs["10"] >= eigs["18"]).apply(lambda x: int(x))
eigs["10>21"] = (eigs["10"] >= eigs["21"]).apply(lambda x: int(x))

eigs["10>15 cum prob"] = eigs["10>15"].cumsum()/(eigs.index+1)
eigs["10>18 cum prob"] = eigs["10>18"].cumsum()/(eigs.index+1)
eigs["10>21 cum prob"] = eigs["10>21"].cumsum()/(eigs.index+1)

fig, ax = plt.subplots()

ax.plot(eigs.index[:1000],
        eigs["10>15 cum prob"].as_matrix()[:1000],
        label=r"$R = [0,15]\times [0,\frac{100}{15}]$", lw=1)
ax.plot(eigs.index[:1000],
        eigs["10>18 cum prob"].as_matrix()[:1000],
        label=r"$R = [0,18]\times [0,\frac{100}{18}]$", lw=1)
ax.plot(eigs.index[:1000],
        eigs["10>21 cum prob"].as_matrix()[:1000],
        label=r"$R = [0,21]\times [0,\frac{100}{21}]$", lw=1)

ax.spines["top"].set_visible(False)
ax.spines["right"].set_visible(False)
ax.spines["left"].set_linewidth(0.5)
ax.spines["right"].set_linewidth(0.5)
ax.xaxis.tick_bottom()
ax.yaxis.tick_left()
ax.set_xlabel("n")
ax.set_ylabel("Cumulative proportion")
ax.title.set_text(r"Subspectral mean of $R$ and $[0, 10]^2$ at $n$")
ax.locator_params(axis="x", nbins="3")
ax.locator_params(axis="y", nbins="3")
ax.legend(frameon=False, loc="best", fontsize=10)
fig.savefig("subspec_decay1000.pdf", extension="pdf")

plt.show()
\end{minted}
\end{quotation}

\section{Finite element code used in Chapter 4}
We record here the finite element code used to compute eigenvalues in Chapter 4. The finite element method approximates eigenvalues in a regular domain $\Omega$ by first approximating $\Omega$ with a piecewise linear domain $P$, triangulating $P$, and approximating $H^1(P)$ with a finite-dimensional space spanned by piecewise-polynomial functions supported on the triangles of $P$.

For an introduction to the finite element method and its use in solving partial differential equations and approximating eigenvalues, see books by Fix-Strang \cite{Fix2008} and Brenner-Scott \cite{Brenner2008}, or survey papers by Boffi-Gardini-Gastaldi \cite{Boffi2012}, and Melenk-Babuska \cite{Melenk1996}, and references therein.

The following module uses first-order finite elements to numerically approximate the eigenvalues and eigenfunctions of a bounded piecewise-linear Euclidean domain.

\begin{quotation}
\begin{minted}[fontsize=\footnotesize]{python}
from scipy import linalg as lin
from scipy.sparse.linalg import eigsh
import numpy as np
import itertools
import time


# assemble the finite element matrices from a mesh
# v = list of vertices
# f = list of triangles
# 
# mesh format: 
#       {"vertices": np.array([[v0, v1], [v2, v3],...]),
#        "triangles": np.array([[0, 1, 2], [3, 4, 5], ...]),
#        "segments": np.array([[0, 1], [2, 3]])}
#
# returns the two matrices in the finite-element eigenvalue equation
#
# Lx = uMx
#
# L is the matrix of the weak Laplacian, whose elements are inner products of
# gradients of the finite elements
# M is the matrix of the L^2 inner product, whose elements are the L^2 inner
# products of the finite elements
#
# NB this is defined with *piecewise-linear* elements over *triangular* meshes.
# For something more sophisticated, you're going to have to go to the
# professionals at Deal.II, FEnics, or PyDec.
def assembleMatrices(tri):

    # list the vertices
    v = tri["vertices"]
    # list the triangles
    f = tri["triangles"]

    # number of vertices
    n = len(v)

    # initialize the matrices for the Laplacian and the inner product 
    # dimensions are (no. vertices)x(no. vertices)
    # Laplacian
    L = np.zeros( (n,n) )

    # L^2 inner product
    M = np.zeros( (n,n) )

    # now loop over each triangle in the mesh and add the submatrix
    # corresponding to the face
    for t in f:
        # set up a dict that remembers which vertex is zeroth, first, second
        d = {t[0]:0,t[1]:1,t[2]:2}

        # vertices of the triangle
        v0 = np.array(v[t[0]])
        v1 = np.array(v[t[1]])
        v2 = np.array(v[t[2]])

        # area of t
        area = np.abs( np.cross(v2-v0, v1-v0) )/2.0

        # barycentric embedding matrix for the triangle
        A = np.array([[v0[0], v1[0], v2[0]],
                     [v0[1], v1[1], v2[1]],
                     [1., 1., 1.]])

        # invert to find the coordinates of the gradients of the
        # elements
        B = lin.inv(A)
        #g = B[:,:2]

        # now iterate through the 9 entries of 
        for i,j in itertools.product(t,t):

            # first, the matrix for the Laplacian
            # (it ends up being the area times the inner product
            # of the two gradients)
            L[i,j] += area*B[d[i],:2].dot(B[d[j],:2])

            # now, the matrix for the L^2 inner product

            # treat the diagonal differently ...
            if i == j:
                M[i,j] += 2.*area / (12.)
            # ... than the off-diagonal
            else:
                M[i,j] += 2.*area / (24.)

    # return the finished product
    return L, M

# this finds *all* possible eigenvalues using a *dense* solver
# avoid --- it is very expensive for large meshes and is 
def eigenvalues(L, M):
    ev, ef = lin.eigh(L,M)
    return ev, ef

# use sparse arnoldi solver to find first n eigenvalues of Lx = uMx
def sparseEigs(L, M, n=15):

    # hand it to the black box!
    evals,evecs = eigsh(L, n, M, sigma=0.01, which="LM")

    # return the product
    return evals,evecs

# wrap it all up
# this method takes in a mesh and a number n
# and returns the first n eigenvalues/vectors of that mesh
# as a bonus, it also prints how long it took to stdout
def findEigs(mesh, n, bc='Neumann', verbose=False):
    # start the timer
    start = time.time()

    # build the matrices
    if bc is 'Neumann':
        L,M = assembleMatrices(mesh)
    else:
        print("Sorry! Not valid boundary conditions!")

    # find the eigenvalues and eigenvectors
    evals,evecs = sparseEigs(L,M,n)

    # stop the timer
    finish = time.time()

    # print the results
    if verbose:
        print("             time: " + str(finish - start))

    # return the finished product
    return evals,evecs
\end{minted}
\end{quotation}

\section{Code to generate Chapter 4 figures}

The following code estimates the first fifty Neumann eigenvalues of a random quadrilateral or pentagon.

\begin{quotation}
\begin{minted}[fontsize=\footnotesize]{python}
import numpy as np
import FE
import triangle


def generate_quadrilateral():
    ths = np.random.random(4)*(np.pi/2.)+np.array([0., np.pi/2.,
                                                   np.pi, 3*np.pi/2.])
    rs = np.array([np.random.random(4)*5.0 + 0.5]*2)
    return (rs*np.array([np.cos(ths), np.sin(ths)]))

def generate_pentagon():
    q = generate_quadrilateral()
    return np.concatenate((np.array([[1.], [0.]]), q), axis=1)

def random_mesh(nsides=5, area=0.001):
    if nsides==5:
        sklt = {"vertices": generate_pentagon().T,
                "triangles": np.array([[0,1,2], [0,2,3], [0,3,4]])}
    elif nsides==4:
        sklt = {"vertices": generate_quadrilateral().T,
                "triangles": np.array([[0,1,2], [0,2,3]])}
    m = triangle.triangulate(sklt, "ra"+str(area))
    return sklt, m

def findRandomEigs(num_to_find=50, num_sides=4, A=0.01):
    dat = []
    for j in range(10):
        sklt, mesh = random_mesh(num_sides, area=A)
        eigs = FE.findEigs(mesh, n=num_to_find)[0]
        dat.append((sklt, mesh, area(sklt), eigs))
    return dat
\end{minted}
\end{quotation}

This is the code to generate Figure \ref{random_annuli}. 

\begin{quotation}
\begin{minted}[fontsize=\footnotesize]{python}
N = 5
fig, axarr = plt.subplots(N, 2, sharey='col')
fig.set_size_inches(4, 5)
for k in range(N):
    if k == 0:
        inner_radius = 0.
    else:
        inner_radius = 0.8*np.random.random()
    s = generate(inner_radius=inner_radius)
    mesh = triangle.triangulate(s, "pqa0.001")
    A = np.pi*(1. - inner_radius**2)
    eigs = FE.findEigs(mesh, n=75)[0]
    ax1 = axarr[k, 0]
    
    for seg in s["segments"]:
        pt0 = s["vertices"][seg[0]]
        pt1 = s["vertices"][seg[1]]
        ax1.plot([pt0[0], pt1[0]], [pt0[1], pt1[1]], color="k")

    ax1.set_xlim((-1.05, 1.05))
    ax1.set_ylim((-1.05, 1.05))
    ax1.spines["top"].set_visible(False)
    ax1.spines["right"].set_visible(False)
    ax1.spines["bottom"].set_visible(False)
    ax1.spines["left"].set_visible(False)
    ax1.axis('off')
    ax1.set_aspect('equal')
    
    xs = np.arange(0., 200, 0.1)
    y1s = cum_dist(xs, eigs)
    y2s = A*xs/(4*np.pi)
    ax2 = axarr[k, 1]
    ax2.spines["top"].set_visible(False)
    ax2.spines["right"].set_visible(False)
    ax2.spines["left"].set_linewidth(0.5)
    ax2.spines["bottom"].set_linewidth(0.5)
    if k < N-1:
        ax2.set_xticks([])
        ax2.spines["bottom"].set_visible(False)
    else:
        ax2.xaxis.tick_bottom()
    ax2.yaxis.tick_left()
    ax2.locator_params(axis="x", nbins="3")
    ax2.locator_params(axis="y", nbins="3")
    ax2.plot(xs, y1s)
    ax2.plot(xs, y2s)
fig.tight_layout()
fig.savefig("random_annuli.pdf", extension="pdf")
plt.show()
\end{minted}
\end{quotation}

This is the code to generate Figure \ref{random_pentagons}.

\begin{quotation}
\begin{minted}[fontsize=\footnotesize]{python}
N = 4
fig, axarr = plt.subplots(N, 2)
for i in range(min(len(dat), N)):
    data = dat[i]
    sklt = data[0]
    A = data[2]
    eigs = data[3]
    max_eig = np.max(eigs)
    xs = np.arange(0, 50, 0.1)
    weyl = 4*np.pi*(np.arange(len(eigs))/A)
    ax0 = axarr[i,0]
    ax1 = axarr[i,1]
    ax1.spines["top"].set_visible(False)
    ax1.spines["right"].set_visible(False)
    ax1.spines["bottom"].set_visible(False)
    ax1.spines["left"].set_visible(False)
    ax1.axis('off')
    V = sklt['vertices'].take(range(6), axis=0, mode='wrap')
    ax1.plot(V[:,0], V[:,1], color='k')
    ax1.set_aspect('equal')
    ys = cum_dist(xs, eigs)
    weyl_func = (A/(4.*np.pi))*xs
    ax0.plot(xs, ys)
    ax0.plot(xs, weyl_func)
    ax0.spines["top"].set_visible(False)
    ax0.spines["right"].set_visible(False)
    ax0.xaxis.tick_bottom()
    ax0.yaxis.tick_left()
    ax0.locator_params(axis="x", nbins="2")
    ax0.locator_params(axis="y", nbins="3")
fig.set_size_inches(2, 6)
fig.savefig("random_pentagons.pdf", extension="pdf")
\end{minted}
\end{quotation}


\backmatter
\bibliographystyle{plain}
\bibliography{dissertation}

\end{document}